\RequirePackage{ifpdf}
\newif\iffinalversion
\finalversiontrue
\ifpdf
  \iffinalversion
    \documentclass[final,pdftex,a4paper,notitlepage,british]{article}
  \else
    \documentclass[pdftex,a4paper,notitlepage,british]{article}  
  \fi
\else
  \iffinalversion
    \documentclass[final,a4paper,notitlepage,british]{article}
  \else
    \documentclass[a4paper,notitlepage,british]{article}  
  \fi
\fi

\date{22nd May 2008}
\def\currentversion{6}

\usepackage{bm}
\iffinalversion\else
 \ifpdf
  \usepackage[dvipdf,dvipsnames]{xcolor}
 \else
   \usepackage[dvips,dvipsnames]{xcolor}
 \fi
\fi

\ifpdf
\iffinalversion
\usepackage[usenoshowkeys,usenoborder,usenocolorlink,usenosourcespecials,usehyperref, useenumitem]{load_packages}
\else
\usepackage[useshowkeys,usenoborder,usecolorlink,usenosourcespecials,usehyperref, useenumitem]{load_packages}
\fi
\else
\iffinalversion
\usepackage[usenoshowkeys,usenoborder,usenocolorlink,usenosourcespecials,usehyperref, useenumitem]{load_packages}
\else
\usepackage[useshowkeys,usenoborder,usecolorlink,usenosourcespecials,usehyperref, useenumitem]{load_packages}
\fi
\fi
\iffinalversion\else
\fi

\let\bold\undefined
\input{borceux/Diagram}
\makeatletter
\usepackage{tocloft}
\usepackage{mystyle}
\author{Antongiulio Fornasiero%
\thanks{Universit\`a~di~Pisa,
Dipartimento~di~Matematica~``L.~Tonelli'',
Largo~Bruno~Pontecorvo~5,
56127~Pisa,~Italy}
\thanks{Corresponding author. E-mail: antongiulio.fornasiero@googlemail.com}
\and Marcello Mamino%
\thanks{Scuola~Normale~Superiore,
Piazza~dei~Cavalieri~7,
56126~Pisa,~Italy}
\thanks{m.mamino@sns.it}}
\iffinalversion
\title{Arithmetic of Dedekind cuts of ordered~Abelian~groups}
\else
\title{Arithmetic of Dedekind cuts of ordered~Abelian~groups v.~\currentversion}
\fi

\ifpdf % We are running pdftex using hyperref or nohyperref
  \hypersetup{%
    pdfsubject=%
     {Studies the set of cuts of ordered Abelian groups as an algebraic system},
    pdfkeywords=%
     {Dedekind cut, dom, ordered group},
    pdfauthor=%
     {Antongiuluio Fornasiero and Marcello Mamino}
  }
\fi

\allsectionsfont{\sffamily}
\begin{document}
\maketitle
\setlength{\emergencystretch}{1em}
%if \emergencystretch value is positive, TeX will try a pass that allows \emergencystretch worth of extra stretchability to the spaces in each line.
%default 0
\iffinalversion
\else
%\show\sloppy
%\tolerance 9999\emergencystretch 3em\hfuzz .5\p@ \vfuzz \hfuzz 
 \setlength{\overfullrule}{5.0pt}
 %black box for overfull lines
 \setlength{\hfuzz}{7pt}
 %TeX considers a \hbox overfull if the excess width of the box is larger than \hfuzz
 %default 0.1pt
 %\hbadness=5000
 %%If the \badness of a \hbox is less than or equal to \hbadness, the box is acceptable
 %default:1000
\fi
\vspace{-2ex}
\begin{abstract}
\sf
We study Dedekind cuts on ordered Abelian groups.
We introduce a monoid structure on them, and we characterise, via a suitable
representation theorem, the universal part of the theory of such structures.
\end{abstract}

\indent\textsl{MSC:} 06F05; 06F20\\
\indent\textit{Key words:} Dedekind cut; \Dom; Ordered group

\begin{center}
%\begin{minipage}{\tocwidth}
\begin{small}
\tableofcontents
\end{small}
\end{center}

%\def\width{order\xspace}
%% order of an element: x^ = x - x
\section*{Introduction}
In this paper, we consider Dedekind cuts on linearly ordered Abelian groups.
Given such a group $\vg$, call $\vgh$ the set of cuts on $\vg$.
This set is naturally endowed with an order and a minus.
The interesting fact is that there are two non-equivalent ways of defining the sum of two cuts $\Lambda$ and $\Gamma$, which we call the \emph{left} sum (or simply the sum) $\Lambda + \Gamma$ and the \emph{right} sum $\Lambda \ar \Gamma$, which are also definable in terms of each other and the minus.
The resulting structure on the set $\vgh$ is an ordered monoid, in which the cut $\zO:=\cut{(-\infty,0]}{(0,+\infty)}$ is the neutral element.
However, the cancellation law does not hold; specifically, $\Lambda + (-\Lambda) \neq \zO$ in general.

After some preliminaries on ordered sets in \S\ref{SEC:ORDERS},
in \S\ref{SEC:GROUP-BASIC} we determine some basic properties of the \Langnd-structure $\vgh$, where $\Lang$ is the signature  $\struct(\leq,\zO,+,-)$.
%We then take a more abstract approach, and i
In \S\ref{SEC:DOMS} we take a more abstract approach and introduce the notion
of \textbf{double ordered monoids (\doms)}, which are \Langnd-structures satisfying some basic universal axioms, which are true both for every ordered Abelian group and for the set of Dedekind cuts of such groups.
We then show that many properties follow from these axioms alone.
In particular, the concept of \intro{\width{}} of an element turns out to be of fundamental importance:
%There is a natural action of $\vg$ on $\vgh$.
for every $\Lambda \in \vgh$, the invariance group of $\Lambda$ is defined as
\[
\IG \Lambda := \set{\mu \in \vg: \mu + \Lambda = \Lambda}
\]
(where $+$ represents the natural action of $\vg$ on $\vgh$),
and the \width of $\Lambda$ is equal to $\ig \Lambda := (\IG \Lambda)^ +$
(where, for every $Z \subseteq \vg$, $Z^+ \in \vgh$ denotes the upper edge of~$Z$).
Since $\ig \Lambda = \Lambda \ar (- \Lambda)$,
the \width can be defined in the language \Langnd,
and it measures to what extent the cancellation law fails in~$M$.

In \S\ref{SUBSEC:TYPE} we show that \doms can be classified into three types, according to a simple rule.
For instance, for every ordered Abelian group $\vg$, the group itself is of the first type,
while $\vgh$ will be of either of the second or the third type, according to whether $\vg$ is discrete, or densely ordered.

\begin{comment}
Moreover, in \S\ref{SEC:SIGNATURE} we introduce the useful notion of \intro{signature} of elements.
For instance, every $\Lambda \in \Razhat$ is of the form $\lambda^+$, $\lambda^-$, $\mu$ or $\pm \infty$, where $\lambda \in \Raz$ and $\mu \in \Real \setminus \Raz$.
Then, $\sig{+\infty} = \sig{\lambda^+} = +1$, $\sig{-\infty} = \sig{\lambda^-} = -1$, and $\sig{\mu} = 0$.
We show that the signature can be defined on any \dom, and give some rules that help in the computation of $\sig{\Lambda +  \Gamma}$ from $\sig{\Lambda}$ and $\sig{\Gamma}$.
\end{comment}

For every \dom~$M$, let $\Mgroup := \set{x \in M: \ig x = \zO}$,
the set of those elements with a trivial invariance group.
$\Mgroup$ modulo a suitable equivalence relation is an ordered Abelian group, which we denote with~$\G{M}$,
% In \S\ref{SUBSEC:ASS-GROUP} we define an equivalence relation $\Mequiv$ on $\Mgroup$ in such a way that $\GM := \Mgroup/\Mequiv$ is an ordered Abelian group, the group associated to~$M$.
% In particular, for every ordered Abelian group~$H$, the equivalence relation
% $\Mequiv$ identifies $\lambda^+$ with $\lambda^-$,
and $\G{\Dedekind H}$ is the Cauchy completion of~$H$.
\begin{comment}
If $M$ is of the third type, the elements of $\Mgroup$ of signature $0$ form a subgroup $\HM$ of $\GM$, and we have that, if $H$ is a densely ordered Abelian group, then $\collapse{\Dedekind H} = H$.
\end{comment}

Many variants of the basic constructions on groups
(quotients, direct products, etc.)
turn out to have
useful equivalents in the context of \doms:
we examine some of them in \S\ref{SEC:CONSTRUCTIONS}.

In \S\ref{SEC:EMBEDDING} we prove the main theorem (Thm.~\ref{Munipart}):
the universal part of the theory cuts of ordered Abelian groups is given precisely by the axioms of \doms plus the condition $-\zO < \zO$
(i.e. if a universal sentence for cuts is true, then it can be proven using the axioms for \doms alone, plus $\zO < \zO$);
%which yields the main theorem (Thm.~\ref{Membedding}):
moreover, every \dom
% of the second or third type
satisfying the additional condition $-\zO < \zO$
is a \subdom of $\vgh$ for some ordered Abelian group~$\vg$.
In conclusion, a \dom is nothing else than a substructure of $\vgh$ or $\vgt$ (defined in \ref{EX:G-TILDE}) for some ordered Abelian group~$\vg$,
and the axioms of \doms characterise the class of such substructures.

In \S\ref{SEC:AXIOMS} we study the independence of the axioms for \doms we gave in \S\ref{SEC:DOMS},
and we also give some alternative axiomatisations.

In \S\ref{SEC:VALUATION} we examine the generalisations of classical concepts of valuation theory from ordered Abelian groups to \doms.
We also study phenomena that are peculiar to \doms, namely \emph{strong valuations}, which have only trivial counter-parts on groups.

We will now explain some of the motivation for studying Dedekind cuts on ordered Abelian groups.
On one hand we think that such objects are quite natural, and deserve to be examined for their own sake. Moreover, the theory turns out to be more complex than one could think at first sight, but still manageable.
On the other hand, the knowledge of the arithmetic rules of Dedekind cuts is necessary in the study of many ``practical'' contexts.
For instance, let $\KF$ be a valued field, with value group~$G$:
if one wants to build an additive complement to the valuation ring of~$\KF$, one needs to study the cuts on~$G$.
In this context, the first author needed to prove
statements like
Corollary~\ref{COR:D-K-M-J}
\begin{comment}
\[
\forall j,j',k,m,d\in\Nats i,j,k<m \wedge j + j' = m + d \\
\rightarrow \\
(x - j y) + (x' - j' y) + (m y - k y) \leq (x - x') - (d - k) y.
\]
\end{comment}
when he undertook the study of \doms.

If $\KF$ is an ordered field, the usefulness of the study of cuts on the additive and multiplicative groups of $\KF$ has already been recognised, for instance they are the theme of~\cite{GONSHOR:1985,TRESSL:2005,KUH?22,GULDENBERG:2004}.
%However, in the context examined in the aforementioned articles the groups studied are divisible, and the theory becomes quite different from the general case.
%In the study of the integer part of $\Kf$ it t
%For instance, the first author needed to prove
%Corollary~\ref{COR:D-K-M-J} when he undertook the study of \doms, for the cuts on the value group of a valued field~$\KF$.

The article should be understandable to everybody with some basic knowledge of algebra, except for \S\ref{SEC:EMBEDDING}, where some acquaintance with model theory is required.

We wish to thank A.~Berarducci, S.~and F.-V.~Kuhlmann, and M.~Tressl for many useful discussions on the topic of the article.
\begin{proviso*}
Unless we say otherwise, all orders will be linear, and all groups Abelian.
\end{proviso*}
\section{Dedekind cuts of ordered sets}\label{SEC:ORDERS}
Let $O$ be a (linearly) ordered set.
A subset $S$ of $O$ is \intro{convex} if for every $\lambda$, $\lambda' \in S$ and $\gamma \in O$, if $\lambda \leq \gamma \leq \lambda'$, then $\gamma \in S$.
It is \intro{initial} (resp. \intro{final}) if for every $\lambda \in S$ and $\gamma \in O$, if $\gamma \leq \lambda$ (resp. $\gamma \geq \lambda$), then $\gamma \in S$.

A \intro{cut} $\cut{\Lambda^L}{\Lambda^R}$ of $O$ is a partition of $O$ into two subsets $\Lambda^L$ and $\Lambda^R$, such that, for every $\lambda \in \Lambda^L$ and $\lambda' \in \Lambda^R$, $\lambda < \lambda'$.
We will denote with $\Dedekind{O}$ the set of cuts $\Lambda:=\cut{\Lambda^L}{\Lambda^R}$ on the ordered set $O$ (including $-\infty := \cut{\emptyset}{O}$ and $+\infty := \cut{O}{\emptyset}$).

In this \Wpartuno, unless specified otherwise, small Greek letters $\gamma, \lambda, \dotsc$ will denote elements of $O$, capital Greek letters $\Gamma,\Lambda, \dotsc$ will denote elements of $\Oh$.

Given $\gamma\in O$, 
\[%\begin{aligned}
\gamma^- := \cut{(-\infty, \gamma)}{[\gamma, +\infty)} \quad \text{and} \quad
\gamma^+ := \cut{(-\infty, \gamma]}{(\gamma, +\infty)}
%\end{aligned}
\]
are the cuts determined by $\gamma$.
Note that $\Lambda^R$ has a minimum $\lambda\in O$ iff $\Lambda = \lambda^-$.
Dually, $\Lambda^L$ has a maximum $\lambda\in O$ iff $\Lambda = \lambda^+$.

To define a cut we will often write $\Lambda^L := S$ (resp.\ $\Lambda^R := S'$), meaning 
that $\Lambda$ is defined as $\cut{S}{O \setminus S}$ (resp.\ $\Lambda := \cut{O \setminus S'}{S'}$)
when $S$ is an initial subset of $O$ (resp.\ $S'$ is a final subset of $O$).
For instance, the above definition of $\gm$ and $\gp$ can be written $(\gm)^L := (-\infty
,\gamma)$ and $(\gp)^L := (-\infty, \gamma]$.

The \intro{ordering} on $\Oh$ is given by $\Lambda\leq \Gamma$ if $\Lambda^L\subseteq\Gamma^L$ (or, equivalently, $\Lambda^R\supseteq\Gamma^R$).\\
To simplify the notation, we will write $\gamma < \Lambda$ as a synonym of $\gamma\in\Lambda^L$, or equivalently $\gm < \Lambda$, or equivalently $\gp \leq \Lambda$.
Similarly, $\gamma > \Lambda$ if $\gamma\in\Lambda^R$, or equivalently $\gp > \Lambda$.
Hence, we have $\gm < \gamma < \gp$.

An ordered set $O$ is \intro{complete} if for every $S \subseteq O$, the least upper bound and the greatest lower bound of $S$ exist.
\begin{remark}\label{REM:COMPLETE-SET}
If $O$ is an ordered set, then $\Oh$ is complete.
\end{remark}
%%\begin{proof}
%% Let
%% \[
%% \Lambda^L := \bigcup_{\Gamma \in S}\Gamma^L.
%% \]
%% We claim that $\Lambda$ is the \lub of $S$.
%% In fact, by definition, $\Lambda \geq \Gamma$ for every $\Gamma \in S$.
%% Moreover, if $\Theta$ is an upper bound for $S$, then $\Theta^L$ contains $\Gamma^L$ for every $\Gamma \in S$.
%% Hence $\Theta^L \supseteq \Lambda^L$, namely $\Theta \geq \Lambda$.
%% Similarly, one proves that $S$ has a \glb.
%%\end{proof}
%%
The ordering induces a topology on $O$, where a basis of open sets is the family of open intervals $(a,b)$, as $a < b$
vary in $O$.
By~\cite[Theorem~X.20]{BIRKHOFF:1967}, an order is complete iff it is compact; hence, $\Oh$ is compact.

Given a subset $S\subseteq O$, the \intro{upper edge} of $S$, denoted by $\Sp$, is the smallest cut $\Lambda$ such that $S \subseteq \Lambda^L$.
Similarly, $\Sm$, the \intro{lower edge} of~$S$, is the greatest $\Lambda \in \Oh$ such that $S \subseteq \Lambda^R$.
Note that $\Sp = -\infty$ iff $S$ is empty, and $\Sp = +\infty$ iff $S$ is unbounded.

If $Z \subseteq O$ (resp.\ $Z \subseteq \Oh$) we will denote by $\sup Z$ 
the least upper bound of $Z$ in~$O$ (resp.\ in $\Oh$), provided that it exists.

Note that that $\Sp = \sup\set{\gp:\gamma\in S}$,
%\footnote{The supremum on the \rhs is taken in the ordered set~$\Oh$.}
and $\Sp > \gamma > \Sm$ for every $\gamma\in S$.
Moreover, $\Lambda = (\Lambda^L)^+ = (\Lambda^R)^-$.
Note that, for every $\gamma \in O$, $\gp = \mset{\gamma}^+$, and
$\gm = \mset{\gamma}^-$.

An ordered set $O$ is \intro{densely ordered} if for every $x < x'\in O$ there exists $y \in O$ such that $x < y < x'$.
It is \intro{discrete} if it is discrete as a topological space in its order topology.

Assume that $O$ is a densely ordered set.
Then we can define an equivalence relation $\Mequiv$ on $\Oh$ by $\gm \Mequiv \gp$ for every $\gamma \in O$.
The ordering on $\Oh$ induces an ordering on the set of residues $\Oh/\Mequiv$, which is complete and densely ordered.
Since any complete and densely ordered set is connected, $\Oh/\Mequiv$ is connected.
\section{Dedekind cuts of ordered groups}\label{SEC:GROUP-BASIC}
Let $\vg$ be a (linearly) ordered (Abelian) group.
Note that a non-trivial group $\vg$ is discrete iff there is a minimal positive element. Otherwise, it is densely ordered.

Given $S, S' \subseteq \vg$ and $\gamma \in \vg$, define
\[\begin{aligned}
\gamma + S &:= \set{\gamma + \sigma: \sigma \in S}\\
S + S' &:= \set{\sigma + \sigma': \sigma \in S, \sigma' \in S'}.
\end{aligned}
\]
Given $\Lambda, \Gamma \in \vgh$, their (left) \intro{sum} is the cut
\[
\Lambda + \Gamma := (\Lambda^L + \Gamma^L) ^+;
\]
%$\Lambda + \Gamma$ having as left elements the set
\ie,
$
(\Lambda + \Gamma)^L = \set{\lambda + \gamma : \lambda \in \Lambda^L, \gamma \in \Gamma^L}
$.

$(\vgh,+)$ is an Abelian monoid: \ie, the addition is associative and commutative, and $0^+$ is its neutral element.
However, it does not obey the cancellation law:
i.e.\ there exist $\Lambda$, $\Lambda'$ and $\Gamma \in \vgh$ such that $\Lambda + \Gamma = \Lambda' + \Gamma$, but $\Lambda \neq \Lambda'$;
for instance, take $\Gamma = -\infty$, and $\Lambda$, $\Lambda'$ any cuts.

$(\vgh, \leq)$ is a complete linearly ordered set; moreover, $(\vgh, +, \leq)$ is an ordered monoid, namely if $\alpha \leq \beta$, then $\alpha + \gamma \leq \beta + \gamma$.

We can also define right addition by
\[
\Lambda \ar \Gamma := (\Lambda^R + \Gamma^R)^-;
\]
\ie,
$
(\Lambda \ar \Gamma)^R = \set{\lambda + \gamma: \lambda > \Lambda, \gamma >\Gamma}
$.

\begin{remark}
$\Lambda + \Gamma \leq \Lambda \ar \Gamma$.
\end{remark}
$(\vgh, \ar, \leq)$ is also an ordered Abelian monoid, with neutral element~$0^-$.
The map $\phi^+$ (resp. $\phi^-$) from $\struct(\vg, \leq, 0, +)$ to $\struct(\vgh, \leq, 0^+, +)$ (resp. to $(\vgh, \leq, 0^-, \ar)$) sending $\gamma$ to $\gp$ (resp. to $\gm$) is a homomorphism of ordered monoids.

Given $\gamma \in \vg$, we write
\[
\gamma + \Lambda := 
\cut{\gamma + \Lambda^L}{\gamma + \Lambda^R} = 
\cut{\set{\gamma + \lambda: \lambda\in\Lambda^L}}{\set{\gamma + \lambda': \lambda' \in \Lambda^R}}.
\]
One can verify that $\gamma + \Lambda = \gp + \Lambda = \gm \ar \Lambda$.

Consider the anti-automorphism $-$ of $(\vg, \leq)$ sending $\gamma$ to $-\gamma$.
It induces an isomorphism (with the same name $-$) between $(\vgh, \leq, +)$ and $(\vgh, \geq, \ar)$,
sending $\Lambda$ to  $\cut{-\Lambda^R}{-\Lambda^L}$.
Hence, all theorems about $+$ have a dual statement about $\ar$.
\begin{remark}
$-(\gp) = (-\gamma)^-$ and
$-(\gm) = (-\gamma)^+$ for all $\gamma \in \vg$.
\end{remark}

\begin{definizione}
Given $\Lambda, \Gamma\in\vgh$, define their (right) difference $\Lambda - \Gamma$ in the following way:
\[
\Lambda - \Gamma := (\Lambda^R - \Gamma ^L)^+;
\]
\ie,
$
(\Lambda - \Gamma)^R = \set{\lambda - \gamma : \lambda > \Lambda, \gamma < \Gamma}
$.
\end{definizione}
Note that $\Lambda - \Gamma$ is not equal to $\Lambda + (-\Gamma)$ in general. 
\begin{enumexamples}
\item If $\vg = \Raz$, we have $3$ kinds of cuts in $\vgh$ (besides $\pm\infty$): rational cuts of the form $\gp$ (\eg $0^+$), rational cuts of the form $\gm$ (\eg $0^-$), and irrational cuts (\eg~$\sqrt 2$).
\item If $\vg = \Real$, we have only $2$ kinds of cuts in $\vgh\setminus\set{\pm\infty}$: cuts of the form $\gp$ and cuts of the form $\gm$.
\item If $\vg = \Zed$, all cuts in $\vgh\setminus\set{\pm\infty}$ are of the form $\gp = (\gamma + 1)^-$.
\item An important source of counterexamples is the ordered group $\Zedtwo$ (the localisation of $\Zed$ at the prime ideal $(2)$):
it is the subgroup of $\Raz$ of fractions with odd denominator.
\end{enumexamples}
\begin{lemma}\label{LEM:MINUS}
For all $\Lambda, \Gamma, \Theta \in \vgh$ we have:
\[\begin{aligned}
\Lambda - \Gamma     &= \Lambda \ar (-\Gamma),\\
(\Lambda - \Gamma)^L &= \bigcap_{\gamma < \Gamma}\Lambda^L -\gamma, \\
(\Lambda - \Gamma)^R &= \bigcup_{\gamma < \Gamma}\Lambda^R -\gamma.
\end{aligned}
\]
Moreover, 
\[\begin{array}{l@{\,}l@{\,}l@{\,}l}
\Lambda + \Gamma  &= \set{\lambda + \gamma: \lambda < \Lambda, \gamma < \Gamma}^+ &=
\sup \set{\lambda + \Gamma: \lambda < \Lambda} &=
\sup \set{\Lambda + \gamma: \gamma < \Gamma};\\
\Lambda - \Gamma &= \set{\lambda - \gamma: \lambda > \Lambda, \gamma < \Gamma}^- &=
\inf\set{\lambda - \Gamma : \lambda > \Lambda} &= 
\inf \set{\Lambda - \gamma: \gamma < \Gamma} =\\
&= \set{\alpha: \alpha + \Gamma \leq \Lambda}^-.
\end{array}\]
Finally,
$
\Lambda - (\Gamma + \Theta) = (\Lambda - \Gamma) - \Theta = (\Lambda - \Theta) - \Gamma.
$
\end{lemma}
Note that $\Lambda + \Gamma \leq \inf\set{\Lambda + \Gamma': \Gamma' > \Gamma}$, but equality does not hold in general, even when $\vg$ is densely ordered.
For instance, take $\vg = \Zedtwo$, $\Lambda = -\unmezzo$, $\Gamma = \unmezzo$.

One can also define the left difference $\Lambda \dl \Gamma$ as:
\[
\Lambda \dl \Gamma := (\Lambda^L - \Gamma^R)^+;
\]
\ie,
$
(\Lambda \dl \Gamma)^L = \set{\lambda - \gamma : \lambda < \Lambda, \gamma > \Gamma}
$.
It is easy to see that $-(\Lambda + \Gamma) = (-\Lambda) - \Gamma$, and that
$-\Lambda = 0^- - \Lambda = 0^+ \dl \Lambda$.
\begin{remark}
$\Lambda - \Gamma \geq \Lambda \dl \Gamma$.
\end{remark}
%% \begin{proof}
%% $\Lambda - \Gamma = \Lambda \ar (-\Gamma) \geq \Lambda + (-\Gamma) = \Lambda \dl \Gamma$.
%% \end{proof}
%%
\begin{remark}
If $\alpha + \beta = \gamma$, then
\[\begin{aligned}
\alpha^+ + \beta^+   &= \gp, & 
\alpha^+ \ar \beta^+   &\geq \gp,\\
\alpha^+ + \beta^-   &= \gm, &
\alpha^+ \ar \beta^- &= \gp,\\
\alpha^- + \beta^-   &\leq \gm,&
\alpha^- \ar \beta^- &= \gm,\\
\Lambda + \gp &= \Lambda \ar \gm =\Lambda + \gamma.
\end{aligned}\]
\end{remark}
It can happen that $\alpha^- + \beta^- < \gm$, and similarly $\alpha^+ \ar \beta^+ > \gp$.
For instance, take $\vg = \Zed$, and $\alpha$, $\beta$ any integers.
%%
%\subsection{Basic inequalities}
%%
\begin{lemma}\label{LEM:CUTS-AXIOMC}
$\Lambda < \Gamma$ iff $\Lambda - \Gamma < 0^+$.
\end{lemma}
\begin{proof}
\begin{itemize}
\item[$\Rightarrow)$]
If $\Lambda < \alpha < \Gamma$, then $0 = \alpha - \alpha > \Lambda - \Gamma$.
\item[$\Leftarrow)$]
If $0 = \lambda - \gamma$ for some $\lambda > \Lambda$, $\gamma < \Gamma$, then $\Lambda < \lambda = \gamma < \Gamma$.
\qedhere
\end{itemize}
\end{proof}
\begin{definizione}
The ordered group $\vg$ acts on $\vgh$ via the map $\Lambda \mapsto \gamma + \Lambda$.
Given $\Lambda \in \vgh$, the \intro{invariance group} of $\Lambda$
is $\IG{\Lambda} := \set{ \gamma \in \vg: \gamma + \Lambda = \Lambda}$, \ie
its stabiliser under the above action.
It is easy to see that $\IG\Lambda$ is a convex subgroup of~$\vg$.
Define the \intro{\width} of $\Lambda$ to be
$\ig\Lambda := {\IG\Lambda}^+ \in \vgh$.
The set of \widths of $\vgh$ is the set of cuts of the form $\ig\Lambda$,
as $\Lambda$ varies in~$\vgh$.
\end{definizione}
\begin{remark}
For every $\Lambda \in \vgh$, $\ig \Lambda = \Lambda - \Lambda > 0$.
Moreover, $\ig{\ig\Lambda} = \ig \Lambda$.
Besides, if $\vh$ is a submonoid of~$\vg$ and $\Lambda := {\vh}^+$,
then $\ig\Lambda = \Lambda$.
If moreover $\vh$ is a subgroup of~$\vg$, then $\IG\Lambda$ is the convex hull of~$\vh$; therefore, if $\vh$ is a convex subgroup of~$\vg$, then
$\IG\Lambda = \vh$.
%$\vh$ is a convex subgroup of~$\vg$, and $\Lambda := {\vh}^+$, then
%$\IG\Lambda = \vh$ and $\ig\Lambda = \Lambda$.%%
%\footnote{%
%If $\vh$ is not convex (but only a subgroup), then $\IG\Lambda$ is the convex hull of~$\vh$, and $\ig\Lambda =\Lambda$ is still valid.}
%If $\vh$ is just a submonoid of $\vg$ (not necessarily covex), nevertheless $\ig\Lambda =\Lambda$ holds.
%Moreover if $\vh$ is a subgroup,  then $\IG\Lambda$ is the convex hull of~$\vh$.}
Hence, there is a canonical bijection between convex subgroups of $\vg$ and
\widths of $\vgh$, sending $\vh$ to $\vh^+$
(whose inverse maps $\Lambda$ to $\IG \Lambda$).
\end{remark}
The above remark helps explain our choice of definition $\Lambda - \Gamma = \Lambda \ar (- \Gamma)$.
\begin{enumexamples}[Let $1\in\vg$ be a positive element, and $\Omega := \set{n\cdot 1:n\in\Nat}^+$.]
\item $\Omega + \Omega %= \Omega \ar \Omega 
= \Omega - \Omega = \Omega$.%%
\footnote{It can happen that $\Omega \ar \Omega > \Omega$.}
Moreover, $\Omega \dl \Omega = - \Omega$.
\item\label{EX:OMEGA-DISCRETE}
If $1$ is the immediate successor of $0$, then 
\[\begin{array}{r@{\,=\,}c@{\,=\,}l}
0^+ & 1^-    & \cut{(-\infty, 0]}{[1,+\infty)},\\
0^- & (-1)^+ & \cut{(-\infty, -1]}{[0,+\infty)}.
\end{array}\]
Moreover, $\underbrace{0^+ + \dotsb + 0^+}_{n\text{ times}} = 0^+$, while
\mbox{$\underbrace{0^+ \ar \dotsb \ar 0^+}_{n\text{ times}} = n^-$}.
\end{enumexamples}
\subsection{Group extensions}
\label{SUBSEC:EXTENSION}
Let $\vgp$ be an ordered group, and $\vg$ be a subgroup of~$\vgp$.
In this subsection we investigate the relationship between elements of $\vgp$ and cuts of~$\vg$.
The content of this subsection will not be used in the rest of the article; however, we hope that the former will help to clarify the latter.

In the rest of this subsection, the capital Greek letters will denote cuts of $\vg$; the small Greek letters, elements of $\vg$; and the small Latin letters, elements of~$\vgp$.
\begin{definizione}
For every $\Lambda \in \vgh$ and $x \in \vgp$,
we define
\[\begin{aligned}
\Lambda \leq x  &\text{ if }  \forall \lambda \in G\ \lambda < \Lambda \rightarrow \lambda < x,\\
\Lambda > x  &\text{ if } \Lambda \nleq x;
\end{aligned}\]
and similarly for $\Lambda \geq x$ and $\Lambda < x$.
We also say that $x$ satisfies $\Lambda$ (or that $x$ fills~$\Lambda$), and write $x \models \Lambda$,%
\footnote{We use the symbols $\models$ from model theory, because 
we can regard a cut as a set of formulae;
moreover, $x \models \Lambda$ iff $x$ satisfies the corresponding formulae $\lambda' < x < \lambda''$, as $\lambda'$ and $\lambda''$ vary in $\vg$, with $\lambda' < \Lambda < \lambda''$.}
iff
\[
%\begin{array}{rcl}x \models \Lambda &\text{iff} & 
\Lambda \leq x \text{ and } \Lambda \geq x.
\]
%Finally, we write
%\[
%\Lambda < x  \text{ iff } \Lambda \ngeq x.
%\end{array}
%\]
\end{definizione}
\begin{remark}
The following are equivalent:
\begin{enumerate}
\item $\Lambda < x$;
\item $
\exists \lambda \in \vg\ \Lambda < \lambda \leq x$;
\item $\Lambda \leq x$ and $x \not\models \Lambda$.
\end{enumerate}
Moreover, if $x \in \vg$, then $\Lambda < x$ in $\vgp$ iff $\Lambda < x$ in $\vg$.
\end{remark}
Note that if $x \models \Lambda$, then $x \in \vgp \setminus \vg$.
\begin{lemma}
\begin{enumerate}
\item If $\Lambda \leq x \leq \Gamma$, then $\Lambda \leq \Gamma$.
\item If $x < \Lambda \leq y$, then $x < y$.%
\footnote{$x \leq \Lambda \leq y$ does not imply that $x \leq y$.}
\item If $x \leq \Lambda \leq \Gamma$, then $x \leq \Gamma$.
\item If $x \leq y \leq \Lambda$, then $x \leq \Lambda$.
\item $x \geq \Lambda$ iff $-x \leq - \Lambda$.
\item If $x \leq \lambda < \Lambda$, then $x < \Lambda$.%
\footnote{$x < y \leq \Lambda$ does not imply that $x < \Lambda$.}
\end{enumerate}
\end{lemma}
\begin{proof}
Easy.
\end{proof}
\begin{lemma}
Let $\Lambda \leq x$ and $\Gamma \leq y$.
Then, $\Lambda + \Gamma \leq x + y$.
If moreover $x \models \Lambda $ and $y \models \Gamma$, then
\[\begin{array}{r@{\ \leq\ }c@{\ \leq\ }l}
\Lambda + \Gamma   & x + y   & \Lambda \ar \Gamma;\\
\Lambda \dl \Gamma & x - y  & \Lambda - \Gamma.
\end{array}\]
If moreover $\Lambda = \Gamma$ \rom(and $x, y \models \Lambda$\rom),
then $\abs{x - y} \leq \ig\Lambda$.\\
If instead $\Lambda < x$ and $\Gamma < y$, then $\Lambda + \Gamma \leq \Lambda \ar \Gamma < x + y$.%%
\footnote{It might happen that $\Lambda < x$, $\Gamma \leq y$, but $x + y \Lambda + \Gamma $.}
\end{lemma}
\begin{proof}
The first and last inequalities can be done by direct computation.
The others are direct consequences of the first one.
\end{proof}
\begin{lemma}
If $\Gamma \leq x$, then $\lambda + \Gamma \leq \lambda + x$.
If $\Gamma < x$, then $\lambda + \Gamma < \lambda + x$.
\end{lemma}
\begin{proof}
Trivial.
\end{proof}
\begin{lemma}
Assume that there exists $x_0 \in \vgp$ such that $x_0 \models \Lambda$.
Define
\[\begin{aligned}
S_L &:= \set{\alpha \in \vg: \exists x \in \vgp\ x \models \Lambda \et x \geq \alpha};\\
S_R &:= \set{\alpha \in \vg: \exists x \in \vgp\ x \models \Lambda \et x \leq \alpha};\\
T &:= \set{\alpha \in \vg: \forall x, y \in \vgp\ x,y \models \Lambda \rightarrow y - x \leq \alpha}.
\end{aligned}\]
Then, $\Lambda = S_L^+ = S_R^-$, and $\ig \Lambda = T^-$.
Moreover,
\[%\begin{aligned}
S_L = \set{ \alpha \in \vg: \alpha < x_0}, \qquad S_R = \set{\alpha \in \vg: \alpha > x_0}.
%\end{aligned}
\]
\end{lemma}
\begin{proof}
%Give proof.
First, we shall prove that $\alpha < x_0$ iff, for all $x \models \Lambda$, $\alpha < x$.
The ``if'' direction is obvious.
For the other direction, let $x \models \Lambda$, and suppose, for contradiction, that $x \leq \alpha < x_0$.
Since $x \geq \Lambda$, we have that $x > \beta$ for all $\beta < \Lambda$, and therefore $\alpha > \Lambda$.
Similarly, since $x_0 \leq \Lambda$, we have that $\alpha < \Lambda$, a contradiction.

%The following claim is obvious.
\begin{claim}
$\lambda < S_L^+$ iff $\lambda < x_0$.
%The following are equivalent:
%\begin{enumerate}
%\item $\lambda < S_L^+$;
%\item there exists $x \models \Lambda$ such that $\lambda < x$;
%\item for all $x \models \Lambda$ we have $\lambda < x$;
%\item $\lambda < x_0$.
%\end{enumerate}
\end{claim}
Let $\lambda < S_L^+$.
Then, there exists $\alpha \in S^L$ such that $\lambda \leq \alpha$.
Hence, $\lambda \leq x$ for every $x \models \Lambda$, and in particular $\lambda \leq x_0$.
Since we cannot have $\lambda \models \Lambda$, we conclude that
$\lambda < x_0$.
% for all $x \models \Lambda$.
Conversely, if $\lambda < x_0$, then, by definition, $\lambda \in S^L$, and therefore $\lambda < S_L^+$.

Let us prove now that $\Lambda = S_L^+$.
If $\lambda < \Lambda$, then $\lambda < x_0$; since $x_0 \models \Lambda$, we have that $\lambda < S_L^+$, and therefore $\Lambda \leq S_L^+$.\\
If $\lambda < S_L^+$, then $\lambda < x_0$.
Since $x_0 \models \Lambda$, we have that $\lambda < \Lambda$, and therefore $S_L^+ \leq \Lambda$.

$\Lambda = S_R^-$ is dual to $\Lambda = S_L^+$.

It remains to prove that $\ig\Lambda = T^-$.
The following claim is obvious.
\begin{claim}\label{CL:T}
$\alpha < T^-$ iff there exist $x,y \models \Lambda$ such that $y - x \geq \alpha$.
\end{claim}

Let $\alpha < \ig \Lambda$.
Since $\ig \Lambda$ is the upper bound of a group (specifically, the group~$\IG \Lambda$), we also have that $2 \alpha < \ig \Lambda$.
%This means that $\alpha < \lambda'' - \lambda'$, for every $\lambda' < \Lambda < \lambda''$.
We want to prove that $\alpha < T^-$, which would imply that $\ig \Lambda \leq T^-$.
By Claim~\ref{CL:T}, the thesis is equivalent to $\alpha \leq y - x$, for some $x, y \models \Lambda$.
Assume not.
Hence, $\alpha > 0$, and $x_0 - \alpha$ and $x_0 + \alpha$ do not satisfy~$\Lambda$.
Therefore, there exist $\lambda'$ and $\lambda'' \in \vg$, such that
\[
x_0 - \alpha \leq \lambda' < \Lambda < \lambda'' \leq x_0 + \alpha.
\]
Thus, $\lambda'' - \lambda' \leq 2 \alpha$, implying that $2 \alpha > \ig \Lambda$, a contradiction.

Conversely, let $\alpha > \ig \Lambda$.
Then, $\alpha > 0$.
Moreover, there exist $\lambda', \lambda'' \in \vg$, such that
$\lambda' < \Lambda < \lambda''$, and $\alpha \geq \lambda'' - \lambda'$.
%We have to prove that $\alpha > T^-$.
%Since $x_0 + \alpha \geq \lambda' + \alpha \geq \lambda'$, we have that
%$x_0 + \alpha > \Lambda$.
%Similarly, $x_0 - \alpha < \Lambda$.
Let $x , y \in \Lambda$, with $x \leq y$.
Since $\lambda' < x \leq y < \lambda''$, we have that
$\lambda'' - \lambda' > y - x$.
Hence, $\alpha > y - x$ for every $x,y \models \Lambda$.
Thus, by Claim~\ref{CL:T}, $\alpha > T^-$.
Therefore, $\ig \Lambda \geq T^-$.
\end{proof}
Assume that there exist $x_0, y_0 \in \vgp $ such that $x_0 \models \Lambda$ and $\Gamma \models y_0$.
We could try to define the sum of $\Lambda$ and $\Gamma$ using $x_0$ and $y_0$.
More precisely, we could define $z_0 := x_0 + y_0$, consider the cut $\Theta$ induced by $z_0$ on~$\vg$, and define $\Lambda \ad \Gamma := \Theta$.
However, $\Theta$~depends on our choice of $x_0$ and~$y_0$: different choices would produce different cuts, in general.
The canonical choice would be: define
\[\begin{aligned}
V   &:= \set{x + y: x \models \Lambda \et y \models \Gamma} \subseteq \vgp,\\
V_L &:= \set{\alpha \in \vg: \forall z \in V\ \alpha < z },\\
V_R &:= \set{\alpha \in \vg: \forall z \in V\ \alpha > z};
\end{aligned}\]
set $\Lambda \mathbin{\dot +^L} \Gamma := V_L^+$, 
$\Lambda \mathbin{\dot +^R} \Gamma := V_R^-$.
The next proposition implies that $\dot +^L = +$ and $\dot +^R = \ar$.
\begin{proposition}
Assume that there exist $x_0, y_0 \in \vgp $ such that $x_0 \models \Lambda$ and $y_0 \models \Gamma$.
Define $V$, $V_L$ and $V_R$ as above.
%\[\begin{aligned}
%V   &:= \set{x + y: x \models \Lambda \et y \models \Gamma} \subseteq \vgp.\\
%V_L &:= \set{\alpha \in \vg: \forall z \in V\ \alpha < z }\\
%V_R &:= \set{\alpha \in \vg: \forall z \in V\ \alpha > z}.
%\end{aligned}\]
Then, %$z \in V$ iff $\Lambda + \Gamma \leq z \leq \Lambda \ar \Gamma$.
%Moreover,
$\Lambda + \Gamma = V_L^+$, and $\Lambda \ar \Gamma = V_R^-$.
\end{proposition}
Note that if $z \in V$, then $\Lambda + \Gamma \leq z \leq \Lambda \ar \Gamma$.
However, the opposite implication might be false.
For instance, take $\vg = \Zed$, $\vg' = \frac{\Zed}{2}$, $\Lambda = \Gamma = 0^+$, and $z = 1/2$.
\begin{proof}
Let $\alpha < \Lambda + \Gamma$.
Hence, there exist $\lambda < \Lambda$ and $\gamma < \Gamma$ such that $\alpha \leq \lambda + \gamma$. If $x \models \Lambda$ and $y \models \Gamma$, then 
$\lambda < x$ and $\gamma < x$.
Therefore, $\alpha \leq \gamma + \lambda < x + y$.
Thus, $\alpha < V_L^+$.

Conversely, let $\alpha > \Lambda + \Gamma$.
Define $z_0 := x_0 + y_0$.
Suppose, for contradiction, that $\alpha < z$ for every $z \in V$, and let $v := z_0 - \alpha \geq 0$.
Let $x := x_0 - v$ and $y ;= y_0 - v$.
If $x \models \Lambda$, then $\alpha = x + y_0 \in V$, a contradiction.
Similarly, we reach a contradiction if $y \models \Gamma$.
Hence, there exist $\lambda$ and $\gamma$ such that $x \leq \lambda < \Lambda$ and $y \leq \gamma < \Gamma$.
Define $\theta := \alpha - (\lambda + \gamma)$.
\begin{claim}
$\lambda + \theta > \Lambda$ and $\gamma + \theta > \Gamma$.
\end{claim}
If, for instance, $\lambda + \theta < \Lambda$, then
$\alpha + \gamma = \lambda + \gamma + \theta < \Lambda+ \Gamma$, absurd.

Therefore, $\lambda + \gamma + 2 \theta > \Lambda \ar \Gamma$.
Thus,
\[
v = z_0 - \alpha < \lambda + \gamma + 2 \theta - \alpha = \alpha + \theta - \alpha = \theta,
\]
and we have $v < \theta$.
On the other hand,
$
x_0 + y_0 - 2 v = x +y \leq \lambda + \gamma = \alpha - \theta$,
thus
$\alpha + 2 v \geq z_0 + \theta = \alpha + v + \theta$,
and therefore $v \geq \theta$, absurd.
%We have $0 \leq \alpha \leq v$ by definition.
\end{proof}
%%
\begin{comment}
%%%%%%
\begin{proof}
We shall prove only that,
if $\Lambda + \Gamma \leq z \leq \Lambda \ar \Gamma$, then $z \in V$.
%there exist $x \models \Lambda$ and $y \geq \Gamma$.
The other assertions are easy consequences of this one.

Let $z_0 := x_0 + y_0 \in V$.
\Wlog, we can assume that $z \leq z_0$.
Hence, we have that $\Lambda + \Gamma \leq z \leq z_0$.
Let $v := z_0 - z \geq 0$, and $x := x_0 - v \in \vgp$.

\begin{claim}
$x \models \Lambda$.
\end{claim}
Suppose not.
Then, since $x \leq x_0 \models \Lambda$, we have that $x < \Lambda$.
Namely, there exists $\lambda' \in \vg$ such that $x \leq \lambda' < \Lambda$.
Hence,
\[
\Lambda + \Gamma \leq z = x + y_0 \leq \lambda' + y_0
\leq \lambda' + \Gamma < \Lambda + \Gamma,
\]
a contradiction.

Therefore, $z = x + y_0 \in V$, and we are done.
\end{proof}
%%%%%
\end{comment}
%%
\begin{enumexamples}
\item Let $\vg := \Raz$ and $\vgp := \Real$.
Every $r \in \Real \setminus \Raz$ determines a cut $\Theta$ on~$\Raz$, and $\Theta$ is the unique cut such that $r \models \Theta$.
\item Let $\vg := \Zed$ and $\vgp := \Real$.
Let $\Lambda = \Gamma = 0^+$.
Note that $x \geq \Lambda$ iff $x > 0$, and that $x \models \Lambda$ iff $0 < x < 1$.
Moreover, $V$~is the interval $(0,2) \subseteq \Real$.
Since $\Lambda + \Gamma = 0^+$ and $\Lambda \ar \Gamma = 2^-$, the proposition is verified.
\end{enumexamples}
%%

%% Sezione sulle estensioni di gruppi

%\part{\Doms}
\section[\Doms]{\Doms: Double ordered monoids}\label{SEC:DOMS}
\subsection{Basic definitions and facts}
\begin{definizione}
Let $\lang$ be the language $(\leq,\zO,+,-)$, where $\leq,\zO,+$, and $-$ are symbols for a binary relation, a constant, a binary function, and a unary function respectively.
A \intro{\predom{}} is an $\lang$-structure $\struct(M,\leq,\zO,+,-)$, where $\struct(M,\leq,\zO,+)$ is an Abelian ordered monoid with neutral element $\zO$,%
\footnote{I.e., $+$ is associative and commutative, with neutral element $\zO$, and $\leq$ is a linear order
satisfying $z\leq z' \rightarrow z + t \leq z' + t$ for every $z,z',t\in M$.}
and $\mathord{-} : M \to M$ is an anti-automorphism of $\struct(M,\leq)$ such that $-(-x) = x$.
For $x,y\in M$ define:
\[\begin{aligned}
x \ar y & := -\Pa{(-x) + (-y)},\\
x - y   & := x \ar (-y) = -\Pa{(-x) + y},\\
\zD     & := -\zO,\\
\ig x   & := x - x,\\
\abs x  & := \max(x, -x),\\
x \dl y & := x + (-y).
\end{aligned}\]
$\ig x = x - x$ is called the \intro{\width{}} of $x$.
\end{definizione}
%% M2: TODO: correggere!
\begin{remark}\label{REM:PREDOM}
Let $M$ be a \predom and $x$, $y$, $z \in M$. Then,
\[\begin{aligned}
x \ar \zD   &= x - \zO  = x;\\
\zD - x     &= -x;\\
-\ig x      &= x \dl x;\\
%\item\label{EN:DOM-4P} 
\ig{-x}     &= \ig x;\\
- (x - y)   &= y \dl x;\\
x\ar (y\ar z) &= (x\ar y)\ar z;\\
x - (y + z) &= (x - y) -z\\
            &= x \ar (-y) \ar (-z);\\
\text{if }x \leq \zD \text{ th}&\text{en } y \ar x \leq y\text{;}\\
\text{if }x \geq \zO \text{ th}&\text{en } y - x \leq y\text{;}\\
\text{if } x + y < x + z \text{ or } x \ar y < x \ar z, \text{ th}&\text{en } y < z\text{;}\\
\text{if } y \leq x \leq z, \text{ th}&\text{en } z - y \geq \ig x\text{.}
\end{aligned}\]
\end{remark}
Since $x \ar (y - z) = (x \ar y)- z$, in the following we often will write
$x \ar y - z$ for any of the above equivalent expressions.
Similarly, we will drop the parentheses in $x + (y + z)$ and $x \ar (y  \ar z)$.
\begin{remark}
If $M$ is a \predom, then $M^{\text{dual}} := \struct(M, \geq, \zD, \ar, -)$ is also a \predom, the \intro{dual} of $M$.
Hence, any theorem about \predoms has a dual theorem (the corresponding theorem for $M^{\text{dual}}$).
\end{remark}
\begin{definizione}
A \intro{\dom} $M$ is a \predom satisfying the following axioms: for every $x$, $y \in M$:
\begin{enumerate}[label={\textbf{M}\Alph*.}, ref={\textbf{M}\Alph*}]
\item\label{AX:M1} $\zD \leq \zO$;
\item\label{AX:M4} $\abs x \geq \zO$;
\item\label{AX:M5} $x \geq y$ iff $x - y \geq \zO$.
%$x < y$ iff $x - y < \zO$.
\end{enumerate}
\end{definizione}
Note that the above axiomatisation is universal.
\begin{remark}
For a \predom, Axiom~\ref{AX:M4} is equivalent to the fact that the interval $(\zD, \zO)$ is empty
(if $\zD > \zO$, the interval $(\zD, \zO)$ is empty by default).
\end{remark}
%% \begin{proof}
%% \begin{itemize}
%% \item[$\Rightarrow$)]
%% If, for contradiction, $\zD < x <\zO$, then $\zO > -x > \zD$, hence $\abs x < \zO$.
%% \item[$\Leftarrow$)] If $x \geq \zO$, then $\abs x \geq x \geq \zO$.
%% If $x < \zO$, then $x \leq \zD$, hence $-x > \zO$,
%% and $\abs x \geq x \geq \zO$. \qedhere
%% \end{itemize}
%% \end{proof}
%%
\begin{remark}\label{REM:M5}
For a \predom, Axiom~\ref{AX:M5} is equivalent to (the universal closure of) any of the following:
\begin{enumerate}
\item\label{AX:M5-2} 
$x < y$ iff $x - y < \zO$;
%$x \geq y$ iff $x - y \geq \zO$;
\item\label{AX:M5-3} $x \dl y > \zD$ iff $x > y$;
\item\label{AX:M5-4} $x - y < z$ iff $x < y + z$;
\item\label{AX:M5-5} $x - y \geq z$ iff $x \geq y + z$.
\end{enumerate}
\end{remark}
\begin{proof}
Applying \eqref{AX:M5-2} twice,
%Axiom~\ref{AX:M5}, 
we get $x - y < z$ iff $(x - y )- z = x - (y+z)< \zO$ iff $x < y + z$.
Therefore, \eqref{AX:M5-2} implies~\eqref{AX:M5-4}.
The rest is easy.
\end{proof}
Hence, in the following we will refer to any of the aforesaid
equivalent forms as Axiom~\ref{AX:M4}, \ref{AX:M5} respectively.
\begin{proviso*}
For the rest of the article, $M$~will be a \dom and $\vg$ a (linearly) ordered (Abelian) group, unless we say otherwise.
%% We recall that $\vgh$ is the set of cuts of $\vg$.
\end{proviso*}
\begin{remark}
The axioms~\ref{AX:M1}, \ref{AX:M4} and \ref{AX:M5} are self-dual.
Therefore, if $M$ is a \dom, then its dual is also a \dom.
Hence, any theorem about \doms has a dual theorem.
\end{remark}
\begin{proof}
%The only non-trivial axiom is the dual of~\ref{AX:M5}:
%~$\text{\ref{AX:M5}}^{\text{dual}}$:
For instance, the dual of Axiom~\ref{AX:M5} is
$
x \leq y \text{ iff } x \dl y \leq \zD
$, \ie
$
x \leq y \text{ iff } y - x \geq 0
$,
which is axiom~\ref{AX:M5} itself.
%For instance, the fact that Axiom~\ref{AX:M5} is self-dual is contained in Remark~\ref{REM:M5}.
%For instance, the dual of Axiom~\ref{AX:M5} is Remark~\ref{REM:M5}(\ref{AX:M5-3}).
\end{proof}
\begin{enumexamples}
\item\label{EX:GROUP}
Any ordered group is a \dom, with $\zD = \zO$ and $x + y = x \ar y$.
\item We have seen that $\vgh$ is a \dom (cf.\ Lemma~\ref{LEM:CUTS-AXIOMC}).
\item\label{EX:TRIVIAL} The trivial models. Let $N$ be an ordered set with minimum $\zO$.
Define $M$ as the disjoint union of two copies of $N$, \ie\ $M := -N \sqcup N$, with the reversed order on $-N$, and the rule $-N < N$, and the minus defined in the obvious way.
Define
\[
x + y := 
\begin{cases}
  x & \text{if } \abs x > \abs y,\\
  y & \text{if } \abs x < \abs y,\\
  \min(x,y) & \text{if } \abs x = \abs y.
\end{cases}\]
It is easy to see that $M$ is a \dom.
Moreover, $\zD < \zO$, and  for every $x,y\in M$, 
\[
x \ar y = 
\begin{cases}
  x & \text{if } \abs x > \abs y,\\
  y & \text{if } \abs x < \abs y,\\
  \max(x,y) & \text{if } \abs x = \abs y,
\end{cases}\]
and $\ig x = \abs x$.
\item\label{EX:TRIVIAL-2} The above example can be modified by identifying $\zO$ with $\zD$.
The resulting structure (with the same definition of operations and relations) is also a \dom, satisfying $\zD = \zO$ and $\ig x = \abs x$.
\item\label{EX:G-TILDE}
Define $\vgt$ to be the disjoint union of $\vg$ and $\vgh$,  with order and operations extending the ones on $\vg$ and $\vgh$ in the way defined in \S\ref{SEC:GROUP-BASIC}.
Note that $\zO$, the neutral element of~$\vgt$, is the neutral element of~$\vg$.
Moreover, $\vgt$ is a \dom with $\zD = \zO$.
\item\label{EX:GROUP-DELTA} Let $\struct(\vg, \leq, 0, +)$ be an ordered group, and use $\md$ to denote its minus operation.
Fix  $\delta \in \vg$. Define $M := \vg$, with the same order, neutral element and plus. However, define $-x := (\md x) + \delta$.
Then, $x \ar y = -(-x + -y) = \md(\md x + \md y + 2 \delta) + \delta = x + y \md \delta$.
Moreover, $x - y = x \md y$, and $\ig x = 0$, and $-0 = \delta$.
Therefore, $M$ is a \predom.
Moreover, if $\delta > 0$, $M$ satisfies
axioms~\ref{AX:M4} and \ref{AX:M5}, while if $\delta < 0$, $M$ satisfies axioms~\ref{AX:M1}
and~\ref{AX:M5}.
If $\delta = 0$, we retrieve Example~\ref{EX:GROUP}.
Note that we have $0 \ar 0 = \md\delta$.
Hence, if $\delta > 0$, we have $0 \ar 0 < 0$.

Moreover, if $\delta < 0$, then $M$ is a \dom iff the interval $(\delta, 0)$ is empty, \ie\ iff $\delta = -1$, where $1$ is the minimal positive element of $\vg$.
For instance, if $\vg = \Zed$ (and $\delta = -1$), then $M = \Zedhat \setminus \set{\pm \infty}$.
In the general case, $M$ is the subset of $\vgh$ of elements with invariance group $\{0\}$ (again, if $\delta = -1$).
\end{enumexamples}
%%Examples~\ref{EX:G-TILDE} and~\ref{EX:GROUP-DELTA}
%%
\begin{definizione}
Given $n \in \Nat$ and $x,y \in M$, define
\begin{align*}
x - n y & := (\dotso(x \underbrace{- y) -y) - \dotsb) - y}_{n \text{ times}}, &%\quad
x + n y & := x \underbrace{+ y + \dotsb + y}_{n \text{ times}}.\\
\intertext{In particular, $x - 0 y = x + 0 y = x$.
%% M2: TODO: Dice di non restringersi a Nstar, quindi dice una michiata...
Moreover, given $n \in \Nats$, define}
n x    &:= \underbrace{x + \dotsb + x}_{n\text{ times}}, &%\quad
(-n)x  &:= -(n x) = \underbrace{-x - \dotsb -x}_{n \text{ times}}.
\end{align*}
\end{definizione}

It is easy to see that $(x - n y) - m y = x - (n+m) y$, and $(x + n y) + m y =\linebreak[2] x + (n+m) y$.
In general, $n(-x) \neq -(n x) = (-n)x$.
Take for instance $M = \wideDed{\Zedtwo}$, $x = \unmezzo$, $n = 2$.
However, if $\vg$ is divisible or $n = 1$, then $n(-\Gamma) = - (n\Gamma)$ for every $\Gamma \in \vgh$.
\begin{proposition}\label{PROP:DOM}
Let $w$, $x$, $y$, $z \in M$. Then,
\begin{enumerate}
\item\label{EN:DOM-17} $\ig x \geq \zO$; % ex axiom M3
\item\label{EN:DOM-18} $\zO \ar \zO \geq \zO$ and $\zD + \zD \leq \zD$;
\item\label{EN:DOM-19} $x + y \leq x \ar y$; % ex axiom M2
\item\label{EN:DOM-4} $(x + y) - y \geq x \geq (x - y) + y$;
\item\label{EN:DOM-20} $x - y = \max\set{z \in M: y + z \leq x}$;%
\footnote{G~.Birkhoff~\cite[XIV\S5]{BIRKHOFF:1967} calls such an element the residual of $x$ by $y$.}
\item\label{EN:DOM-4PP} $((x + y) - y) + y = x + y$ and 
$((x - y) + y) -y = x - y$;
\item\label{EN:DOM-5} $x + \ig x = x = x - \ig x$;
\item\label{EN:DOM-5P} $y > \ig x$ iff $x + y > x$;
\item\label{EN:DOM-5PP} $\ig x \leq \abs x$;
\item\label{EN:DOM-7} $x - y \geq \zO$ and $y - x \geq \zO$ iff $x = y$;
\item\label{EN:DOM-9} $(x + y) \ar z \geq x + (y \ar z)$;
\item\label{EN:DOM-10} $(x + y) \ar z \ar w \geq (x \ar z) + (y \ar w)$, and
 $(x + y) - (z + w) \geq (x -z) + (y - w)$;
\item\label{EN:DOM-10P} %$(x + y) \ar (z + w) \geq x + z + (y \ar w)$\\
$(x \ar y) + z + w \leq (x + z) \ar (y + w)$;
\item\label{EN:DOM-10PP} if $x + y < x \ar z$, then $y \leq z$;
%Quondam Axiom 6
\item\label{EN:DOM-8} if $x + y < x \ar y$, then $\ig x = \ig y$;
\item\label{EN:DOM-13P} if $x, y < \zO$, then $x \ar y < \zO$;
\item\label{EN:DOM-13PP} if $x < z$ and $y < w$, then $x \ar y < z + w$,
 and \emph{a fortiori} $x + y < z + w$;%
\footnote{From $y < w$, we can only conclude that $x \ar y \leq x + w$ (it is the contrapositive of~\ref{EN:DOM-10PP}), but not even that $x + y < x \ar w$.
For instance, take $\vg = \Zedtwo$, $x = \unmezzo$, $y = 0^-$, $w = 0^+$;
then, $x + y = x \ar w = \unmezzo$.}
%in particular, if $z' < x < z$, then $z \dl z' > \ig x$;
in particular, if $t < z < x$, then $x \dl t > \ig z$ (by taking $y:=-z$ and $w:=-t$);
\item\label{EN:DOM-16} $\ig x + \ig x = \ig x$;
\item\label{EN:DOM-16PP} if $2x = x$, then $\ig x = \abs x$;%
\footnote{The converse is not true. Take for instance $M = \Zedhat$ and $x = 0^-$.}
\item\label{EN:DOM-16P} $\igigx = \ig x$;
\item\label{EN:DOM-8P} $\ig{x + y} = \ig{x \ar y} = \max(\ig x, \ig y)$;
\item\label{EN:DOM-11} %if $\ig x = \ig y$, then
$x \ar y \leq (x + y) \ar \ig x$ \rom(and the same for $\ig y$\rom);
\item\label{EN:DOM-11P} if $\ig x > \ig y$, then $x \ar \ig y = x$;
%\item\label{EN:DOM-11PP} if $\ig x \neq \ig y$, then  $x + y = x \ar y$;
\item\label{EN:DOM-12} if $y > x$, then $y \geq x \ar \ig y$;
\item\label{EN:DOM-12PP} the intervals $(x, x \ar \zO)$ and $(x + \zD, x)$ are empty;
\item\label{EN:DOM-14} if $x, y \leq \ig z$, then $x + y \leq \ig z$;
\item\label{EN:DOM-15} if $x, y < \ig z$, then $x \ar y < \ig z$.
\end{enumerate}
\end{proposition}
%We have already proved the above proposition in the case $M = \vgh$ for some ordered group $\vg$.
\begin{proof}
\begin{prooflist}
\item[\ref{EN:DOM-17}] If, by contradiction, $x - x < \zO$, then, by Axiom~\ref{AX:M5}, $x < x$, absurd.
\item[\ref{EN:DOM-18}] If, by contradiction, $\zO \ar \zO < \zO$, then, by Axiom~\ref{AX:M5}, $\zO < \zD$, contradicting Axiom~\ref{AX:M1}.
%$\zO < \zO + \zO = \zO$, absurd.
The other inequality is the dual one.
\item[\ref{EN:DOM-19}] If, by contradiction, $x + y > x \ar y$, then $x \ar y -  (x + y) < \zO$, and it follows $\ig x \ar \ig y < \zO$.
The conclusion follows from~\ref{EN:DOM-17} and~\ref{EN:DOM-18}.
\item[\ref{EN:DOM-4}] If, by contradiction, $(x + y) - y < x$, then,
by Axiom~\ref{AX:M5}, $x + y < x + y$, absurd.
Similarly for the other inequality.
\item[\ref{EN:DOM-20}] By Axiom~\ref{AX:M5}, $y + z \leq x$ iff $y \leq x - z$.
\item[\ref{EN:DOM-4PP}] By~\ref{EN:DOM-4}, the pair of maps $x \mapsto x \pm y$ forms a Galois connection between the ordered set $M$ and its dual.
The conclusion is true for any such correspondence \cite[Theorem~IV.5.1]{MACLANE:1998}.
%\footnote{See S. MacLane \textit{Categories for the Working Mathematician}, Theorem~IV.5.1.}
More in detail, applying \ref{EN:DOM-4} twice, we get
\[
((x \underbrace{+ y) - y}_{\text{simplify}}) + y  \geq x + y  \geq ((x + y) \underbrace{- y) + y}_{\text{simplify}}.
\]
The other equality is the dual one.
\item[\ref{EN:DOM-5}]
By~\ref{EN:DOM-17}, $x + \ig x \geq x$.
If, by contradiction, $x < x + \ig x $, then, %by~\ref{EN:DOM-3P},
by Axiom~\ref{AX:M5}, $\ig x < \ig x$, absurd.
The other equality is the dual one.
\item[\ref{EN:DOM-5P}] %By~\ref{EN:DOM-3P},
By Axiom~\ref{AX:M5},
$x - x < y$ iff $x < x + y$.
\item[\ref{EN:DOM-5PP}] Since $\ig{-x} = \ig x$, \wloG we can assume that $x = \abs x $, hence, by Axiom~\ref{AX:M4}, $x \geq \zO$.
Therefore, $\ig x = x - x$ which is $\leq x = \abs x$ by Remark~\ref{REM:PREDOM}.
\item[\ref{EN:DOM-7}] Immediate from Axiom~\ref{AX:M5}.
\item[\ref{EN:DOM-9}] If, by contradiction, $(x + y) \ar z < x + (y \ar z)$,
then $(x + y) \ar z - x < y \ar z$; however \ref{EN:DOM-4} implies
$((x+y)-x)\ar z \geq y\ar z$, which is absurd.
\item[\ref{EN:DOM-10}] If, by contradiction, $(x + y) \ar z \ar w < (x \ar z) +
(y \ar w)$, then $(x + y) \ar w - (y \ar w) \ar z < x \ar z$, hence, by
Remark~\ref{REM:PREDOM}, $(x + y) \ar w - (y \ar w) < x$, thus  by
Axiom~\ref{AX:M5} $(x + y) \ar w <
x + (y \ar w)$, contradicting \ref{EN:DOM-9}.
\item[\ref{EN:DOM-10P}] Dual of \ref{EN:DOM-10}.
\item[\ref{EN:DOM-10PP}]
Assume, for contradiction, that $y > z$, \ie\ $z - y < \zO$.
The hypothesis is equivalent to $(x + y) - (x \ar z) < \zO$,
i.e.\ $(-x - y) + (x \ar z) > \zD$. Since by \ref{EN:DOM-10P} $(-x - y) + (x \ar z)\leq z - y - \ig x$, we have 
by Axiom~\ref{AX:M4} $z - y - \ig x \geq \zO$.
% which, by \ref{EN:DOM-19} and Axiom~\ref{AX:M4}, implies $z - y - \ig x \geq \zO$. 
Therefore, by \ref{EN:DOM-17}, $z - y \geq \zO$, absurd.
\item[\ref{EN:DOM-8}] Assume, by contradiction, that $\ig x < \ig y$.
Then, $\ig x - \ig y < \zO$.
Moreover, $x + y < x \ar y$ implies that $(x \ar y) \dl (x + y) > \zD$, thus $(x \ar y) + (-x - y) \geq \zO$ (we used Axiom~\ref{AX:M4}).
Therefore, by~\ref{EN:DOM-10}, $(x - x) \ar (y \dl y) \geq \zO$, hence $\ig x - \ig y \geq \zO$, a contradiction.
\item[\ref{EN:DOM-13P}] Since $y < \zO$, $y \leq \zD$, thus $x \ar y \leq x < \zO$.
\item[\ref{EN:DOM-13PP}] The hypothesis is equivalent to $x - z < \zO$ and $y - w < \zO$, which, by~\ref{EN:DOM-13P}, 
implies $(x - z) \ar y - w < \zO$, \ie\ $(x \ar y) - (z + w) < \zO$, which is equivalent to the conclusion.
\item[\ref{EN:DOM-16}] Since $\ig x \geq \zO$, $\ig x + \ig x \geq \ig x$.
If, by contradiction, $\ig x + \ig x > \ig x$, then $x + \ig x + \ig x > x$, contradicting~\ref{EN:DOM-5}.
\item[\ref{EN:DOM-16PP}] 
Suppose, for contradiction, that $2x = x$, but $\ig x < \abs x$ (by~\ref{EN:DOM-5PP}).
W.l.o.g.\ we may assume $x\geq\zO$, because if not then $x':=-x$ satisfies the same hypothesis as $x$
(in particular $2x'=x'$ is a consequence of $x'\ar x'=x$ since $x' \leq 2x' \leq x'\ar x'$).
Then, $- x < \ig x < x$, hence $\ig x - x < x + \ig x $, thus $x - 2x < \ig x$, absurd.
\item[\ref{EN:DOM-16P}] $\igigx \leq \ig x$ by \ref{EN:DOM-5PP} and~\ref{EN:DOM-17}.
If, by contradiction, $\ig x - \ig x < \ig x$, then $\ig x < \ig x + \ig x$, contradicting~\ref{EN:DOM-16}.
\item[\ref{EN:DOM-8P}]
By~\ref{EN:DOM-10}, $\ig{x + y} = (x + y) - (x + y) \geq \ig x + \ig y \geq \max(\ig x,\ig y)$.
\Wlog, $\ig x \geq \ig y$.
By~\ref{EN:DOM-16}, $\ig x + \ig y = \ig x$.
Assume, for contradiction, that $\ig x + \ig y < \ig{x + y}$.
Then, by~\ref{EN:DOM-5P},
\[\begin{aligned}
x + \ig{x + y} &> x,\\
y + \ig{x + y} &> y.
\end{aligned}\]
Therefore, by~\ref{EN:DOM-13PP}, $x + y + 2(\ig{x + y}) > x + y$,
contradicting~\ref{EN:DOM-5}.
%~\ref{EN:DOM-16} and
The proof that $\ig{x \ar y} = \ig x$ is similar.
\item[\ref{EN:DOM-11}] If, by contradiction, $x \ar y > (x + y) \ar x - x$, then $y > (y + x) - x$, contradicting~\ref{EN:DOM-4}.
\item[\ref{EN:DOM-11P}] By hypothesis, $\ig y - \ig x < \zO$.
If, by contradiction, $x \ar \ig y > x$, then $x - (x \ar \ig y) < \zO$, \ie $x \ar \Pa{(-x) \dl \ig y} < \zO$.
Therefore, $(x - x) \dl \ig y < \zO$, \ie $\ig x \dl \ig y < \zO$, \ie $\ig y - \ig x \geq \zO$, absurd.
%\item[\ref{EN:DOM-11PP}] \Wlog, $\ig x > \ig y$, hence $\ig{x + y} = \ig x > \ig y$.
%Thus, by~\ref{EN:DOM-11} and~\ref{EN:DOM-11P}, $x + y \leq x \ar y \leq (x + y) \ar \ig y = x + y$.
\item[\ref{EN:DOM-12}]
If, by contradiction, $y < x \ar \ig y$, then $y < y \ar (x - y)$, hence $x - y \geq \zO$, contradicting $x < y$.
\item[\ref{EN:DOM-12PP}] Immediate from~\ref{EN:DOM-12}.
\item[\ref{EN:DOM-14}] Immediate from~\ref{EN:DOM-16}.
\item[\ref{EN:DOM-15}] By hypothesis, $x - \ig z < \zO$ and $y - \ig z < \zO$.
Hence, $ (x - \ig z) \ar y - \ig z < \zO$, \ie\ $(x \ar y) -(\ig z + \ig z) < \zO$, hence $(x \ar y) - \ig z < \zO$, and the conclusion follows.
\qedhere
%\end{enumerate}
\end{prooflist}
\end{proof}
\begin{corollary}\label{COR:D-K-M}
Let $d$, $k$, $m\in\Nat$, with $k < m$, and $x$, $y \in M$.
Then,
\[
\bigl( x - (m + d) y \bigr) + (m y - k y) \leq x - (d + k) y.
\]
\end{corollary}
\begin{proof}
If we define $x' := x - d y$, we see that it is enough to treat the case $d = 0$.
By Proposition~\ref{PROP:DOM}(\ref{EN:DOM-10}),
%Lemma~\ref{LEM:SUM-EXCHANGE},
\[
(x - m y) + (m y - k y) \leq (x - k y) \ar (m y \dl m y) = (x - k y) - \ig y = 
x - k y. \qedhere
\]
\end{proof}
In general, we do not have equality: take for instance $M = \Razhat$, $x = 0^+$, $y = 0^-$, $k = 0$, $m = d = 1$.
Then, the left hand side is equal to $0^-$, while the \rhs is equal to $0^+$.
Another counter-example, this time with $k >0$: take $M = \wideDed{\Zedtwo}$, $x = y = \unmezzo \in \Raz \setminus \Zedtwo$, $d = 0$, $k = 1$, $m = 2$.
Then, the \lhs is equal to $0^-$, while the \rhs is $0^+$.
With the same $M$, $x$, and $y$, we could also take $d = m = 1$ and $k = 0$.
%% TODO: delete!
\begin{corollary}\label{COR:D-K-M-J}
Let $x,x',y\in M$, $j,j',k,m,d \in \Nats$ such that $j,j',k < m$, and $j + j' = m + d$. Then,
\[
(x - j y) + (x' - j' y) + (m y - k y) \leq (x - x') - (d - k) y.
\]
\end{corollary}
\begin{proof}
The \lhs is less or equal to $\Pa{(x + x') - (m + d)y} + (my - ky)$, hence the conclusion is immediate from Corollary~\ref{COR:D-K-M}.
\end{proof}
\subsection{Sub-\doms and \domhoms}
\begin{definizione}
An element $x \in M$ is a \intro{\width{}} element (or \width for short) if it is equal to its \width (\ie, $\ig x = x$).
Define $\orders := \Orders{M}$ as the set of \widths of $M$:
it is an ordered subset of~$M$, and $\ \ig{}\ $ is a surjective map from $M$ to $\orders$.
Note that $\zO$ is in $\orders$.

Given $a \in \orders$, define
\[%\begin{aligned}
\eqdom{M}{a} := \set{x \in M : \ig x = a}, \quad
\gqdom{M}{a} := \set{x \in M : \ig x \geq a},
%\gdom{M}{a}  &:= \set{x \in M : \ig x > a},\\
%\lqdom{M}{a} &:= \set{x \in M : \ig x \leq a},\\
%\ldom{M}{a}  &:= \set{x \in M : \ig x < a}.
%\end{aligned}
\]
and similarly for $\gdom{M}{a}$, $\lqdom{M}{a}$ and $\ldom{M}{a}$.
We shall write $\Mgroup$ instead of $\eqdom{M}{\zO} = \set{x \in M: \ig x = \zO}$.
More generally, for any $S\subseteq \orders$, define
\[
M^{\{S\}} := \set{ x \in M: \ig x \in S}.
\]
\end{definizione}
\begin{definizione}
A subset~$A$ of a \predom~$M$ is \intro{symmetric} if $-A = A$.
It is a \intro{\qsubdom }of~$M$ if it is symmetric and $A + A \subseteq A$.
% and $-A \subseteq A$.
If moreover $\zO \in A$, then $A$ is a \intro{\subdom }of~$M$, and $M$ is a \intro{\overdom }of~$A$.
\end{definizione}
\begin{definizione}
A function $\phi: M \to M'$ between two structures $\struct(M, \leq, +, -)$
and $\struct(M', \leq, +, -)$ is a \intro{\qdomhom} if $\phi$ preserves the structure, \ie $\phi(x + y) = \phi(x) + \phi(y)$, $\phi(-x) = - \phi(x)$, and if $x \leq y$ then $\phi(x) \leq \phi(y)$.

A function $\phi$ between two {\domstruct}s $M$ and $M'$ is a \intro{\domhom}
if it is a \qdomhom and $\phi(\zO) = \zO$.

The \intro{kernel} of such a \domhom~$\phi$ is
\[
\ker (\phi) := \set{x \in M: \phi(x) = \zO}.
\]
\end{definizione}
Note that a \subdom of a \dom is indeed a \dom (because the axiomatisation of \doms is universal) and that the corresponding inclusion map is a \domhom.
Moreover, a subset of a \predom is a \qsubdom iff the corresponding inclusion map is a \qdomhom.
\begin{enumexamples}
\item Let $M$ and $M'$ be two ordered groups, considered as \doms with the group minus.
A function $\phi: M \to M'$ is a \domhom iff $\phi$ is a homomorphism of ordered groups, and a \subdom of $M$ is the same as a subgroup.
%\item\label{EX:QSUBDOM} The notion of \qsubdom is important in the following case:
%let $\vh$ be a convex subgroup of $\vg$, and let $F$ be the quotient $\vg/\vh$.
%Then, the quotient map $\vg \to F$ induces a bijection between the cuts $\Lambda$ of $\vg$ such that $\vh \subseteq \IG\Lambda$, and the cuts of~$\vf$ (cf. Lemma~\ref{LEM:QUOTIENT-GROUP}).
%Via this map, $\vfh$ can be identified with a \qsubdom of $\vgh$.
\item Let $M := \Zedhat \setminus\set{\pm\infty}$, and consider $\Zed$ as a \dom (with the group minus).
The map from $\Zed$ into $M$ sending $x$ to $x^+$ is an isomorphism of ordered groups, but it is not a \domhom, since it does not preserve the minus (cf. Example~\ref{EX:GROUP-DELTA}).
\item\label{EX:M-INFTY} Let $-\infty$ and $+\infty$ be two elements not in $M$.
Define the \dom $M^{\infty} := M \sqcup \set{-\infty, + \infty}$ with the operations and relation extending the ones on $M$,
and, for every $x \in M$,
\[\begin{aligned}
-(+\infty)  &= -\infty,\\
-\infty     &< x < + \infty,\\
-\infty + x &= -\infty, \\
+\infty + x &= +\infty, \\
(-\infty) + (+\infty) &= -\infty.
\end{aligned}\]
It is easy to see that $M^\infty$ is indeed a \dom, and $M$ is a convex sub\dom of $M^\infty$.
\end{enumexamples}
\begin{corollary}
For every $a \in \orders$, the set $\set{y \in M: \abs y \leq a}$
is a convex \subdom of~$M$.
Moreover, if $a > \zO$, then the set $\set{y \in M: \abs y < a}$ is also a convex \subdom of~$M$.%
\footnote{If $a = \zO$, it is the empty set.}
\end{corollary}
\begin{proof}
By Proposition \ref{PROP:DOM}(\ref{EN:DOM-14}, \ref{EN:DOM-15}).
\end{proof}
Let $\vh$ be a convex subgroup of $\vg$.
The quotient group $\vg/\vh$ is, in a canonical way, an ordered group, with the definition $\gamma/\vh \leq \lambda /\vh$ if there exists $\eta \in \vh$ such that $\gamma \leq \eta + \lambda$.
Equivalently, $\gamma /\vh < \lambda/\vh$ if for every $\eta\in \vh$, $\gamma < \eta + \lambda$.

Fix a convex subgroup $H$ of $\vg$, and denote by $\pi: \vg \to \vg/H$ the quotient map.
\begin{lemma}\label{LEM:QUOTIENT-GROUP}
Let $\Lambda \in \vgh$ such that $H \subseteq \IG\Lambda$.
Then,
\[
\pi^{-1}(\pi\Lambda^L) = \Lambda^L \quad \text{and} \quad
\pi^{-1}(\pi\Lambda^R) = \Lambda^R.
\]
Hence, the map $\pi$ induces a bijection  between the cuts $\Lambda$ of $G$ such that $H \subseteq \IG\Lambda$, and the cuts of $\vg / H$.
Denote by $\Lambda/H$ the cut of $\vg/H$ induced by $\Lambda$.

Moreover, if $\gamma \in \vg$ and $\Gamma$ is a cut such that $H \subseteq \IG\Gamma$, we have:
\begin{enumerate}
\item $\Lambda \gtrless \gamma$ iff $\Lambda /H \gtrless \gamma/H$;
\item $\Lambda \gtreqqless \Gamma$ iff $\Lambda /H \gtreqqless \Gamma/H$;
\item $(\Lambda + \Gamma)/H   = \Lambda /H + \Gamma /H$;
\item $(\Lambda \ar \Gamma)/H = \Lambda /H \ar \Gamma /H$;
\item $(\Lambda \dl \Gamma)/H = \Lambda /H \dl \Gamma /H$;
\item $(\Lambda - \Gamma)/H   = \Lambda /H - \Gamma /H$;
\item $(\gamma + \Gamma)/H   = \gamma /H + \Gamma /H$.
\end{enumerate}
\end{lemma}
\begin{proof}
Let $\gamma \in \pi^{-1}(\pi\Lambda^L)$. %and $\gamma > \Lambda^L$.
Then, there exists $\theta \in H$ such that $\gamma - \theta \in \Lambda^L$,
\ie\ $\gamma \in \theta +\Lambda^L$, \ie\ $\gamma < \theta + \Lambda$.
However, since $\theta \in \IG \Lambda$, $\theta + \Lambda = \Lambda$, thus $\gamma < \Lambda$. Similarly for $\gamma \in \pi^{-1}(\pi\Lambda^R)$.

The first point is a consequence of the above, and the second of the first.
%In the same way we can prove that $\Lambda \gtreqless \Gamma$ iff $\Lambda \gtreqless \Gamma/H$.

Let us prove, for instance, the third point (the others are similar).\\
Let $\alpha /H < (\Lambda + \Gamma) /H$.
Hence, $\alpha < \Lambda + \Gamma$, so  $\alpha = \lambda + \gamma$ for some $\lambda < \Lambda$, $\gamma < \Gamma$, thus $\alpha /H = \lambda  /H + \gamma  /H < \Lambda /H + \Gamma /H$.
Therefore, $(\Lambda + \Gamma) /H \leq \Lambda /H + \Gamma /H$.\\
Let $\alpha /H < \Lambda /H + \Gamma /H$.
Hence, $\alpha /H = \lambda /H + \gamma /H$ for some $\lambda < \Lambda$, $\gamma < \Gamma$, thus $\alpha + \theta = \lambda + \gamma$ for some $\theta \in H$, \ie $\alpha + \theta < \Lambda + \Gamma$, so $\alpha < \Lambda + \Gamma$.
Therefore, $(\Lambda + \Gamma) /H \geq \Lambda /H + \Gamma /H$.
\end{proof}

\begin{remark}\label{REM:PI-HAT}
Let $\pi: \vg \to \vk$ be a surjective homomorphism between two ordered groups.
By the above lemma, $\pi$ induces a map $\pihat: \vkh \to \vgh$, via the formula
\[
\pihat(\Lambda) := 
\Pa{\pi^{-1}(\Lambda^L)}^+ = \Pa{\pi^{-1}(\Lambda^R)}^-
= \cut{\pi^{-1}(\Lambda^L)}{\pi^{-1}(\Lambda^R)}.
\]
Moreover, $\pihat$ is an injective \qdomhom.
Finally, $\Dedekind{}$ is a contravariant functor between the categories of
ordered (Abelian) groups with surjective homomorphisms,
and \doms with injective \qdomhoms.
\end{remark}
Note that if $\pi$ is not surjective, then $\pihat$ can be still defined by the same formula, but it will not preserve the sum.
\begin{corollary}\label{COR:MGROUP}
$\Mgroup$ is a \subdom of~$M$.
More generally, for any $S\subseteq \orders$,
$M^{\{S\}}$ is a \qsubdom of $M$.
Moreover, if $\zO \in S$, then $M^{\{S\}}$ is a \subdom of~$M$.
\end{corollary}
\begin{enumexamples}
\item Let $a$ be a \width.
Then, by Corollary~\ref{COR:MGROUP}, $\gqdom{M}{a}$ and $\gdom{M}{a}$ are \qsubdoms of $M$.
By Proposition~\ref{PROP:DOM}(\ref{EN:DOM-11P} and \ref{EN:DOM-5PP}), $\gqdom{M}{a}$ is actually a \dom, with neutral element $\ig a$.
However, the only case when $\gqdom{M}{a}$ is a \subdom of $M$ is when $a = \zO$, and in that case $\gqdom{M}{a} = M$.
\item Fix a \width $a$.
Then, $\lqdom{M}{a}$ and $\ldom{M}{a}$ are \subdoms of $M$ (unless $a = \zO$, because in that case the second set is empty).
\end{enumexamples}
\begin{corollary}\label{COR:KERNEL}
Let $\phi: M \to N$ be a \domhom.
Then, $\ker (\phi)$ is a convex submonoid of $M$.
If moreover $\zD = \zO$ in $N$, then  $\ker (\phi)$ is a convex sub\dom of $M$.
Moreover, if $y \in N$, then $\phi^{-1}(y) + \ker (\phi) \subseteq \phi^{-1}(y)$ and $\phi^{-1}(y) - \ker (\phi) \subseteq \phi^{-1}(y)$.
\end{corollary}
Note that if $\zD < \zO$ (in $N$), neither $\ker(\phi)$ nor $\ker(\phi) \cup - \ker(\phi)$ are, in general, sub\doms of $M$.

See \S\ref{SEC:MORPHISM-QUOTIENT} for more on \domhoms.
\subsection{Type of \doms}
\label{SUBSEC:TYPE}
\begin{definizione}
%Define $\zD := -\zO$.
$M$ is of the \intro{first type} if $\zD = \zO$.
It is of the \intro{second type} if $\zD + \zD < \zD$.
It is of the \intro{third type}  if $\zD + \zD = \zD < \zO$.
\end{definizione}
\begin{remark}
The classification of \doms into first, second and third type is a partition of the class of \doms (namely, every \dom has exactly one type by \ref{PROP:DOM}(\ref{EN:DOM-18})).

Moreover, $M$ is of the same type as its dual, and if $N$ is a \subdom of~$M$, then $M$ and $N$ are of the same type.
\end{remark}
\begin{enumexample}
%Let $\vg$ be an ordered group.
The \doms $\vg$ and $\vgt$ (see Example~\ref{EX:G-TILDE}) are of the first type.
Moreover, $\vgh$~is of the second type iff
%$\vg^{>0} := \set{x\in\vg: x > 0}$ has a minimum,%
$\vg$ has a minimal positive element,%
\footnote{This is equivalent to $\vg$ being discrete and non-trivial.}
otherwise it is of the third type.
\end{enumexample}
\begin{lemma}\label{LEM:FIRST-TYPE}
Assume that $M$ is of the first type.
%$\zO = \zD$.
Then $x - \zD = x = x +\zD$.
Moreover, if $\ig x = \zO$, then $x + (-x) = \zO$.
\end{lemma}
\begin{proof}
$x - \zD = x -\zO = x$.
Moreover, $x + \zD = x + \zO = x$.\\
Finally, if $\ig x = \zO$, then $x + (-x) = - \ig x = \zD = \zO$.
\end{proof}
\begin{lemma}\label{LEM:SECOND-TYPE}
Assume that $M$ is of the second type, and let $x \in M$ such that $\ig x = 0$.
%$\zO \ar \zO > \zO$.
\begin{enumerate}
\item
%If $\ig x = \zO$, then 
$x - \zD > x > x + \zD$;
\item moreover, $x + (-x - \zD) = \zO$.
\item For any $y \in M$, $(y - \zD) + \zD = y = (y + \zD) - \zD$.
\end{enumerate}
\end{lemma}
\begin{proof}
%The hypothesis is equivalent to $2\zD < \zD$.
\begin{enumerate}
\item If, for contradiction, $x - \zD \leq x$, then $x - (x - \zD) \geq \zO$, \ie $x \ar (-x + \zD) \geq \zO$.
By Proposition~\ref{PROP:DOM}(\ref{EN:DOM-11}), the \lhs is less or equal than $(x \dl x + \zD) -\zD$, therefore
$(-\ig x + \zD) - \zD \geq \zO$, \ie $2\zD -\zD \geq \zO$.
Thus, $2\zD \geq \zD$, hence $2\zD = \zD$, absurd.

The second inequality is the dual of the first one.
\item 
$x + (-x - \zD) \leq - \ig x - \zD = \ig \zD = \zO$.
If, for contradiction, $x + (-x - \zD) < \zO$, then 
$x + (-x - \zD) \leq \zD$, \ie $x \leq \zD - (- x -\zD) = x + \zD$, contradicting the previous point.
\item
Let us prove the first equality.
If $\ig y > \zO$, it is obvious, since  in that case $y + \zD = y - \zD = y$.
Otherwise, by the previous inequalities, $y - \zD > y \geq  (y - \zD) +\zD$.
Moreover, there can be only one element in the interval $[(y - \zD) + \zD, y - \zD)$, hence $(y - \zD) +\zD  = y$. 
The other equality is the dual.
\qedhere
\end{enumerate}
\end{proof}
Therefore, we have proved that, in the case when $M$ is of the first or second type,
%either $\zO = \zD$ or $\zD + \zD < \zD$,
$\Mgroup = \set{x\in M: \ig x = \zO}$ is an ordered group.
Note that in the first type case the group inverse of $x$ (with $\ig x = \zO$) is precisely $-x$.
In the second type case, instead, the inverse of $x$ is not $-x$, but $-x - \zD$, which is strictly greater than $-x$.
Moreover, if $M$ is of the second type, then $\Mgroup$ is discrete, with minimal positive element $\zO \ar \zO$.
%In these two cases, we will say that $M$ is a \intro{\quasigroup}.
%%
\subsection{Associated group and multiplicity}
\label{SUBSEC:ASS-GROUP}
\begin{definizione}
Define $\FP: M \to M$ as $\FP(x) := (x + \zD) - \zD$.
Define moreover $x \Mequiv y$ if $\FP(x) = \FP(y)$, and $\eqclass{x}$ to be the equivalence class of~$x$.

Define also $\FM : M \to M$ as $\FM(x) := (x - \zD) + \zD$.
\end{definizione}
Note that $\FM$ is the dual operation of $\FP$, namely $\FM(x) = -\FP(-x)$.
Hence, all theorems about $\FP$ have a dual theorem about $\FM$, and we will usually prove (and often also state) only one of the two forms. 
\begin{lemma}\label{LEM:F}
\begin{enumerate}
%\item\label{EN:F-DUAL} $\FM$ is the dual operation of $\FP$, namely $\FM(x) = -\FP(-x)$;
\item\label{EN:F-1} $x \leq \FP(x) \leq x - \zD$;
therefore, $\FM(x) \leq \FP(x)$;
\item\label{EN:F-1P} if $x \leq y$, then $\FP(x) \leq \FP(y)$ and $\FM(x) \leq \FM(y)$;
\item\label{EN:F-2} if $\ig x > \zO$, then $\FP(x) = x$;
\item\label{EN:F-3} $\ig{\FP(x)} = \ig x$;
\item\label{EN:F-4} the interval $\bigl( \FM(x) , \FP(x) \bigr)$ is empty;
\item\label{EN:F-5} $\FP\Pa{\FP(x)} = \FP(x)$;
\item\label{EN:F-11} if $M$ is of the first or second type, then $\FP$ is the identity;
\item\label{EN:F-5P} $\FP\Pa{\FM(x)} = \FP(x)$;
\item\label{EN:F-6} $\Mequiv$ is an equivalence relation;
\item\label{EN:F-7} if $\ig x > \zO$, then $\eqclass{x} = \mset{x}$;
\item\label{EN:F-8}
%there are at most $2$ elements $y$ and $y'$ in $\eqclass{x}$;
%moreover, $y = \FP(y')$ \rom(or $y' = \FP(y)$\rom);
the only elements of $\eqclass{x}$ are $\FM(x)$ and $\FP(x)$ (which might coincide); in particular, all elements of $M$ are either of the form $\FP(y)$ or of the form $\FM(y)$;
\item\label{EN:F-8P}
the equivalence relation induced by $\FP$ is the same as the one induced by~$\FM$; namely, $\FP(y) = \FP(y')$ iff $\FM(y) = \FM(y')$;
\item\label{EN:F-8PP}
$x \Mequiv y$ iff $x - \zD = y - \zD$.
\setcounter{saveenum}{\value{enumi}}
\end{enumerate}
\vspace{1ex}
In general, if $x \Mequiv x'$ and $y \Mequiv y'$, then
\begin{enumerate}
\setcounter{enumi}{\value{saveenum}}
\item\label{EN:F-9} $x + y \Mequiv x' + y'$;
\item\label{EN:F-10} $- x \Mequiv - x'$.
%\item\label{EN:F-12} $-x - \zD \equiv -x' - \zD$.
\end{enumerate}
\end{lemma}
\begin{proof}
\ref{EN:F-1}, \ref{EN:F-1P}, \ref{EN:F-2}, \ref{EN:F-3} and \ref{EN:F-6} are obvious.
\begin{prooflist}
\item[\ref{EN:F-4}] Assume for contradiction that $\FM(x) < y < \FP(x)$.
By~\ref{EN:F-1} and Proposition~\ref{PROP:DOM}(\ref{EN:DOM-12PP}), $y = x$.
Therefore, by Proposition~\ref{PROP:DOM}(\ref{EN:DOM-13PP}), $\FM(x) \ar x < x + \FP(x)$, \ie\ $\FM(x) \ar \ig x < \FP(x)$.
Thus, by~\ref{EN:F-1}, $(x + \zD) \ar \ig x < (x + \zD) - \zD$, hence $ \ig x < -\zD$, absurd.
\item[\ref{EN:F-5}] Immediate from Proposition~\ref{PROP:DOM}(\ref{EN:DOM-4PP}).
%Since $(y - \zD) + \zD \leq y$, $F(F(x)) = ((((x + \zD) - \zD) + \zD) - \zD) \leq (x + \zD) -\zD = F(x)$.
\item[\ref{EN:F-11}] Trivial from Lemmata~\ref{LEM:SECOND-TYPE} and~\ref{LEM:FIRST-TYPE}.
\item[\ref{EN:F-5P}] If $M$ is of the first or second type, the conclusion follows immediately from~\ref{EN:F-11}.
Otherwise,
\[
\FP\Pa{ \FM(x) } = \Pa{(x - \zD) + 2 \zD} - \zD = \Pa{(x - \zD) + \zD} - \zD =
x - \zD \geq \FP(x).
\]
%If, for contradiction, $\FP\Pa{ \FM(x)} < \FP(x)$, then
%$\Pa{(x - \zD) + 2 \zD} - \zD < (x + \zD) - \zD$.
%Therefore, $(x - \zD) + 2 \zD < x + \zD$, namely $\Pa{(x - \zD) + 2 \zD} - \zD < x$.
\item[\ref{EN:F-7}] Immediate from~\ref{EN:F-2} and~\ref{EN:F-3}.
\item[\ref{EN:F-8}] From~\ref{EN:F-5} and~\ref{EN:F-5P} we see that $\FM(x)$ and $\FP(x)$ are both in $\eqclass{x}$.
Let $y, y' \in\eqclass{x}$.
If $y' > y$, then $\FP(y) = \FP(y') \geq  y' > y$, hence, by~\ref{EN:F-4}, $y' = \FP(y)$.
Therefore, $\FM(y') = \FM\Pa{\FP(y)} = \FM(y)$, and as before we can conclude that $y = \FM(y')$.
\item[\ref{EN:F-8P}] Immediate from~\ref{EN:F-8}.
\item[\ref{EN:F-8PP}]
If $\ig x > \zO$, then $x - \zD = x$, and therefore $x - \zD = y - \zD$ iff $x = y$.
Thus, we can conclude using~\ref{EN:F-2}.
Suppose instead that $\ig x = \zO$.
If $M$ is of the first or second type, then $\Mgroup$ is a group.
Therefore, $x -\zD = y - \zD$ iff $x = y$, and the conclusion follows from~\ref{EN:F-11}.
Finally, if $M$ is of the third type, then 
\[
\FP(x) \leq x - \zD \leq \FP(x) - \zD = (x + \zD) - 2 \zD = (x + \zD) - \zD = \FP(x),
\]
hence $\FP(x) = x - \zD$, and we are done.
\item[\ref{EN:F-9}] It suffices to prove the case when $y' = y$.
Moreover, if $x' \neq x$, by~\ref{EN:F-8}, we can assume that $x' = \FP(x)$.
Since $x' \geq x$, $\FP(x' + y) \geq \FP(x + y)$.
Then
\begin{multline*}
%F(x + y)  &= (x + y + \zD) - \zD, \\
\FP(x' + y) = (((x + \zD) \overbrace{ - \zD) + y}^{switch} + \zD) -\zD \leq \\ 
\leq (((x + y + \zD) -\zD) + \zD) -\zD = \FP(\FP(x + y)) = \FP(x + y).
\end{multline*}
\item[\ref{EN:F-10}] %There are $2$ cases.
If $M$ is of the first or second type, then the conclusion is immediate from ~\ref{EN:F-11}.
%$F$ is the identity function, hence $\eqclass{x} = \set{x}$, and the conclusion is trivial.

%Otherwise, $2\zD = \zD < \zO$.
Otherwise, if $M$ is of the third type, then, as before, we can assume that $x' = \FP(x) > x$.
Hence, $\FP( - x') \leq \FP( - x)$.
Let $z := - x$.
Then, 
\begin{multline*}
%F(-x)  &= F(z) = (z + \zD) - \zD, \\
\FP(-x') = ((z - \zD) + 2\zD) - \zD = ((z - \zD) + \zD) - \zD =\\
= \FP( z - \zD) \geq \FP(z) = \FP(-x).
\qedhere
\end{multline*}
\end{prooflist}
\end{proof}
\begin{lemma}\label{LEM:QUOTIENT}
Let $\sim$ be an equivalence relation on $M$ such that for every $x$, $x'$, $y\in M$ such that $x \sim x'$:
\begin{enumerate}
\item $-x \sim - x'$;
\item $x + y \sim x' + y$;
\item if $x \leq y \leq x'$, then $y \sim x$.
\end{enumerate}
Then, $M/{\sim}$ inherits an \domstruct from $M$.
Moreover, $M/{\sim}$ is a \dom, and the quotient map $M \to M/{\sim}$ is a \domhom.
\end{lemma}
\begin{proof}
Trivial checks.
%considering that Axiom~\ref{AX:M5} is equivalent to Proposition~\ref{PROP:DOM}(\ref{EN:DOM-2}).
\end{proof}
By Lemmata~\ref{LEM:F} and~\ref{LEM:QUOTIENT}, the \domstruct on $M$ induces a well-defined \domstruct on the quotient $\Mquot$, and the quotient map $\pi$ sending $x$ to $\eqclass{x}$ is a \domhom.%
\footnote{See Lemma~\ref{LEM:DOMHOM} and Remark~\ref{REM:QUOTIENT} for more on \domhoms and equivalence relations.}
%In particular, when such homomorphisms are injective, the fact that when not injective the image is of the first type.
%%
\begin{lemma}
If $M$ is of the first or third type, then $\Mquot$ is of the first type.
If $M$ is of the second type, then $\Mquot$ is also of the second type.
%$\Mquot$ is a \quasigroup.
Therefore, the set $\G{M} := \zerogroup{(\Mquot)}$ is an ordered group.
\end{lemma}
\begin{proof}
If $M$ is of the first or second type, then the quotient map is the identity, and the conclusion is trivial.
Otherwise, $2\zD = \zD < \zO$ (in $M$).
We have to prove that $-\eqclass\zO = \eqclass\zO$, \ie\ that $\FP(\zD) = \zO$.
In fact, $\FP(\zD) = 2\zD -\zD = \zD - \zD = \zO$.
\end{proof}
\begin{definizione}\label{DEF:MULTIPLICITY}
The ordered group $\G{M}$ defined above is the \intro{group associated} to~$M$.
The cardinality of $\eqclass{x}$ (as a subset of~$M$)
%the equivalence class of $x$,
is the \intro{multiplicity} of~$x$.
For every $x \in \Mgroup$, $x$~is a \intro{simple point} iff it has multiplicity $1$, otherwise it is a \intro{double point}.
\end{definizione}
Note that the multiplicity of $x$ is either $1$ or $2$, and it can be $2$ only if $M$ is of the third type.
Moreover, if $M$ is of the second type, then $\G{M}$ is discrete.
In \S\ref{SEC:EMBEDDING} we will show that under some additional conditions, there is an important correlation between $M$, $\GT{M}$ (see Example~\ref{EX:G-TILDE}) and $\wideDed{\GM}$.
\subsection{Signature}\label{SEC:SIGNATURE}
In this section, $M$ will be a \dom of the third type, unless we explicitly say otherwise.
\begin{deflemma}\label{LEM:SIG}
%Let $x\in M$, such that $\ig x = \zO$.
Let $x\in \Mgroup$.
There are exactly $3$ distinct possible cases:
\begin{enumerate}
\item $x + \zD = x < x - \zD$;
\item $x + \zD = x = x - \zD$;
\item $x + \zD < x = x - \zD$.
\end{enumerate}
In the first case, we say that the signature of $x$ is $\sig x := -1$, in the second case the signature is $0$, in the third $+1$.
Moreover,
\begin{enumerate}
\item $\sig x = 0$ iff $x$ is a simple point;%
\footnote{See Definition~\ref{DEF:MULTIPLICITY}.}
\item $\sig x = -1$ iff $x$ is a double point
and $x < \FP(x)$, iff $x - \zD > x$;
\item $\sig  x = 1$ iff $x$ is a double point
and $x = \FP(x)$, iff $x + \zD < x$;
\item $\sig x \geq 0$ iff $x - \zD = x$;
\item $\sig x \leq 0$ iff $x + \zD = x$;
\item $\sig x = 0$ iff $x + \zD = x - \zD$, and $\sig x \neq 0$ iff $x + \zD < x - \zD$;
\item $\sig \zO = 1$, $\sig \zD = -1$;
\item $\sig{\FP(x)} \geq 0$ and $\sig{\FM(x)} \leq 0$.
\end{enumerate}
\end{deflemma}
\begin{proof}
Easy.
\end{proof}
\begin{enumexample}
If $M = \vgh$ for some dense group $\vg$, and $\Lambda \in \Mgroup$, then $\sig\Lambda = 1$ iff $\Lambda = \lambda^+$ (for some $\lambda\in\vg$), $\sig \Lambda = -1$ iff $\Lambda = \lambda^-$, $\sig \Lambda = 0$ otherwise.
\end{enumexample}
\begin{proposition}[Signature rule]\label{PROP:SIG-RULE}
Let $x, y\in \Mgroup$.
Then, $\sig{-x} = - \sig x$.
Moreover, define $\alpha := \sig x$, $\beta := \sig y$, and $\gamma := \sig{x + y}$. Then:
\begin{enumerate}
\item\label{EN:SIGN:++} if $\alpha = \beta = 1$, then $\gamma = 1$;
\item\label{EN:SIGN:--} if $\alpha = \beta = -1$, then $\gamma = -1$;
\item\label{EN:SIGN:0?} if $\alpha \leq 0$, then $\gamma \leq 0$;
\item\label{EN:SIGN:+?} if $\alpha = 1$ and $\beta \geq 0$, then $\gamma \geq 0$;
\item\label{EN:SIGN:+0} if $\alpha = 1$ and $\beta = 0$, then $\gamma = 0$;
\item\label{EN:SIGN:-0} if $\alpha = -1$ and $\beta = 0$, then $\gamma = 0$;
\item\label{EN:SIGN:00} if $\alpha = \beta = 0$, then $\gamma \leq 0$;
\footnote{Both $\gamma = 0$ and $\gamma = -1$ are possible.
For instance, take $\vg = \Zedtwo$: if $x = y = \unmezzo$, then $x + y = 1^-$;
if $x = y = \unquarto$, then $x + y = \unmezzo$.}
\item\label{EN:SIGN:+-} if $\alpha = 1 $ and $\beta = -1$, then $\gamma = -1$.%
\end{enumerate}
\end{proposition}
\begin{proof}
\begin{proofenum}
\item Since $x + \zD < x$ and $y + \zD < y$, then, by Proposition~\ref{PROP:DOM}(\ref{EN:DOM-13PP}), $x + \zD + y + \zD < x + y$, \ie\ $x + y + \zD < x + y$, and the conclusion follows.
\item Since $x - \zD > x$ and $y - \zD > y$, then, by Proposition~\ref{PROP:DOM}(\ref{EN:DOM-13PP}), $x - \zD \ar y - \zD < x + y$, thus $(x + y) - \zD < x + y$, and the conclusion follows.
\item Since $x + \zD = x$, then $(x + y) + \zD = (x + \zD) + y = x + y$.
\item
\makeatletter
\refstepcounter{i@claim}
%\begin{claim}
\label{CLAIM:SIGN1}
\textit{Claim \thei@claim.}
\makeatother
$(x + \zD) \ar y = x + y$.
\par
%\end{claim}
Since $x + \zD < x$, them, by Proposition~\ref{PROP:DOM}(\ref{EN:DOM-10PP}), $(x + \zD) \ar y \leq  x + y$.
Moreover, by Proposition~\ref{PROP:DOM}(\ref{EN:DOM-9}), $(x + \zD) \ar y \geq x + (y \ar \zD) = x + y$.

Therefore, $(x + y) -\zD = (x + \zD) \ar y - \zD$, and, since $\beta \geq 0$, the latter is equal to $(x + \zD) \ar y = x + y$.
\item %If $\alpha = 1$, the conclusion follows 
Immediate from~\ref{EN:SIGN:0?} and~\ref{EN:SIGN:+?}.
\item \ref{EN:SIGN:0?} implies that $\gamma \leq 0$.
Moreover, the dual of %we can deduce, as in 
Claim~\ref{CLAIM:SIGN1} implies that $(x - \zD) + y = x \ar y$.
Suppose, for contradiction, that $(x + y) - \zD > x + y$.
Therefore, $x + y < x \ar y - \zD = x \ar y = (x -\zD) + y$.
Let $z := x - \zD$.
We have $z \Mequiv x$, hence $x + y \Mequiv z + y$, thus $\sig z = 1$ and $\sig{z + y} = 1$, contradicting~\ref{EN:SIGN:+0}.
\item Immediate from~\ref{EN:SIGN:0?}.
\item By~\ref{EN:SIGN:0?}, $\gamma \leq 0$.
By hypothesis, $x + \zD < x$ and $y < y -\zD$, hence $x + y + \zD< x + (y - \zD)$.
The latter is less or equal to $(x + y) -\zD$, and we deduce that $\gamma \neq 0$.
\qedhere
\end{proofenum}
\end{proof}
We will now define the signature of $x$ in the case when $M$ is not of the third type, or $x \notin \Mgroup$.
\begin{definizione}
Let $M$ be a \dom, $y \in M$, and $x \in \Mgroup$.
If $M$ is of the second type, then $\sig x := \infty$.
If $M$ is of the first type, then $\sig x := \firstsign$.
Note that $y \in \gqdom{M}{\ig y}$.
Define $\sig y$ as the signature of $y$ in $\gqdom{M}{\ig y}$.
\end{definizione}
Note that if $y \notin \Mgroup$, then $\sig y \neq \firstsign$.
Equivalently, if $\ig y > \zO$, then $\gqdom{M}{\ig y}$ cannot be of
the first type.
Moreover, if $\ig y > \ig x$, then $\sig{x + y} = \sig{y}$.

The above definition of signature coincides with Tressl's definition in~\cite{TRESSL:2005}, except that we introduced the $\firstsign$ symbol for \doms of the first type, which he does not treat in his article.
\begin{definizione}
Let $M$ be a \dom of the third type.
Let $\twoM$ be the set of $x\in \Mgroup$ with multiplicity~$2$,
and $\HM$ be the image of $\twoM$ in $\GM$ via the quotient map $\pi$.
\end{definizione}
\begin{lemma}\label{LEM:HM}
$\twoM$ is a \subdom of $M$, and $\HM$ is a subgroup of $\GM$.
Moreover, for every $x\in\twoM$ and $y \in M$, $x + y \in \twoM$ iff $y\in\twoM$ iff $x \ar y \in \twoM$.
\end{lemma}
\begin{proof}
Immediate consequence of Lemma~\ref{LEM:SIG} and Proposition~\ref{PROP:SIG-RULE}.
\end{proof}
\subsection{Proper and trivial \doms}
\begin{definizione}
$M$ is a \intro{proper \dom{}} if for every $x < y \in M$ there exists $z \in \Mgroup$  such that $x \leq z \leq y$.

It is \intro{strongly proper} if it is proper, and for every $y \in \notmzero$ there exists $x \in \Mgroup$ such that $\zO < x < \ig y$.
%\todo{Give some equivalent definition of strongly proper.}
\end{definizione}
\begin{remark}
A \dom $M$ is strongly proper iff for every $x < y \in M$ such that $\ig y > \zO$ there exists $z \in \Mgroup$ such that $x < z < y$.
\end{remark}
\begin{lemma}\label{LEM:AR}
Let $x$, $y\in M$.
Then, $x \ar y \leq (x \ar \ig x) + (y \ar \ig y)$.
\end{lemma}
\begin{proof}
By Proposition~\ref{PROP:DOM}, $x \ar y \leq (x + y) \ar \ig x$.
If $\ig x \neq \ig y$, then $x \ar y = x + y$, and we are done.
Otherwise, define $M' := \gqdom{M}{\ig x}$: 
$x$ and $y$ are in $M'$, and $\ig x = \ig y =\zO_{M'}$.
By splitting into cases according to the type of $M'$ we can easily prove the lemma.%
\footnote{Problem: give an easy proof of the lemma without distinguishing the various cases.}
\end{proof}
\begin{definizione}\label{DEF:LAMBDA}
%Let $\vg := \GM$.
Let $\vl$ be an ordered group containing $\vg:=\GM$.

For every $a \in \notmzero$, define $\Lk(a) \in \vlh$ by 
\[
\Lk(a) := \set{x \in \vl: \exists y \in \Mgroup\ x \leq [y] \et y < a}^+.
\]
%\[
%\Lambda(a) := \cut{\set{[x] \in \GM: x < a}}{\set{[x] \in \GM: x > a}}.
%\]
Moreover, we will write $\Lambda$ for $\Lambda_{\vg}$.
\end{definizione}
\begin{lemma}\label{LEM:LAMBDA}
Let $M$ be a proper \dom, $\vg := \GM$, and $\vl$ be an ordered group containing $\vg$.
Assume that for every $a \in \notmzero$ there is no $x\in\vl$ such that, for all $y,y'\in\Mgroup$, $y<a<y$ implies $[y] < x < [y']$.
%\[
%\cut{\set{x \in \vl: \exists y \in \Mgroup\ x \leq [y] \et y < a}}
%{\set{x' \in \vl: \exists y' \in \Mgroup\ x' \geq [y'] \et y' > a}}
%\]
%is a cut of $\vl$.%
%\footnote{Namely, $\sup\set{x \in \vl: \exists y \in \Mgroup\ x \leq [y] \et y < a} = \inf\set{x' \in \vl: \exists y' \in \Mgroup\ x' \geq [y'] \et y' > a}$.
%Note that such cut must be equal to $\Lk(a)$.}
Then, the map $\Lk: \notmzero \to \vlh$ is an injective \qdomhom.
\end{lemma}
Note that if $M$ is proper, then $\vl = \GM$ satisfies the hypothesis of the lemma.
\begin{proof}
$\Lk$ preserves the minus, because
\begin{multline*}
\Lk(-a) = \set{x \in \vl: \exists y \in \Mgroup\ x \leq [y] \et y < -a}^+=\\
=\set{-x' \in \vl: \exists y' \in \Mgroup\ x' \geq [y'] \et y' > a}^+ =\\
=-\set{x' \in \vl: \exists y' \in \Mgroup\ x' \geq [y'] \et y' > a}^- =
-\Lk(a).
\end{multline*}
$\Lk$ is injective and preserves the order, because if $a < a' \in \notmzero$, then there exists $y_0 \in \Mgroup$ such that $a < y_0 < a'$.
Moreover,
\[
\set{x \in \vl: \exists y \in \Mgroup\ x \leq [y] \et y < a} \subset
\set{x \in \vl: \exists y \in \Mgroup\ x \leq [y] \et y < a'},
\]
and $[y_0]$ is not in the former set.
Hence, $\Lk(a) < [y_0] < \Lk(a')$.

Finally, $\Lk$ preserves the sum.
%To simplify the proof, we will assume that $M$ is of the third type (the lemma is true in general, and the proof is a modification of the present one, along the lines of the proof of Lemma~\ref{LEM:LAMBDA}).
Let $a$ and $a'\in A$.
Then,
\begin{multline*}
\Lk(a) + \Lk(a') =
\set{x + x': \exists y, y' \in \Mgroup\ x \leq [y] \et y < a \et x' \leq [y'] \et y' < a'}^+ \leq\\
\leq \set{x'' \in \vl: \exists y'' \in \Mgroup\ x'' \leq [y''] \et y'' < a + a'}^+ =
\Lk(a + a').
\end{multline*}
Assume, for contradiction, that $\Lk(a) + \Lk(a') < x_0 < \Lk(a + a')$ for some
$x_0 \in \vl$.
Then, for every $x$, $x'\in \vl$ such that there exist $y$, $y' \in \Mgroup$ satisfying $x \leq [y] < [a]$ and $x' \leq [y'] < [a']$, we have $x + x' < x_0$.
Moreover, $x_0 < \Lk(a + a')$, \ie\ there exists $y'' \in \Mgroup$ such that $x_0 \leq [y''] < [a + a']$.

Hence, for every $y$, $y' \in \Mgroup$ such that $y < a$ and $y' < a'$, we have $y + y' < y'' < a + a'$.
Since $y'' < a + a'$, we have $y'' - a < a'$.
Moreover, both $a'$ and $y'' - a$ are in $\notmzero$.
Hence, there exists $z' \in \Mgroup$ such that $y'' - a < z' < a'$.
Therefore, $y'' - z' < a$.
Define $z := y'' - z'$, $y := z \ar \zO$, and $y' := z' \ar \zO$.
Since $z$ and $z$ are in $\Mgroup$, $z < a$ and $z' < a'$, we have
also $y$, $y'\in\Mgroup$, $y < a$ and $y' < a'$, therefore $y + y' < y''$.
By Lemma~\ref{LEM:AR}, $z \ar z' \leq y + y'$, hence $z \ar z' < y''$.
However, $z \ar z' = y'' - z' \ar z' = y'' \ar \zO \geq y''$, absurd.
\end{proof}
Note that $\Lambda$ is not a \domhom, because either it does not preserve the neutral element, or $\notmzero$ has no neutral element.
\begin{question}
What happens if $M$ is not proper? Does $\Lambda$ still preserve the sum?
\end{question}
\begin{remark}\label{REM:LAMBDA-COFINAL}
Let $M$ be a strongly proper \dom
% such that $\vg := \GM$ is a non-trivial group,%
%\footnote{Namely, $M$ has more than $4$ elements.}
and $\vl$ be an ordered group containing $\vg:=\GM$.
Assume moreover the following hypothesis:
\begin{itemize}
\item[\rom(*\rom)]
For every $x \in \vl \setminus \vg$ and $\varepsilon > 0$ in~$\vg$, there exist $\gamma$ and $\gamma' \in \vg$ such that $\gamma < x < \gamma'$ and $\gamma' - \gamma \leq \varepsilon$.
\end{itemize}
Then, $M$ and $\vl$ satisfy the hypothesis of Lemma~\ref{LEM:LAMBDA}.%
\footnote{\textbf{Question:} Does there exist a proper \dom $M$ with a non-trivial associated group $\vg$ and an super-group $\vl$ of $\vg$, such that $M$ and $\vl$ satisfy the hypothesis of Lemma~\ref{LEM:LAMBDA}, but not~(*)?}
\end{remark}
\begin{proof}
Assume for contradiction that there exist $a \in \notmzero$ and $x_0 \in \vl$ such that
\[
\set{x \in \vl: \exists y \in \Mgroup\ x \leq [y] \et y < a}^+ < x_0 <
\set{x' \in \vl: \exists y' \in \Mgroup\ x' \geq [y'] \et y' > a}^-.
\]
Let $y_0 \in \Mgroup$ such that $\zO < y_0 < \ig a$ (it exists because $M$ is strongly proper).
Since $x_0$ cannot be in $\vg$,  we can find $y$, $y'\in \Mgroup$ such that
$[y] < x_0 < [y']$ and $y' - y \leq y_0$.
Therefore, $y  < a < y'$, thus $y' - y \geq \ig a > y_0$, absurd.
\end{proof}
In \S\ref{SUBSEC:EMBEDDING-COLLAPSE} we will see more on proper \doms,
and some applications of the above results.
\begin{definizione}\label{DEF:DOM-TRIVIAL}
$M$ is a \intro{trivial \dom} if
\[
x + y = %\left\{ \begin{aligned}
\begin{cases}
x & \text{if } \abs x > \abs y,\\
y & \text{if } \abs x < \abs y,\\
\min(x,y) & \text{if } \abs x = \abs y .
%\end{aligned}\right.
\end{cases}\]
\end{definizione}
\begin{remark}
The only trivial \doms are the ones shown in examples~\ref{EX:TRIVIAL} and~\ref{EX:TRIVIAL-2}.
\end{remark}
\begin{remark}
$M$ is trivial iff $\forall x \in M$ $\ig x = \abs x$.
\end{remark}
\begin{remark}\label{REM:DOM-TRIVIAL-UNIQUE}
Let $\struct(M,\leq, \zO, -)$ be a structure, such that $\leq$ is a linear
ordering of $M$, $\zO$ is an element of $M$, and $-:M \to M$ is an
anti-automorphism of $(M,\leq)$, such that $-(-x)=x$.
Assume moreover that $M$ satisfies axioms~\ref{AX:M1} and~\ref{AX:M4}.
Then, there exists a unique binary operation $+$ on $M$ such that $\struct(M,\leq, \zO, +, -)$ is a \emph{trivial} \dom.
\end{remark}
\begin{proof}
$x + y$ is defined as in Definition~\ref{DEF:DOM-TRIVIAL}.
\end{proof}
\begin{enumexamples}
\item\label{EX:M-INFTY-TRIVIAL} If $M$ is a trivial \dom, then $M^\infty$ (defined in Example~\ref{EX:M-INFTY}) is also trivial.
\item\label{EX:FINITE-DOM} For every $n \in \Nats$ there exists exactly one (up to \dom-isomorphisms) \dom with $n$ elements, which we will denote by $\fdom{n}$.
Moreover, once we fix a linear ordering of the set of $n$ elements, there is only one \dom \emph{tout court} which subsumes the given order.
The existence and uniqueness is proved by induction on $n$, starting with the \doms with $1$ and $2$ elements,
and proceeding from $n$ to $n+2$ using Example~\ref{EX:M-INFTY} and Proposition~\ref{PROP:DOM}(\ref{EN:DOM-15}).
Moreover, by Example~\ref{EX:M-INFTY-TRIVIAL}, $\fdom{n}$ must be trivial.
\item\label{EX:PROPER-TRIVIAL}
The finite \doms $\fdom{1}$, $\fdom{2}$, $\fdom{3}$ and $\fdom{4}$ are the only \doms which are both trivial and proper.
\item
If $\vg$ is a densely ordered group, then $\vgt$ is a proper \dom of the first type, but not a strongly proper one.
For any ordered group~$\vg$, the \dom $\vg \times \fdom{4}$ (with lexicographic order and component-wise plus and minus: \cf\ \S\ref{SUBSEC:PRODUCT}) is a proper \dom of the third type, but not a strongly proper one.
On the other hand, every proper \dom of the second type is also strongly proper.
\end{enumexamples}
%%

%% subsection of proper and trivial doms
%%

%%
\section{Constructions on \doms}\label{SEC:CONSTRUCTIONS}
To understand the constructions in this section, the reader is advised to try them in the cases when $M$ is equal either to~$\vg$, or to~$\vgh$, or to~$\vgt$,
for some ordered group~$\vg$.
% a group $\vg$, or the insemination of the \dom $\vgh$ (defined in \S\ref{SUBSEC:INSEMINATION}).
%%
\subsection{Morphisms and quotients}\label{SEC:MORPHISM-QUOTIENT}
We will now study quotients and maps of \doms.
We have seen in Lemma~\ref{LEM:QUOTIENT} that under suitable conditions
an equivalence relation $\sim$ on a \dom $M$ induces a \dom structure on the quotient $M/\sim$.
We will call an equivalence relation $\sim$ on a \dom $M$ a \intro{\domeqrel}
if $\sim$ satisfies the hypothesis of Lemma~\ref{LEM:QUOTIENT}.

We will now prove the analogue for \doms of some basic theorems for groups.
\begin{lemma}\label{LEM:DOMHOM}
Let $\phi: M \to N$ be a \domhom.
Then, $\phi$ is injective iff $\ker(\phi) = \mset{\zO}$.
Moreover, if $\phi$ is not injective, then $\ker(\phi)$ is a \subdom of $M$, and $N$ is of the first type.
\end{lemma}
\begin{proof}
Easy.
\end{proof}
\begin{remark}\label{REM:QUOTIENT}
In particular, if $\sim$ is a non-trivial \domeqrel on $M$, then the quotient map is a \dom-homomorphism with non-trivial kernel, thus $M/\sim$ is a \dom of the first type.
\end{remark}
We will now show that every convex \subdom of $M$ defines an equivalence relation.
\begin{definizione}[Quotients]
Let $N$ be a convex \subdom of $M$.
%We will write $M/N$ instead of $M/\simN$.
Define an equivalence relation $\simN$ on $M$  in the following way:\\
$x \simN y$ if there exist $w_1, w_2 \in N$ such that $y + w_1 \leq x \leq x \ar w_2$.%
\footnote{Actually, we do not need $N$ to be convex to define $\simN$. However, the equivalence relation defined by $N$ is the same as the one defined by the convex closure of~$N$.}\\
We will write $M/N$ instead of $M/\simN$.
\end{definizione}
\begin{lemma}
Let $N$ be a convex \subdom of $M$.
The binary relation $\simN$ on $M$ is indeed a \domeqrel on $M$.
 Finally, the quotient map $\pi: M \to M/N$ is a surjective \dom-homomorphism, with kernel $N$.
\end{lemma}
\begin{proof}
Easy.
\end{proof}
\begin{lemma}\label{LEM:FIRST-DOM-HOMOMORPHISM}
%\todo{What is the name of the analogue theorem for groups? First homomorphism theorem?}
Let $P$ be a \dom of the first type, $\phi: M \to P\,$ be a \domhom,
and $N := \ker(\phi)$.
Then, the map $\phibar :M/N \to P$ sending the equivalence class $[x]_N$ to $\phi(x)$ is a well-defined injective \domhom.
\end{lemma}
\begin{proof}
If $x \simN y$, then $y + w_1 \leq x \leq y \ar w_2$
for some $w_1$, $w_2\in N$.
Therefore,
\[
\phi(y) = \phi(y + w_1) \leq \phi(x) \leq \phi(y \ar w_2) = \phi(y) \ar \zDM[P]  = \phi(y).
\]
Thus, $\phibar$ is well-defined.

If $\phi(x) = \phi(y)$, then $\phi(x - y) = \phi(y - x) = \zOL$.
Let $w_1 = x \dl y$, $w_2 = x - y$.
Hence, $y + w_1 = x \dl \ig y \leq x$, and $y \ar w_2 = x \ar \ig y \geq x$.
Moreover, $\phi(w_1) = \phi(w_2) = \zOL$, thus $x \simN y$.
Therefore, $\phibar$ is injective.

The fact that $\phibar$ is a \domhom is trivial.
\end{proof}
Therefore, we can identify the image of $\phi$ with $M/{\ker(\phi)}$.
Moreover, there is a bijection between convex \subdoms of $M$ and \domeqrels on $M$, given by the map $N \mapsto \simN$.

\begin{remark}
Let $M$ be a \dom of the third type, $N$ a \dom of the first type, and $\phi: M \to N$ a \domhom.
Then, $\phi$ factors uniquely through $\Mequiv$.
Namely, there exists a unique function $\phidot : \Mquot \to N$ such that $\phidot([x]) = \phi(x)$.
\end{remark}
\begin{proof}
In light of Lemma~\ref{LEM:FIRST-DOM-HOMOMORPHISM}, the conclusion is equivalent to the fact that every \domeqrel $\sim$ such that $M/\sim$ is of the first type is a coarsening of $\equiv$.
Namely, if $x \equiv y$, then $x \sim y$.
Since a \domeqrel is uniquely determined by the equivalence class of $\zO$, it suffices to treat the case when $x \equiv \zO$, that is $x = \zO$ or $x = \zD$.
Since $N$ is of the first type, in both cases $\phi(x) = \zO$.
\end{proof}
\subsection{Compatible families of morphisms}
In the following, when two or more \predoms $M$, $N$ are involved, we will sometimes need to distinguish the zero of $M$ from the one of $N$:
we will then use the notation $\zOM$ and $\zDM$ for the zero of $M$ and its opposite.

Moreover, we will sometimes need to split Axiom~\ref{AX:M5} into two parts:
\begin{enumerate}[label={\ref{AX:M5}(\alph*).}, ref={\ref{AX:M5}(\alph*)}]
\item\label{AX:M5-RIGHT} If $x < y$ then $x - y < \zO$.
\item\label{AX:M5-LEFT}  If $x - y < \zO$ then $x < y$.
\end{enumerate}
For, convenience, we will give a name to the following axiom of \predoms:
\begin{enumerate}[label={\textbf{P}\Alph*.}, ref={\textbf{P}\Alph*}]
\item\label{AX:PA} If $z < z'$, then $z + t \leq z' + t$.
\end{enumerate}
\begin{remark}\label{REM:MC-EQ}
For a \predom, Axiom~\ref{AX:M5-LEFT} is equivalent to $\forall x\ \ig x \geq \zO$.
\end{remark}
\begin{definizione}[Extensible group]
For every \dom $M$, we will denote by $\EG{M}$ the following ordered group:
\begin{itemize}
\item $\mset{0}$, if $M$ of the first type;
\item $\GM$, if $M$ is of the second type;
\item $\HM$, if $M$ is of the third type.
\end{itemize}
As a \dom, the minus of $\EG{M}$ is the group minus (which can be different from the one induced by $M$).
Note that the structure $\struct(\EG{M},\leq, 0, +)$ (\ie the additive structure of the group) has a natural embedding in $\struct(M,\leq,0,+)$, sending $x$ to $x^+$, the maximum of $x$ as a subset of~$M$.
\end{definizione}
\begin{notation}
Let $k \in \Orders{M}$.
For $x \in \gqdom{M}{k}$, define $[x]_k$ to be the equivalence class of $x$ in $\gqMk$.
\end{notation}
\begin{definizione}[Compatible families]\label{DEF:COMPATIBLE}
Let $M$ be an arbitrary \dom, $N$~be a \dom of the second or third type, and
$O$ be a final segment of $\Orders{N}$.
Let $\Theta$ be a sequence $\Pa{\theta_i}_{i \in O}$ of \domhoms \st each $\theta_i$ maps $M$ in a subgroup of $\EG{\gqdom{N}{i}}$.
We say that $\Theta$ is a \intro{compatible} family from $M$ to $N$ (or that the the $\theta_i$ are compatible), if for every $k < j \in O$
\begin{equation}\label{EQ:COMPATIBLE}
\theta_j(x) = [j + \theta_k(x)^+]_j.
\end{equation}
\end{definizione}
Note that if $O$ has a minimum $k$, and $\theta_k:M \to \EG{\gqNk}$ is a \domhom,
then there exists a unique compatible family $\Theta \Pa{\theta_i}_{i\in O}$ containing $\theta_k$: each $\theta_j$ is defined
via the formula~\eqref{EQ:COMPATIBLE}.
\begin{lemma}\label{LEM:RK}
Let $k \in \Orders{M}$.
Define the map $\Rk : M \to \gqMk/\Mequiv$ sending $y$ to $[y + k]_k$.
Then, $\Rk$ is a \domhom.
Moreover, the kernel of $\Rk$ is $[-k,k]$.
Besides, if $k > 0$, then $\Rk(M^{<k}) \subseteq \EG{\gqMk}$.
Finally, for every $\Lambda > -\infty$ cut of $\Orders{M}$, the family $\Pa{\Rk\upharpoonright{M^{\{\Lambda^L\}}}}_{k>\Lambda}$ is a compatible family from $M^{\{\Lambda^L\}}$ to $M^{\{\Lambda^R\}}$.
\end{lemma}
\begin{proof}
Easy.
\end{proof}
We will see in the rest of this section various applications of the above construction.
\subsection{Gluing}
\begin{definizione}[Gluing of \doms]
Let $M$, $N$, and $O$ be as in Definition~\ref{DEF:COMPATIBLE}, and $\Theta =\Pa{\theta_i}_{i\in O}$ be a compatible family of homomorphisms from $M$ to $N$.
We define $\cruxMTN$ as the set $M \sqcup \NO$, endowed with the following structure on the language $\struct(\leq, 0, +, -)$:
\begin{itemize}
\item the zero of $\cruxMTN$ is $\zOM$;
\item the order, the sum and the minus  extend the ones on $M$ and $N$;
\item if $x \in M$ and $y \in N$, then $x + y = y + x = \theta_{\ig y}(x)^+ \an  y \in N$, and $x < y$ iff \mbox{$\theta_{\ig y}(x)^+ \leq_N y$}.
\end{itemize}
If $O$ has a minimum $k$, we will write $\cruxMtkN$ for $\cruxMTN$, where $\Theta$ is the unique compatible family containing $\theta_k$.
\end{definizione}
Note that the hypothesis of $N$ being of the second or third type is necessary only when $O = \Orders{N}$.
\begin{lemma}\label{LEM:CRUX}
$\cruxMTN$ is a \dom of the same type as $M$.
Moreover, $M$ is a convex \subdom and $N$ a \qsubdom of $\cruxMTN$ via the inclusion maps.
Besides, $\cut{\Orders{M}}{O}$ is a cut of $\Orders{\cruxMTN}$.
Finally, $M = (\cruxMTN)^{\{< O\}}$ and $N^{\{O\}} = (\cruxMTN)^{\{O\}}$.
\end{lemma}
\begin{proof}
The correctness of our definition (\ie the well-definedness of $\leq$), as
well as the monoid properties, the involutiveness of $-$ and axioms~\MMA
and~\MMB are trivial.
%The anti-monotonicity of $-$ is an easy consequence of \MMC, therefore we only have to deal with the monotonicity of $+$ and the Axiom~\MMC.

The proof of monotonicity for $+$ is nothing but a boring enumeration of
cases: we will prove, as an example, that $a \leq b \rightarrow a + c \leq b + c$, when $a \in \MN$ and $b,c \in \MM$.
For $a \leq b$ implies, by definition, $a < \Mpthet{\Meord{a}}{b}$, and
$a + \Mpthet{\Meord{a}}{c} < \Mpthet{\Meord{a}}{b} + \Mpthet{\Meord{a}}{c}$ (since we have $\leq$ by the monotonicity of $+$ on $\MN$, and, if by contradiction
equality holds, then we will have $a = \Mplus{[a]_{\ig a}}$, but in that case $a < \Mpthet{\Meord{a}}{b}$ implies $[a]_{\ig a} < \theta_{\ig a}(b)$);
now $\Mpthet{\Meord{a}}{b} + \Mpthet{\Meord{a}}{c} = \Mpthet{\Meord{a}}{b + c}$,
and, substituting the right side, the previous inequality becomes
$a + \Mpthet{\Meord{a}}{c} < \Mpthet{\Meord{a}}{b+c}$, which is just the definition of $a + c \leq b + c$.

Let us prove that $-$ is anti-monotone: \ie\ that if $a < b$, then $-a > -b$.
We will treat only the case when $a \in M$ and $b \in N$ (the other cases are either trivial or similar to this one).
By definition, $\theta_{\ig b}(a)^+ \leq_N b$.
We must prove that $\theta_{\ig b}(-a)^+ >_N -b$.
Since $\theta_{\ig b}$ is a \domhom, $\theta_{\ig b}(-a) = -_G\theta_{\ig b}(a)$, where $-_G$ is the group minus of~$\G{\gqNb}$.
Moreover,
\[ \Pa{-_G\theta_{\ig b}(a)}^+ = 
\begin{cases}
-(\theta_{\ig b}(a)^-) &\text{if } \gqNb \text{ is of the third type},\\
-(\theta_{\ig b}(a)^+) + 1_{\ig b} &\text{if } \gqNb \text{ is of the second type},
\end{cases}\]
where $1_{\ig b}$ is the minimal positive element of $\gqNb$.
Hence, we must prove that $\theta_{\ig b}(a)^- <_N b$ if $\gqNb$ is of the third type, or that $(\theta_{\ig b}(a)^+) - 1_{\ig b} <_N b$ if $\gqNb$ is of the third second type, which are both trivial.

For Axiom~\MMC, we will prove the equivalence between $a < b$ and $a - b < 0$
assuming $a \in M$ and $b \in N$: the other cases are similar, or simpler.
By the anti-monotonicity of minus, $a < b$ is equivalent to $-b > -a$, \ie 
%$\Mpthet{\Meord{b}}{a} \leq b$, which in turn is equivalent to 
$-\Mpthet{\Meord{b}}{-a} <_N b$.
By Axiom~\MMC for $\MN$, we have $-\Mpthet{\Meord{b}}{-a} <_N b \leftrightarrow
-\Mpthet{\Meord{b}}{-a} - b <_N \Mpthet{\Meord{b}}{0}$, and the right side is just the expanded definition of $a - b < 0$.
\end{proof}
\begin{enumexample}
\label{EX:CRUX-UNION}
Let $\Lambda > -\infty$ be a cut of $\Orders{M}$.
Then, $M$ is canonically isomorphic to $\crux {M^{\{\Lambda^L\}}}{\Theta}{M^{\{\Lambda^R\}}}$, where
$\Theta := \Pa{\Rk\upharpoonright{M^{\{\Lambda^L\}}}}_{k>\Lambda}$.
\end{enumexample}
\subsection{Insemination}\label{SUBSEC:INSEMINATION}
\begin{definizione}[Insemination]
Let $M$ be a \dom of the third type, $P$~a be subgroup of~$\HM$,
and $\iota:P\hookrightarrow\HM$ be the immersion.
%Define $\Theta = \mset{\iota}$, where $\iota$ is the immersion of $P$ in $\HM$.
Define the \intro{insemination} of $M$ at $P$ as
\[
\insML := \crux{P}{\iota}{M}.
\]
\end{definizione}
\begin{remark}
$\insML$ is a \dom of the first type, whose underlying set is $M \sqcup P$.
Moreover, $M$ is a \qsubdom of $\insML$.
Besides, for every $y \in M$ and $[x] \in P$,
\[\begin{aligned}{}
[x] + \FP(y) &= \FP(x \am y),\\
[x] + \FM(y) &= \FM(x \am y),\\
[x] < y &\leftrightarrow \FP(x) \leq y,\\
[x]> y  &\leftrightarrow \FM(x) \geq y.
\end{aligned}\]
Finally, the map $\pibar: \insML \to \Mquot$ sending $[x]\in P$ to itself and $y \in M$ to $[y]$ is a \domhom with kernel $\set{\zDM, 0_G, \zOM}$.
\end{remark}
\begin{example}
If $\vg$ is a densely ordered group, then $\insem{\vgh}{\vg} = \vgt$.
\end{example}

\subsection{Union}\label{SEC:UNION}
\begin{definizione}[Union]\label{DEF:UNION}
Let $M$ be a \dom, $k$~be a \width of $M$ such that $k > \zO$, and $N$ be a \overdom of $\gqMk$.%
\footnote{Remember that this implies that the zero of $N$ is $k$.}
%For every $i \in \Orders{N}$, define the map $\theta_i: M \to \EG{\gqNi}$ sending $x$ to $[(x \am k) \an i]_i$.
%Let $\Theta = \Pa{\theta_i}_{k \leq i \in \Orders{N}}$.
Let $\theta_k: \lMk \to \EG{N}$ the map sending $x$ to $[x \am k]_k$.
Define 
\[
\unionMN := \crux{\lMk}{\theta_k} N.
\]
\end{definizione}
\begin{remark}\label{REM:UNION}
The universe of $\unionMN$ is $N \sqcup \lMk$ = $N \cup \lqMk$.
Moreover, for every $x \in N$ and $y \in \lMk$,
\[\begin{aligned}
x < y &\leftrightarrow x \leq_N y \dm k,\\
x > y &\leftrightarrow  x \geq_N y \am k,\\
x + y &:= x \an (y \am k).
\end{aligned}\]
Finally, the inclusion map from $\lqMk$ to $\unionMN$ is a \domhom, and the inclusion map from $N$ to $\unionMN$ is a \qdomhom.
\end{remark}
The following remark is a special case of Example~\ref{EX:CRUX-UNION}.
\begin{remark}
Let $M$ be a \dom, and $k > \zO$ be a \width of~$M$.
Then, $M = \union{M}{\gqMk}{k}$.
\end{remark}
%%%%%%
\subsection{Products}\label{SUBSEC:PRODUCT}
\begin{remark}\label{REM:FIRST-GROUP}
Let $M$ be a \dom of the first type, $y$, $y' \in M$, $x \in \Mgroup$.
Then, $x + y = x \ar y$, and $x + y \gtreqless x + y'$ iff $y \gtreqless y'$.
\end{remark}
\begin{proof}
$x \in \Mgroup$ means that $\ig x = \zO$. Thus,
\[
x + y \leq x \ar y \leq (x + y) \ar \ig x = (x + y) \ar \zO = x + y,
\]
because $\zO = \zD$.
Moreover, if $x + y = x + y'$, then $x \ar y - x = x \ar y' - x$, \ie\ $y \ar \zO = y' \ar \zO$, thus $y = y'$.
If instead $x + y < x + y'$, then obviously $y < y'$.
\end{proof}
\begin{definizione}[Fibered product of \doms]\label{DEF:FIBERED-PRODUCT}
Let $N$ be a \dom with a minimum $\muN$ \rom(and a maximum $\nuN := - \muN$\rom), $M$~be a \dom of the first type, and $\Agroup$ be a \subdom of~$\Mgroup$.
Define the \domstruct $M \fprod N$ in the following way.
\begin{itemize}
\item The universe of $M \fprod N$ is
$\Pa{ (M \setminus \Agroup) \times \mset{\muN}} \sqcup \Pa{\Agroup \times N}$.
It is a subset of $M \times N$.
\item The order is the lexicographic one, with $M$ more important than $N$:
namely, $(x,y) \leq (x',y')$ iff $x < x'$ or $x = x'$ and $y \leq y'$.
\item The zero is the pair $(\zOM, \zON)$.
\item The minus is defined component-wise, namely 
\[
-(x, y) := \begin{cases} (-x, -y)   & \text{if } x \in \Agroup;\\ 
                         (-x, \muN) & \text{otherwise}.
           \end{cases}
\]
\item $(x, y) + (x', y') := (x + x', y'')$, where
\[
y'' := \begin{cases}y + y' & \text{if } x + x' \in \Agroup;\\ 
                    \muN   & \text{if } x + x' \notin \Agroup.
\end{cases}\]
\end{itemize}
Define $M \fprod[] N := M \fprod[\Mgroup] N$.
If $\vg$ is a group, then $\vg \fprod[] N = \vg \fprod[\vg] N$.
\end{definizione}
Note that if $\vg$ is a group, then $\vg \fprod[] N$ is the Cartesian product of $\vg$ and $N$, with the lexicographic order and component-wise addition and subtraction.
If $N$ is also a group, then $\vg \fprod[] N$ coincides with the usual product of ordered groups (with lexicographic order).
\begin{lemma}\label{LEM:FIBERED-PRODUCT}
$M \fprod N$ is a \dom, of the same type as $N$.
Moreover, $(x, y) \ar (x', y') = (x \ar x', y''')$, where
\[
y''' := \begin{cases}y \ar y' & \text{if } x + x' \in \Agroup;\\ 
                    \muN   & \text{if } x + x' \notin \Agroup.
\end{cases}\]
%The natural map $M \to M \fprod N$ sending $x$ to $(x, \zON)$ is a homomorphism of ordered monoids, and if $\zON = \zDN$, it is actually a \dom-homomorphism.
Moreover, $N$ can be identified naturally with a convex \subdom of
$M \fprod N$ via the map sending $y$ to $(\zOM, y)$.
Finally, the map $\piuno : M \fprod N \to M$ sending $(x, y)$ to $x$ is a surjective \domhom, with kernel $N$.
\end{lemma}
\begin{proof}
The main fact that makes the above definition work is that $y + \muN = y$ for every $y \in N$.
Hence, if $x \notin \Agroup$ or $x' \notin \Agroup$, then $(x, y) + (x', y') = (x + x', \muN)$.
\begin{proofdescription}
\item[Associativity of sum]
$(x,y) + (x',y') + (x'',y'') = (x + x' + x'', y''')$, where
\[
y''' =
\begin{cases}
\muN & \text{if } x + x' + x'' \notin \Agroup;\\
y + y' + y'' & \text{if } x + x' + x'' \in \Agroup.
%\muN & \text{otherwise.}
\end{cases}
\]
\item[Axiom~\ref{AX:PA}] Let $z = (x, y)$, $z' = (x', y')$, $t = (r,s)$.
Since $z < z'$, then either $x < x'$, or $x = x' \in \Agroup$ and $y < y'$.
If $x + r < x' + r$, we are done.
Otherwise, $x + r = x' + r$.
If, for contradiction, $z + t > z' + t$, then $x + r \in \Agroup$.
%and $y + s > y' + s$.
Hence, $r \in \Mgroup$, thus Remark~\ref{REM:FIRST-GROUP}
implies that $x = x'$.
Therefore, $y + s > y' + s$, hence $y > y'$, absurd.
\item[Axiom~\ref{AX:M1}] $\zD = (\zDM, \zDN) \leq (\zOM, \zON) = \zO$.
\item[Axiom~\ref{AX:M4}] If, for contradiction, $\zD < (x, y) < \zO$, then $\zDM \leq x \leq \zOM$. However, $\zDM = \zOM$, hence $x = \zOM$.
Moreover, $\zDN < y < \zOM$, absurd.
\item[Axiom~\ref{AX:M5-LEFT}]
If $(x, y) - (x, y) < \zO$, then either $x - x < \zOM$ (impossible),
or $x - x = \zOM$ and $y - y < \zON$ (also impossible).
\item[Axiom~\ref{AX:M5-RIGHT}]
Assume, for contradiction, that $(x, y) < (x', y')$, but $(x, y) - (x', y') \geq \zO$.
Then, either $x < x'$, and we have a contradiction, or $x = x'$ and $y < y'$.
Since we have two different elements with the same abscissa $x$, we infer that $x \in \Agroup$, and therefore $\ig x = \zOM$.
Thus, $\zO \leq (x, y) - (x', y') = (\zOM, y - y')$.
Therefore, $\zON \leq y - y'$, a contradiction. \qedhere
\end{proofdescription}
\end{proof}
Note however that the map $\pidue: M \fprod N \to N$ sending $(x,y)$ to $y$ is a not a \domhom, because it does not preserve the sum.
On the other hand, the map from $M$ to $M \fprod N$ sending $x$ to $(x, \zON)$ if $x \in \Agroup$, or to $(x,\muN)$ otherwise, is a \domhom iff $N$ is of the first type.
%Note also that the \dom in Example~\ref{EX:RAZ-2} is nothing else than ${\Raz \fprod[\mset{0}] \fdom{2}}$.
Moreover, $M \fprod \fdom{1} = M$ for any choice of $\Agroup$.
Note also the following: if $x \in \Mgroup\setminus \Agroup$, then $\ig{(x,\zON)} = (\zOM, \nuN) > \zO$.
On the other hand, if $x \in \Agroup$, then $\ig{(x,\zON)} = \zO$.
\begin{example}
Let $A$ be a densely ordered group, and $B$ be an ordered group.
Then, $\Dedekind{A \times B}$ is naturally isomorphic to
$\Pa{\Dedekind A / \Mequiv} \fprod[\collapse A] \Dedekind B$.
\end{example}
\begin{definizione}\label{DEF:PRODUCT}
Let $M$ be a \dom of the first type, $N$ be any \dom, and $\mu$ be a symbol not in $N$.
Define the \domstruct $\MpN$ in the following way:
\begin{itemize}
\item The universe of $\MpN$ is $\Pa{\Mgroup \times N} \sqcup \Pa{\notmzero \times \mset {\mu}}$.
\item The order is the lexicographic one, with the first component  more important than the second one.
\item The zero is the pair $(\zOM, \zON)$.
\item The minus is defined component-wise, namely
\[-(x,y) = \begin{cases}
(-x, -y) & \text{if } x \in \Mgroup,\\
(-x, \mu) & \text{otherwise.}
\end{cases}\]
\item $(x,y) + (x',y') = (x + x', y'')$, where
\[y'' = \begin{cases}
y + y' & \text{if } x, x' \in \Mgroup,\\
\mu    & \text{otherwise.}
\end{cases}\]
\end{itemize}
\end{definizione}
Note that if $N$ has a minimum $\mu = \muN$ the definition above coincides with the one in~\ref{DEF:FIBERED-PRODUCT}.
\begin{lemma}\label{LEM:PRODUCT}
$\MpN$ is a \dom, of the same type as $N$.
Moreover, $(x, y) \ar (x', y') = (x \ar x', y''')$, where
\[
y''' := \begin{cases}
 y \ar y' & \text{if } x+ x' \in \Mgroup;\\ 
 \mu   & \text{if } x + x' \notin \Mgroup.
\end{cases}\]
Moreover, $N$ can be identified naturally with a convex \subdom of
$\MpN$ via the map sending $y$ to $(\zOM, y)$.
Finally, the map $\piuno : M \fprod N \to M$ sending $(x, y)$ to $x$ is a surjective \domhom, with kernel $N$.
\end{lemma}
\begin{proof}
The proof is a verbatim copy of the one of Lemma~\ref{LEM:FIBERED-PRODUCT}, letting $\Agroup := \Mgroup$ and $\muN = \mu$.
\end{proof}
Note that $\MpN$ contains a copy of $N$ for every point of $\Mgroup$.

We can combine the constructions in \S\S\ref{SUBSEC:INSEMINATION} and~\ref{SUBSEC:PRODUCT} in the following way.
Let $M$ be as \dom of the third type, and $N$ be any \dom.
Consider $T := \insem{M}{\HM} \fprod[] N$.
Then, $T$~is a \dom of the same type as~$N$.
It contains a copy of~$N$ inside every pair of double points of~$M$.
\begin{remark}
Let $M$, $\Agroup$, $N$ be as in Definition~\ref{DEF:FIBERED-PRODUCT}.
For every $k > 0 \in \Orders{M}$, define
\[\begin{aligned}
\theta_k &: \Mgroup \fprod N \to \EG{\gqMk}\\
         &  ([a],b) \mapsto \Rk(a) = [a + k]_k.
\end{aligned}\]
Then, $M \fprod N \simeq \crux{\Pa{\Mgroup \fprod N}}{\Theta}{N}$,
where $\Theta = \Pa{\theta_k}_{0 < k \in \Orders{M}}$.
\end{remark}
\begin{enumexample}
The \dom $M^{\infty}$ defined in~\ref{EX:M-INFTY} is equal to
$\fdom{3} \times M$.
\end{enumexample}
\subsection{Collapse}\label{SEC:COLLAPSE}
%%
%Another interesting example is the following one.
\begin{definizione}
Let $M$ be a \dom of the third type, and $P$ a subgroup of~$\GM$.
Hence, $\Mquot$ is a \dom of the first type, and %$\HM$ a subgroup of 
$\zerogroup{\Mquot} = \GM$.
Define $\collML$ (the \intro{collapse} of $M$ at~$P$) as $\prodMLdue$.
Define also $\pinL: \collML \to \Mquot$ as the projection onto the first component.
\end{definizione}
By Lemma~\ref{LEM:FIBERED-PRODUCT}, $\pinL$ is a surjective \domhom, with kernel $\fdom{2}$.
On the other hand, let $\eta: M \to \collML$ be the following map:
\[\eta(x) := \begin{cases}
([x],\zO) & \text{if } [x] \in P \text{ and }\sig x = 1,\\
([x],\zD) & \text{otherwise}.
\end{cases}\]
In general, $\eta$
is \emph{not} a \domhom, because the co-domain is a \dom of the third type,
and we have Remark~\ref{REM:QUOTIENT}.%
\footnote{For instance, take $M = \Razhat$, $P = \Zed$, $x ={\unmezzo}^-$, $x' = y ={\unmezzo}^+$.
Then, $x$ and $x'$ have the same image, but $x + y = 1^-$ and $x' + y = 1^+$  have different images.}
However, things are different when $P = \HM$.
\begin{lemma}
Let $M$ be a \dom of the third type and $H:= \HM$.
Then, $M$~is isomorphic to $\coll{M}{H}$, via the map $\eta$ defined above.
\end{lemma}
\begin{proof}
The only difficult part is showing that $\eta(z + z') = \eta(z) + \eta(z')$ for every $z$, $z' \in M$.
Assume not.
Since $[z] = \pi_H(\eta z)$, and $\pi_H$ is a \domhom, we must have $z + z' \in \twoM$, and $\sig{\eta(z + z')}$ different from $\sig{\eta z + \eta z'}$.
However, by definition of $\eta$, for every $x \in \twoM$, $\sig x = \sig{\eta x}$.
Hence, we must also have $z$ and $z'$ not in $\twoM$.
But then $z + z'$ and $\eta(z) + \eta(z')$ have both signature $-1$, and we are done.
\end{proof}
\subsection{Dedekind cuts}
We will now study $\mh$, the set of cuts of a \dom $M$.
As before, we will need $M$ to be of the first type.
The order and the minus are the obvious ones.
On the other hand, we have $4$ candidates for the plus, all of them making $\mh$ an ordered monoid.
However, only one of them makes $\mh$ a \dom.
\begin{definizione}
Given a \predom $M$, let $\mh$ be the set of Dedekind cuts of $M$.
Endow $\mh$ with a \domstruct, using the following rules: for every $\Lambda = \cut{\Lambda^L}{\Lambda^R}$ and $\Gamma = \cut{\Gamma^L}{\Gamma^R}$ in $\mh$
\begin{description}
\item[Order:] %$\cut{\Lambda^L}{\Lambda^R} \leq \cut{\Gamma^L}{\Gamma^R}$ 
$\Lambda \leq \Gamma$ if
$\Lambda^L \subseteq \Gamma^L$ iff $\Lambda^R \supseteq \Gamma^R$.
\item[Minus:] $-\Lambda := \cut{-\Lambda^R}{-\Lambda^L}$.
\item[Zero:] $\zO := \zOM^+ = \cut{(-\infty,\zOM]}{(\zOM,+\infty)}$.
\item[Plus:] $\Lambda + \Gamma := (\Lambda^L \ar \Gamma^L)^+$.
\end{description}
Moreover, given $x \in M$, define $x + \Lambda := (x \ar \Lambda^L)^+$.%
\footnote{$\cut{x \ar \Lambda^L}{x \ar \Lambda^R}$ is not a cut in general.}
Finally, $x < \Lambda$ if $x \in \Lambda^L$, and $x > \Lambda$ if $x \in \Lambda^R$.
\end{definizione}
\begin{remark}\label{REM:MH-MINUS}
If $M$ is a \predom, then $\mh$ is also a \predom.
Moreover, $\Lambda \ar \Gamma = (\Lambda^R + \Gamma ^R)^-$, and $\Lambda - \Gamma = (\Lambda^R \dl \Gamma^L)^-$.
\end{remark}
\begin{proof}
It is enough to prove that $\Lambda + \zOh = \Lambda$ and Axiom~\ref{AX:PA}.
In fact, $\Lambda + \zOh = (\Lambda^L \ar \zOM)^+ = (\Lambda^L)^+ = \Lambda$.
Moreover, if for contradiction $\Lambda < \Lambda'$, but $\Lambda + \Gamma > \theta > \Lambda' + \Gamma$, then $\theta \leq \lambda \ar \gamma$ for some $\lambda < \Lambda$, $\gamma < \Gamma$.
Therefore, $\lambda < \Lambda'$, thus $\theta < \Lambda' + \Gamma$, absurd.
\end{proof}
By Remark~\ref{REM:COMPLETE-SET}, $\mh$ is complete, and hence compact in its interval topology.
\begin{proposition}
If $M$ is a \dom of the first type, then $\mh$ is also a \dom, satisfying
$\zD = \zOM^- < \zO$.
\end{proposition}
\begin{proof}
Axioms~\ref{AX:M1} and~\ref{AX:M4} are trivial.
\begin{description}
\item[\ref{AX:M5-LEFT}.] If $\Lambda - \Gamma < \zOh$, then $\zO \geq \lambda \dl \gamma$ for some $\lambda > \Lambda$, $\gamma < \Gamma$.
Since $\zOM = \zDM$, we infer that $\lambda  \leq \gamma$, and therefore $\Lambda < \Gamma$.
\item[\ref{AX:M5-RIGHT}.] If $\Lambda < \Gamma$, then there exists $x \in M$ such that $\Lambda < x < \Gamma$.
Therefore, $\zOM \geq x \dl x > \Lambda - \Gamma$, hence $\zOh > \Lambda - \Gamma$. \qedhere
\end{description}
\end{proof}
\begin{remark}
Let $M$ be a \dom of the first type, and $a \in M$;
then $\ig{a^+} \geq \Pa{\ig a}^-$ and $\ig{a^-} \geq \Pa{\ig a}^-$.
\end{remark}
\begin{proof}
By Remark~\ref{REM:MH-MINUS}, $\ig{a^+} = \set{x - a : a < x}^-$, and the conclusion follows. The other inequality is similar.
\end{proof}
Hence, we might have a discontinuity of the function $\mu(x) := a - x$ at the point $x = a$.
More precisely, $\mu(a) = \ig a \geq \zO$, while 
\[
\set{\mu(x):x>a}^+ = \set{a \ar x': x' < -a}^+ = a^+ + (-a)^- = -(\ig{a^+}) \leq \Pa{- \ig a}^+.
\]
Therefore, if $\ig a$ is sufficiently large (more precisely, if there exists $y \in M$ such that $\zO < \ig y < \ig a)$, we can have $\Pa{- \ig a}^+ < \zO^-$.

In this case, let $\Pa{-\ig a}^+ < z < \ig a$,
%\footnote{It suffices that $\ig a > \zO$ for }
and $\Lambda := \Pa{-x}^-$.
Hence, $a \ar x < z < \ig a$ for every $x < -a$, \ie\ $a \ar \Lambda^L < z < a \ar \Lambda^R$.
Therefore, $\cut{a \ar \Lambda^L}{a \ar \Lambda^R}$ is not a cut.
\begin{enumexamples}
\item
We left open the question of what happens if we choose a different definition of plus for $\mh$.
If we want $\zD \leq \zO$, the only other possible definition of plus is
$\Lambda \ad \Gamma := (\Lambda^L + \Gamma^L)^+$.
If $M$ is a \predom, then $\mh$, with this modified plus, is also a \predom.
However, if $M$ is a \dom of the first type, then $\mh$ will satisfy axioms~\ref{AX:M1}, \ref{AX:M4} and~\ref{AX:M5-LEFT}, but not~\ref{AX:M5-RIGHT}.
For instance, take $x \in M$ such that $\ig x > \zOM$, and let $\Lambda = x^-$, $\Gamma = x^+$. Then, $\Lambda - \Gamma = (\ig x)^- \geq \zO$, even though $\Lambda < \Gamma$.
\item
Let $\vg$ be an ordered group.
Therefore, $\vg$ is a \dom of the first type.
The definition of $\vgh$ in \S\ref{SEC:GROUP-BASIC} and the one given above
coincide.
\end{enumexamples}
\subsection{Shift}
\begin{definizione}
Assume that $M$ is a \dom of the second type.
Define $\dms x := \zOM \dm x$.
The \intro{shift} of $M$ is the structure
$\Ms := \struct(M, \leq_M, \zOM, \am, \dms)$.
\end{definizione}
\begin{definizione}\label{DEF:FIRST-SECOND}
Assume that $M$ is a \dom of the first type, with a minimal positive element $1$, such that $1 \dm 1 = \zOM$.
Define $\dms x := 1 \dm x$.
The \intro{shift} of $M$ is the structure
$\Ms := \struct(M, \leq_M, \zOM, \am, \dms)$.
\end{definizione}
If $M$ is of the third type, the shift of $M$ is $M$ itself.
\begin{lemma}
If $M$ is a \dom of the second type, then $\Ms$ is a \dom of the first type, satisfying the condition of Definition~\ref{DEF:FIRST-SECOND}.
If $M$ is a \dom satisfying the condition of Definition~\ref{DEF:FIRST-SECOND}, then $\Ms$ is a \dom of the second type.
In both cases, $x \dms y = x \dm y$  and $\Mss = M$.
Moreover, $\dms x = \dm x$ for every $x \in \notmzero$.
\end{lemma}
\begin{proof}
The conclusion is an immediate consequence of the facts that $\dms$ is an anti-automorphism of $\struct(\Ms, \dms)$, that $\dms(\dms x) = x$, and that $x \dms y = x \dm y$, which are left as an exercise.
\end{proof}

\begin{example}
Let $\vh$ be a discretely ordered group, $1$~be its minimal positive element, 
and $M := \zerogroup{\vhh}$.
%Note that $M$ is isomorphic to $H$ as a group (but with a different minus).
Then, $M = \vh^s$ and $\vh = \Ms$.
\end{example}

\begin{enumexamples}[%
Let $B$ be a densely ordered group, $A$~be a subgroup of~$B$,
and $M$ be a \dom, with a maximum~$\nuM$.
Note that $\bh$ is a \dom of the third type, and that $B = \collapse{\bh}$.
Call $\theta: A \hookrightarrow \collapse{\bh}$ the inclusion map.%
]
\item
$\insem {\bh}{A} = \crux A \theta {\bh}$ is the \subdom of $\widetilde B$ given by $A \sqcup \bh$.
%\item
% = \crux A \theta {\bh}$.
\item
$B \fprod M$ is the \dom (of the same type as~$M$), constructed in this way:
starting from $B$, substitute every point of $A$ with a copy of $M$.
If $b \in B \setminus A$, then $\ig b = (0, \nuM)$ inside $B \fprod M$.
\item $\coll{\bh}{A}$ is the \dom (of the third type) constructed in this way:
starting from~$\bh$, identify $b^+$ with $b^-$, for every $b \in B \setminus A$.
\end{enumexamples}
In the following section we will apply most of the constructions showed in this one.

\section{Embedding \doms in cuts of groups}\label{SEC:EMBEDDING}
\abstractfont{%
The aim of this section is proving that every \dom $M$  can be embedded 
in a \dom of the form either $\widetilde G$ or~$\vgh$, for some ordered (Abelian) group~$\vg$.
We will deal first with proper \doms, and then prove the general case.}

%\subsection[First and second type]
\subsection{Proper \doms of the first and second type}
Any ordered group $\vh$ is, by definition, a subgroup of $\zerogroup{\vht} = \set{ x \in \vht: \ig x = \zO}$;
it is trivial to see that it is actually equal to $\zerogroup{\vht}$.
Let $\tau: \vh\to\zerogroup{\vht}$ be the natural \dom-isomorphism.

Consider the map $\tilde\iota$ from $H$ to $\G{\vht}$, obtained by composing the embedding of $\vh$ in $\vht$ with the quotient map $\pi$ from $\zerogroup{\vht}$ to $\G\vht$.
It is easy to see that $\tilde\iota$ is an isomorphism of ordered groups, and hence we can identify canonically $\vh$ with $\G\vht$.

There exists also a map $\iotahat$ from $\vh$ to $\G{\vhh}$, the composition of the map from $\vh$ in $\vhh$ sending $x$ to $x^+$ with the quotient map $\pi$ from $\zerogroup{\vhh}$ to $\G\vhh$.
It is easy to see that %if $H$ is densely ordered, then 
$\iotahat$ is an injective homomorphism of ordered groups, and hence we can also identify canonically $\vh$ with a subgroup of $\G\vhh$.
%\footnote{
In general, $\iotahat$ is not surjective. For instance, $\G{\Dedekind\Raz} = \Real$.

If $\vh$ is \emph{discrete} and non-trivial, let $1$ be its minimal positive element.
In this case, we must pay attention to the fact that the minus $-$ of $\G\vhh$ induced by the quotient map $\pi$ is not the group minus (let us call it $\md$), but instead $-x  + 1 = \md x $.
Hence, in that case, if we want $\iotahat$ to be a \domhom, we must either substitute $\md$ to the minus on $\G\vhh$,
or $-x + 1$ to the minus on $\vh$.

If instead $\vh$ is \emph{densely ordered}, $\iotahat$ is a \domhom.
Moreover, the group $\G\vhh$ coincides with the completion of $\vh$ via Cauchy sequences; \cf~\cite[\S~V.11]{FUCHS:1963} (see also~\cite{SCOTT:1969} for the completion of ordered fields).

In general, the image of $\Mgroup$ in $\Mquot$ under the quotient map $\pi$ is exactly $\GM$.

Besides, if $f :M \to N$ is a \domhom, then $f(\Mgroup)\subseteq\zerogroup N$.
\begin{thm}\label{THM:M-CANONICAL}
Let $M$ be a proper \dom, and $\vg := \GM$  the group associated to $M$.
If $M$ is of the first or second type, then there exists a unique homomorphism of ordered monoids $\psitilde : M \to \vgt$ such that the following diagram commutes:
%\\
\begin{equation}\label{DIAG:COMM-1}
\begin{minipage}{100pt}\vspace{1.5ex}
\begin{diagram}
\Mgroup  ¤ \Edotar{\psitilde^0} ¤ \zerogroup{\vgt} ¤¤ %%
\Sar\piO ¤ \neeqL\tau ¤¤
\vg ¤¤
\end{diagram}
\end{minipage}
\end{equation}
%\hfill\rom($1.1$\rom).
%
If $M$ is of the second type, then there exists a unique \domhom $\psihat : M \to \vgh$ such that the following diagram commutes:
\begin{equation}\label{DIAG:COMM-2}
\begin{minipage}[c]{100pt}\vspace{1ex}
\begin{diagram}
\Mgroup  ¤ \Edotar{{\psihat}^0} ¤ \zerogroup{\vgh} ¤¤ %%
\Sar\piO ¤ ¤ \Sar\piO ¤¤
\vg ¤ \Emono{\iotahat} ¤ \G{\vgh} ¤¤
\end{diagram}
\end{minipage}
\end{equation}
Moreover, $\psitilde$ and $\psihat$ are injective, and if $M$ is of the first type, then $\psitilde$ is actually a \domhom.
\end{thm}
In the diagrams above, $\psitilde^0$ (resp. ${\psihat}^0$) is the restriction of $\psitilde$ (resp. $\psihat$) to $\Mgroup$;
besides, the maps denoted by $\piO$, which are the restrictions of the quotient maps, are \domhoms;
the map $\iotahat$ becomes a \domhom if the structure of $\vg$ or $\G{\vgh}$ is modified as described above.
%and all maps (except possibly $\psitilde^0$) are \domhoms.%
%\footnote{Taking into consideration the warning before the theorem about $\iotahat$ when $\vg$ is discrete.}
%
\begin{proof} Let $x \in \notmzero := M \setminus\Mgroup$.
For such $x$, $\psitilde(x)$ and $\psihat(x)$ must coincide, and be equal to $\Lambda(x)$ (cf.\ \ref{DEF:LAMBDA}).

It remains to show the existence and uniqueness of the extensions of $\Lambda$ to $M$ in the various cases.
Let $x \in \Mgroup$: we have to define the image of $x$.

In the case when $M$ is of the first or second type, $\piO$ is an isomorphism of ordered groups.
%$\Mgroup$ is isomorphic to $\vg$,  as ordered groups. Let us call $\sigma$ this isomorphism: 
The only possible way to extend $\Lambda$ to $\psitilde$ is by defining $\psitilde(x) := \pi (x) = [x]$.
It is now trivial to see that $\psitilde$ is an injective homomorphism of ordered monoids.
Moreover, when $M$ is of the first type, the group-minus and the \dom-minus
on $G$ coincide; therefore, in that case $\psitilde$ is a \domhom, too.

Since $\psihat$ is a \domhom, $\psihat(\zO) = 0^+$.
%In the case when $M$ is of the second type, 
Moreover, if $M$ is of the first or second type, the quotient map $\pi$ is injective, hence the only value for $\psihat$ that makes the Diagram~\ref{DIAG:COMM-2} commute is $x^+$.
%In the case when $M$ is of the third type, let $y := \eqclass x \in \vg$.
%The only possible values for $\hat\Psi(x)$ that make the diagram commute are $y^-$ and $y^+$.
%There are two cases: the cardinality of ${\eqclass x}$ is either $2$ or $1$.
%
%When $y = \set{x_1, x_2}$, with $x_1 < x_2$, then $\hat\Psi(x_1) = y^-$, $\hat\Psi(x_2) = y^+$.
%
%When $y = \set{x}$, then $\hat\Psi(x) = \hat\Psi(x \ar \zO) = \hat\Psi(x) \ar 0^+ = y^+$.

It is now easy to see that $\psihat$ so defined is indeed an injective \domhom.
\end{proof}
\begin{enumexample}
\label{EX:RAZ-2}
%% Let $M := \Raz^\star \sqcup \fdom{2}$, with order and minus extending the ones of $\Raz$ and $\fdom{2}$, with the following additional rules: for every $x,y,z\in\Raz^\star$
%% \[\begin{aligned}
%% \zO &= \zOdue,\\
%% %-0^+ &= 0^-,\\
%% %0^- & < 0^+, \\
%% x & < \zDdue \text{ iff } x < 0,\\
%% x & > \zOdue \text{ iff } x > 0,\\
%% x + y & = z \text{ iff } x + y = z \text{ in } \Raz^\star \text{ and }  z \neq 0\\
%% x + (-x) & = \zDdue,\\
%% %x + 0^+ & = x,\\
%% x + \zDdue &= x + \zOdue = x.
%% %0^+ + 0^+ &= 0^+,\\
%% %0^+ + 0^- &= 0^-,\\
%% %0^- + 0^- &= 0^-.
%% \end{aligned}\]
Let $M := \Raz \fprod[\mset{0}] \fdom{2}$.
That is, $M$~is the \dom of the third type 
obtained from~$\Raz$ by duplicating the element~$0$.
Note that $M$ is strongly proper, because, for every $x \in M$, $\ig x = \zO$.
However, if $\vg = \G M = \Raz$, $M$ cannot be embedded in $\vgt$ or $\vgh$ in a way that makes the diagrams~\ref{DIAG:COMM-1} or~\ref{DIAG:COMM-2} commute.%
\footnote{However, if we do not insist on the corresponding diagram to commute, we can define an embedding of $M$ into $\Razhat$, by sending $1$ into $\pi$ (or any irrational positive element).}
\end{enumexample}

%%% Local Variables: 
%%% mode: latex
%%% TeX-master: "main"
%%% End: 

\subsection{Proper \doms of the third type}
Given a proper \dom of the third type~$M$, we want to construct a densely ordered group~$\vg$, and an embedding of $M$ into~$\vgh$.
To see where a difficulty of the task lies, and to get an idea of how we proceed in solving it, the reader can try his hand at the following exercise.
\begin{exercise}
Let $M := \Real \fprod[\Zed] \fdom{2}$.
That is, start from the group $\Real$, and obtain $M$ by duplicating all the natural numbers.
Note that $M$ is a proper \dom of the third type.
Find an embedding of $M$ into some~$\vgh$.
\end{exercise}
%%
%\abstractfont{%
%In this  subsection we will deal with proper \doms of the third type.}
\begin{lemma}\label{LEM:4-GROUPS}
Let $\vk$ be an ordered group, and $A$, $B$ and $C$ be subgroups of $\vk$, such that $A \cap B = C$, and for every $0<\varepsilon \in \vk$ and $a \in A$ there exist $b$ and $b' \in B$ such that
\[
a - \frac{\varepsilon}{2} < b < a < b' < a + \frac{\varepsilon}{2}.
\]
Then, if $\vh$ is any ordered group, there exists an injective \domhom
$\psi: A \fprod[C] \vhh \to \wideDed{\Bth}$.
\end{lemma}
We should interpret the inequality $a - \half{\varepsilon} < b$ as either taking place in the divisible hull of $\vk$, or as a shorthand for $2(a - b) < \varepsilon$.
\begin{proof}
To simplify the notation, we will use letters $a$, $a'$,~\dots for elements of~$A$; $b$, $b'$,~\dots for elements of~$B$; $c$, $c'$,~\dots for elements of~$C$; and $h$, $h'$,~\dots for elements of~$\vh$.
\begin{claim}\label{CL:DENSE}
For every $a < a' \in A$ there exists $b \in B$ such that $a < b < a'$.
\end{claim}
Apply the hypothesis with $\varepsilon := a' - a$.

Define a map $\Thm: A \to \bh$ sending $a$ to ${\set{b \in B: b < a}}^+$.
\begin{claim}\label{CL:THM-PLUS}
$\Thm$ is injective and preserves the sum and the order (but not the minus).
\end{claim}
In fact, Claim~\ref{CL:DENSE} implies immediately that if $a < a'$, then $\Thm(a) < \Thm(a')$, hence $\Thm$ is injective and preserves the order.
Moreover,
\[
\Thm(a) + \Thm(a') = {\set{b \in B: b < a}}^+ + {\set{b' \in B: b' < a'}}^+
\leq \Thm(a + a').
\]
If, for contradiction, $\Thm(a) + \Thm(a') < b'' < \Thm(a + a')$, then $b'' < a + a'$ and $b'' > b + b'$ for every $b$ and $b'\in B$ such that $b < a$ and $b' < a'$.
Let $\varepsilon := a + a' - b'' > 0$, and $b$, $b' \in B$ such that $a - \half\varepsilon < b < a$ and $a' - \half\varepsilon < b' < a'$.
Then, $b'' = a + a' - \varepsilon < b + b'$, absurd.

For every $a \in A$ and $\Lambda \in \vhh$, define
\[
\psi(a, \Lambda) := {\set{(b,h) \in \Bth: (b,h) < (a, \Lambda)}}^+,
\]
where the $<$ is the lexicographic order on $\vk \times \vht$.
We have to prove that $\psi$ is indeed an injective \domhom.
\begin{claim}
$\psi(a, \Lambda)$ is the cut
\[
\Cut{\set{(b,h) \in \Bth: (b,h) < (a, \Lambda)}} {\set{(b,h) \in \Bth: (b,h) > (a, \Lambda)}}.
\]
\end{claim}
Assume, for contradiction, that there exists $(b,h) \in \Bth$ such that $(b',h') < (b,h) < (b'',h'')$ for every $(b',h') < (a, \Lambda)$ and $(b'',h'') > (a, \Lambda)$.
Fix $0 < \varepsilon  \in \vk$, and let $b'$, $b'' \in B$ such that 
$a - \half\varepsilon < b' < a < b'' < a + \half\varepsilon$.
Hence, $b' \leq b \leq b''$, thus $\abs{a - b} < \varepsilon$.
We conclude that $a = b \in C$.
Hence, $\Lambda < h < \Lambda$, absurd.

Let $\mu$ be the minimum of $\vhh$.
Therefore, if $a \in C$, then
\begin{multline*}
\psi\Pa{-(a, \Lambda)} =\\
=\Cut{\set{(b,h) \in B \times H: (b,h) < (-a, -\Lambda)}} {\set{(b,h) \in B \times H: (b,h) > (-a, -\Lambda)}} =\\
= - \Cut{\set{(b,h) \in B \times H: (b,h) < (a, \Lambda)}} {\set{(b,h) \in B \times H: (b,h) > (a, \Lambda)}} =\\
= - \psi(a, \Lambda).
\end{multline*}

If instead $ a \notin C$, then
\begin{multline*}
\psi\Pa{-(a,\mu)} = \psi(-a, \mu)= \\
= \Cut{\set{(b,h) \in \Bth: b < -a}}
{\set{(b,h) \in B \times H: b > -a}} = -\psi(a, \mu).
\end{multline*}
Hence, $\psi$ preserves the minus.

Moreover, if $(a, \Lambda) < (a', \Lambda')$, then either $a < a'$, or $a = a' \in C$ and $\Lambda < \Lambda'$.
In the first case, let $b \in B$ such that $a < b < a'$.
Hence, $\psi(a, \Lambda) < (b, h) < \psi(a', \Lambda')$ for every $h \in H$.
In the second case, let $h \in H$ such that $\Lambda < h < \Lambda'$.
Hence, $\psi(a, \Lambda) < (a, h) < \psi(a', \Lambda')$.
Therefore, we conclude that $\psi$ is injective and preserves the order.

Moreover, $\psi$ preserves the zero, because
\[
\psi(0, 0^+) = {\set{(b,h) \in \Bth: (b,h) \leq (0,0)}}^+ = (0, 0)^+.
\]
\begin{claim}\label{CL:SUM-1}
For every $a\in A$ and $c \in C$,
\begin{multline*}
\set{b + c: b < a}^+ = \Thm(a + c) = \\
= \Thm(a) + \Thm(c) = \set{b \in B: b < a}^+ + \set{b' \in B: b' < c}^+. 
\end{multline*}
\end{claim}
It suffices to prove the first equality: the others follow from the definition of $\Thm$ and Claim~\ref{CL:THM-PLUS}.
The fact that the \lhs is less or equal to the \rhs is trivial.
Assume for contradiction that $\set{b + c: b < a}^+ < b'' < \set{b'\in B: b' < a + c}^+$ for some $b'' \in B$.
Hence, for every $b \in B$ such that $b < a$, we have $b + c < b'' < a + c$, \ie\ $b < b'' - c < a$.
But $b'' - c$ is in $B$, and we have an absurd.

It remains to show that $\psi$ preserves the sum, namely
$
\psi(a, \Lambda) + \psi(a', \Lambda') = \psi\Pa{(a,\Lambda) + (a', \Lambda')}
$.
We will make a case distinction.
\begin{itemize}
\item
If $\Lambda = \mu$, while $\Lambda' > \mu$, then
\begin{multline*}
\psi(a, \mu) + \psi(a', \Lambda') =
\set{(b + b', h + h'): b < a \et (b', h') < (a', \Lambda')}^+ = \\
= \set{(b + a', h + h'): b < a \et h' < \Lambda'}^+ =
\set{(b + a', h): b < a}^+.
\end{multline*}
By Claim~\ref{CL:SUM-1}, the latter is equal to 
\[%\begin{multline*}
\set{(b'' , h): b'' < a + a'}+=
%\Thm(a + a') \times \vh =\\
\psi(a + a', \mu) = \psi\Pa{(a, \mu) + (a', \Lambda')}.
\]%\end{multline*}
\item
If $\Lambda = \mu = \Lambda'$, then
\begin{multline*}
\psi(a, \mu) + \psi(a', \mu) =\\
=\set{(b , h) \in \Bth: b < a}^+ +
\set{(b' , h') \in \Bth: b' < a'}^+ =\\
\set{(b + b', h'') \in \Bth: b < a \et b' < a'}^+.
\end{multline*}
By Claim~\ref{CL:THM-PLUS}, the latter is equal to
\[%\begin{multline*}
%\Pa{\Thm(a) + \Thm(a')} \times \vh =
%\Thm(a + a') \times \vh =
\set{(b'',h'') \in \Bth: b'' < a + a'}^+=
\psi(a + a', \mu)= \psi\Pa{(a, \mu) + (a', \mu)}.
\]%\end{multline*}
\item
If $a, a' \in C$, and $\Lambda, \Lambda' > \mu$, then
\begin{multline*}
\psi(a, \Lambda) + \psi(a', \Lambda') =
\set{(a, h) + (a', h'): h < \Lambda \et h'< \Lambda'}^+=\\
\set{(a + a', h''): h''< \Lambda + \Lambda'}^+ =
\psi(a + a', \Lambda + \Lambda').
\qedhere
\end{multline*}
\end{itemize}
\end{proof}
Let $A$ be any ordered group and $C$ be a subgroup of~$A$.
Let $\Raz[\varepsilon]$ be the ring generated by $\Raz$ and a positive infinitesimal element $\varepsilon$.
Let $\vk := \Raz[\varepsilon] \otimes_\Zed A$, with the ordering given by
$a \gg \varepsilon \gg \varepsilon^2 \dotsc $ for every $0 < a \in A$.
Let $B$ be the following subgroup of $\vk$:
\[
B := C + \sum_{n \in \Nats} A \Pa{1 + (-\varepsilon)^n} = C + A(1 - \varepsilon) + A(1 + \varepsilon^2) + \dotsb.
\]
If we identify $A$ with the subgroup $A \cdot 1$ of $\vk$, we have that $A$, $B$, $C$, and $\vk$ satisfy the hypothesis of Lemma~\ref{LEM:4-GROUPS}.
In fact, $A \cap B = C$.
Moreover, if $\varepsilon^{2n}$ is a small positive element of $\vk$ (for some $n \in \Nats$), and
$a \in A$, then $a(1 - \varepsilon^{2n+1})$ and $a(1 + \varepsilon^{2n+2})$ are both in $B$.

Hence, for any ordered group $\vh$, we can embed $A \fprod[C] \vhh$ into $\wideDed{\Bth}$.
Let $M$ be a proper \dom of the third type, $A := \GM$, and $C := \HM$.
Choose $B$ and $\vk$ such that $A$, $B$, $C$ and $\vk$ satisfy the
hypothesis of Lemma~\ref{LEM:4-GROUPS} (for instance,
$\vk := \Raz[\varepsilon] \otimes_\Zed A$ and $B := C + \sum_{n \in \Nats} A \Pa{1 + (-\varepsilon)^n}$ as in the construction above).
Therefore, if we choose $\vh := \mset{0}$, the we can embed
$A \fprod[C] \fdom{2}$ into $\wideDed{B \times \mset{0}} = \bh$.
Moreover, $A \fprod[C] \fdom{2}$ is canonically isomorphic to $\Mgroup$, hence we can embed $\Mgroup$ into $\bh$ via a map $\psi_0$.

The next step is defining an embedding $\notpsi$ of $\notmzero$ in $\bh$, where $B = B(M)$ is the group constructed above.
Let $\sigma: B \to A$ be the map sending $\sum_{i < N} a_i \varepsilon^i$ to $a_0$.
Obviously, $\sigma$ is a surjective homomorphism of ordered groups, hence, by Remark~\ref{REM:PI-HAT}, it induces an injective \qdomhom $\sigmahat: \ah \to \bh$.
Moreover, by Lemma~\ref{LEM:LAMBDA}, there is an injective \qdomhom $\Lambda: \notmzero \to \ah$.
Define $\notpsi: \notmzero \to \bh$ as $\sigmahat \circ \Lambda$.
Hence, $\notpsi$ is also an injective \qdomhom.

Define $\psihat: M \to \bh$ as
\[
\psihat(x) =
\begin{cases}
\psi_0(x) & \text{if } x \in \Mgroup,\\
\notpsi(x) & \text{if } x \in \notmzero.
\end{cases}
\]
The final step is proving that $\psihat$ is indeed an injective \qdomhom.
We will need that $M$ is \emph{strongly proper} (and not simply proper) to do that.
Note that, for every $x \in \Mgroup$,
\[
\psi_0(x) = 
\begin{cases}
\set{b \in B: b < [x]}^+ & \text{if } x = \FM(x),\\
[x]^+ & \text{if } [x] \in C \et x = \FP(x),
\end{cases}\]
or equivalently
\[
\psi_0(x) = 
\begin{cases}
\set{b \in B: b < [x]}^+ =  \set{b \in B: b > [x]}^- & \text{if } [x] \in A \setminus C,\\
[x]^- & \text{if } [x] \in C \et x = \FM(x),\\
[x]^+ & \text{if } [x] \in C \et x = \FP(x).
\end{cases}\]
Moreover, if $y \in \notmzero$, then
\begin{multline*}
\notpsi(y) = \set{b \in B: \sigma(b) < y}^+ =\set{b \in B: \sigma(b) > y}^- = \\
= \set{[x_0] + [x_1] \varepsilon + \dotsb + [x_N] \varepsilon^N \in B: x_0 < y}^+.
%= \inf\set{[x_0] + [x_1] \varepsilon + \dotsb + [x_N] \varepsilon^N \in B: x_0 > y}.
\end{multline*}
Obviously, $\psihat$ preserves the minus, because that is true for both $\notpsi$ and $\psi_0$.

We will now prove that $\psihat$ is injective.
This is a consequence of the following claim.
\begin{claim}
For every $x \in \Mgroup$ and $y \in \notmzero$,
\[\begin{aligned}
\psi_0(x)  & \in \Bgroup, \\
\notpsi(y) & \in \notbgroup.
\end{aligned}\]
\end{claim}
Since $\ig x = 0^+$, $\ig{\psi_0(x)} = \psi_0(\ig x) = \psi_0(0^+) = 0^+$.
Since $\ig y > 0^+$ and $M$ is strongly proper, there exists $x \in \Mgroup$ such that $0 < x < \ig y$.
Let $b := [x](1 + \varepsilon^2) \in B$.
Then, $\sigma(b) = [x] < [\ig y]$, hence $\ig{\notpsi(y)} = \notpsi(\ig y)  > b > 0^+$.

We will now prove that $\psihat$ preserves the order.
It suffices to prove that for every $x \in \Mgroup$ and $y \in \notmzero$, if $x < y$, then $\psi_0(x) < \notpsi(y)$ (the other possibility $x > y$ is proved in a similar way).
Let $b := [x] (1 + \varepsilon^2) \in B$.
Then, $\sigma(b) = [x] < y$, hence $b < \notpsi(y)$.
Moreover, $b > [x]$, thus $b > \psi_0(x)$.
Therefore, $\psi_0(x) < b < \notpsi(y)$.

Finally, we have to prove that $\psihat$ preserves the sum.
It suffices to prove that, for every $x \in \Mgroup$ and $y \in \notmzero$,
\[
\psi_0(x) + \notpsi(y) = \notpsi(x + y).
\]
First, note that
\[
\set{b \in B: \sigma(b) < [x]}^+ < \psi_0(x) \leq \set{b \in B: \sigma(b) \leq [x]}^+.
\]
Hence,
\begin{multline*}
\psi_0(x) + \notpsi(y) \leq
\set{b + b': \sigma(b) \leq [x] \et \sigma(b') < [y]}^+ \leq\\
\leq \set{b'' \in B: \sigma(b'') \leq [x + y] }^+ = \notpsi(x + y).
\end{multline*}
Moreover, $ [x](1 - \varepsilon^{2k + 1})< \psi_0(x) < [x](1 + \varepsilon^{2k})$, for every $k \in \Nats$.
Assume, for contradiction, that there exists $b'' \in B$ such that $\psi_0(x) + \notpsi(y) < b'' < \notpsi(x + y)$.
Let $\gamma := \sigma(b'') \in A$.
Then, $\gamma < [x + y]$.
Let $\lambda := \gamma - [x] \in A$.
Hence, $\lambda + [x] = \gamma < [x] + [y]$, thus $\lambda < [x]$.
Therefore,
\[\begin{aligned}
\lambda(1 + \varepsilon^2) &< \notpsi(y),\\
[x]( 1 - \varepsilon^3) &< \psi_0(x).
\end{aligned}\]
Thus,
\[
\gamma > \lambda(1 + \varepsilon^2) + [x]( 1 - \varepsilon^3) =
(\lambda + [x]) + \lambda\varepsilon^2 - [x] \varepsilon^3 > \lambda + [x] = \gamma,
\]
absurd.
Therefore, we have proved the following lemma.
\begin{lemma}
Let $M$ a strongly proper \dom of the third type.
Then, if $B = B(M)$ is the group constructed above,
%there exists an ordered group $B$ such that 
$M$ can be embedded in $\bh$.
\end{lemma}
What happens if $M$ is a \dom of the third type, which is proper, but not strongly proper (for instance, $M = \fdom{4}$)?
Then, first we embed $M$ in a strongly proper \dom of the same type $M'$, and then we apply the above lemma to embed $M'$ in $\bh$ for a suitable ordered group $B$.
More precisely, let $M$ a proper \dom of the third type.
Let $N := \Mquot \fprod[] \Real$; note that $N$ is a \dom of the first type, and that $\Mquot$ is a \subdom of $N$, via the map sending $y$ to $(y,0)$ if $y \in \Mgroup$, and to $(y,\mu)$ otherwise.
Define
$M' := N \fprod[\Mquot] \fdom{2}$.
\begin{lemma}
The above defined $M'$ is a strongly proper \dom of the third type.
Moreover, the map $\xi : M \to M'$ defined by
\[
\xi(x) = \begin{cases}
([x], \mu, \zD) & \text{if } x \in \notmzero,\\
([x], 0, \zO)   & \text{if } x \in \Mgroup \et \sig x = 1,\\
([x], 0, \zD)   & \text{otherwise}
\end{cases}\]
is a \domhom.
\end{lemma}
Note that we could have used any non-trivial ordered group instead of~$\Real$.
\begin{proof}
Trivial verifications, using Lemma~\ref{LEM:PRODUCT}.
\end{proof}
Hence, we have proved the following lemma:
\begin{lemma}\label{LEM:EMBEDDING-THIRD}
Let $M$ a proper \dom of the third type.
Then, $M$ can be embedded in $\bh$ for some densely  ordered group $B$.
\end{lemma}
\subsection{Abelian extensions of groups}
\abstractfont{%
The definitions and facts in this subsection can be found in any book on homological algebra, and will be used in the next subsection.
We will use~\cite[Ch.~9]{FUCHS:1970-73} as a reference on extensions of groups.%
\footnote{We recall that for us all groups are Abelian.}
The reader can also consult~\cite{MACLANE:1967}.%
}
\begin{definizione}[Factor sets]
Let $A$ and $C$ be (Abelian) groups.
A \intro{factor set} is a map $\factor: C \times C \to A$ such that, for every $x$, $y$, $z \in C$,
\begin{itemize}
\item $\factor[x,y] = \factor[y,x]$;
\item $\factor[x,0] = \factor[0,x] = \factor[0,0] = 0$;
\item $\factor[y,x] + \factor[x,y+z] = \factor[x,y] + \factor[x+y,z]$.
\end{itemize}
Given such a factor set, the \intro{crossed product} of $C$ and $A$
is the group $\crossedCAf$, whose underlying set is $C \times A$,
and whose sum is defined as follows:
\[
(c, a) + (c', a') = \Pa{c + c', a + a' + \factor[c,c']}.
\]
\end{definizione}
The reader can verify that $\crossedCAf$ is indeed an Abelian group, with neutral element $(0,0)$.
Moreover, the maps $\iota: A \to \crossedCAf$ sending $a $ to $(0,a)$ and $\pi: \crossedCAf \to C$ sending $(c,a)$ to $c$ are group-homomorphisms.
Finally, the sequence
\[
0 \arrow A \Ar{\iota} \crossedCAf \Ar{\pi} C \arrow 0
\]
is exact.

Conversely, given any exact sequence of groups
\begin{equation}\label{SEQ:E}
0 \arrow A \Ar{\iota'} B \Ar{\pi'} C \arrow 0, \tag{E}
\end{equation}
a \intro{section} is a map $s : C \to B$ that fixes $0$ and is a right inverse of $\pi'$: \ie, $s(0) = 0$ and for every $c \in C$, $\pi'(s(c)) = c$.
Given such a section $s$, the differential of $s$ is the factor set $\des$ defined as follows: $\des(x,y) = \des(x) + \des(y) -\des(x + y)$.
One can verify that $\des$ is indeed a factor set, and that the map
$\beta: B \to \crossedCAdes$ sending $b$ to $\Pa{\pi'(b), b - s(\pi'(b))}$ is an isomorphism of groups, such that the following diagram commutes:
\begin{diagram}
A    ¤ \Ear{\iota'} ¤ B ¤ \Ear{\pi'} ¤ C  ¤¤%%
\seql ¤    ¤\Sar{\beta} ¤            ¤ \seql ¤¤
A    ¤ \Ear{\iota}  ¤ \crossedCAdes ¤ \Ear{\pi} ¤ C ¤¤ %%
\end{diagram}
Hence, given an exact sequence~\eqref{SEQ:E}, \wloG we can assume that $B = \crossedCAf$ for some factor set $\factor$.

We will need the following proposition (\cite[Proposition~24.6]{FUCHS:1970-73} and \cite[Corollary~III.3.8]{MACLANE:1967}).
\begin{proposition}\label{PROP:MACLANE}
Let $0 \to A \Ar{\iota} B \Ar{\pi} C \to 0$ be an exact sequence of groups, and
$\gamma: C \to C'$ be an injective homomorphism of groups.
Then, there exist an exact sequence 
$0 \to A \Ar{\iota'} B' \Ar{\pi'} C' \to 0$ and an injective group-homomorphism $\beta: B \to B'$
such that the following diagram commutes:
\begin{diagram}
A    ¤ \Emono{\iota} ¤ B ¤ \Eepi{\pi} ¤ C  ¤¤
\seql ¤    ¤\Smono{\beta} ¤      ¤ \Smono{\gamma} ¤¤
A    ¤ \Emono{\iota'}  ¤ B' ¤ \Eepi{\pi'} ¤ C' ¤¤
\end{diagram}
\end{proposition}
By the above considerations, we can assume that $B'$ is of the form $\crossedCpAf$ for some factor set $\factor$.

Assume now that $A$ and $C$ are ordered groups, $B$ is a group, and that we have an exact sequence~\eqref{SEQ:E}.
Then, there exists a unique ordering on $B$ such that all maps in~\eqref{SEQ:E} are homomorphism of ordered groups.
The ordering is defined by:
$b \leq b'$ iff $\pi(b) < \pi(b')$ or $\pi(b) = \pi(b')$ and $b - b' \leq 0$ (in $A$).
In particular, on $\crossedCAf$ the ordering is the lexicographic one.
Moreover, with the ordering defined above, $A$ is a convex subgroup of $B$.

Moreover, in the situation of Proposition~\ref{PROP:MACLANE}, if $C'$ is an ordered group, and $\gamma: C \to C'$ is also a homomorphism of ordered groups, then $\beta: B \to B'$ is also a homomorphism of ordered groups, where $B$ and $B'$ are endowed with the above defined orderings.
Hence, we have proved the following:
\begin{corollary}\label{COR:CROSSED-PRODUCT}
Let $\theta: B \to D$ be a homomorphism of ordered groups.
Then, $B$ can be embedded in the crossed product $\crossed{D}{\ker \theta}{\factor}$ for some factor set $\factor$, where  $\crossed{D}{\ker \theta}{\factor}$ is endowed with the lexicographic ordering.
Namely, the following diagram of ordered groups with exact rows commutes:
\begin{diagram}
0 ¤ \ear ¤ \ker(\theta)    ¤ \emono ¤ B ¤ \eepi ¤ \im(\theta)  ¤ \ear ¤ 0 ¤¤
¤ ¤ \seql ¤    ¤\Smono{\beta} ¤ \Sear{\theta} ¤ \smono ¤¤
0 ¤ \ear ¤ \ker(\theta)    ¤ \Emono{\iota}  ¤ \crossed{D}{\ker \theta}{\factor} ¤ \Eepi{\pi} ¤ D ¤ \ear ¤ 0 ¤¤
\end{diagram}
\end{corollary}
\subsection{General case}
\abstractfont{%
In this subsection we will drop the ``properness'' hypothesis.%
}%
\begingroup
% Nomi

%\def\Mdom{dom\xspace}
%\def\Mdoms{doms\xspace}
\def\Mstrangeprod{crossed product\xspace}

%% Labels

\def\Mthatcorollary{Corollary~\ref{COR:CROSSED-PRODUCT}\xspace}%{COR:CROSSED-PRODUCT }
\def\Membeddingthirdtype{Lemma~\ref{LEM:EMBEDDING-THIRD}\xspace}%{EMBEDDINGTHIRDTYPE }
\def\Membeddingsecondtype{Theorem~\ref{THM:M-CANONICAL}\xspace}%{EMBEDDINGSECONDTYPE }
\def\Mlemmino{\S\ref{SEC:UNION}\xspace}%{Remark~\ref{REM:UNION}\xspace}%{LEMMINO }

%% frasi

\def\Mobda{It is sufficient to prove the \Mlemma when\xspace}
\def\Mst{s.t.\xspace}
\def\Mist{isomorphic to\xspace}
\def\Miff{iff\xspace}
\def\Mclaims{claims\xspace}

%%%\begin{document}

\begin{proposition} \label{Marchesato}
Let $T$ be any first-order theory,
and $S$ be a universal theory (i.e. axiomatised by a set of universal formulae) in the same language~$L$,
such that $T_\forall\vdash S$ (where $T_\forall$ denotes the universal part of $T$).
In this case, the following are equivalent:
\begin{enumerate}
\item $T_\forall = S$;
\item every model of $S$ is a substructure of some model of~$T$;
\item every finitely generated model of $S$ is a substructure of some model of~$T$.
\end{enumerate}
Suppose moreover that any of the above equivalent conditions is satisfied, and that $\Cfam$ is a class of $L$-structures, such that every model of $T$ can be embedded in some structure in~$\Cfam$.
Then,
\begin{enumerate}
\item every model of $S$ can be embedded in some structure in~$\Cfam$;
\item $S$ is the universal part of the theory of the structures in~$\Cfam$.
\end{enumerate}
\end{proposition}
\begin{proof}
Easy.
\end{proof}

\begin{MlemmaE} \label{Municipio}
Let $\MS$ be a \dom: if $\Mord{\MS}$ has a least non zero element
$\Minf$, and
$\Mres{\MS}{\ge\Minf}$ is proper, then $\MS$ is \Mist a \subdom of some proper
\dom~$T$.
\end{MlemmaE}

\begin{MclaimE} \label{Mpulce}
\Mobda $\MS$ is of the first type,
provided that \Mit{if} $\MS$ has a least positive element $1$, \Mit{
then} the constructed $\MT$ has a least positive element too, \Mit{and} $1$ is mapped to it.
\end{MclaimE}

\begin{proof}[Proof of \Mclaim \ref{Mpulce}]
\Mit{If $\MS$ is of the first type} then nothing should be proved.

\Mit{If $\MS$ is of the
second type} then consider a proper \dom $\MT$ \Mst $\Mshif{\MS}$ is \Mist a \subdom
of it. Then we claim: \Mit{$\MS$ is \Mist a \subdom of
$\Mshif{\MT}$},
which is trivial because of our assumption on the least elements, \Mit{
and $\Mshif{\MT}$ is proper}, which is trivial because
$a\Mopstr{-}{\MT}a=0$ \Miff $a\Mopstr{\Mminshft}{\MT}a=0$.

\Mit{If $\MS$ is of the third type} then consider a proper \dom $\MT$ \Mst $\MQuot{\MS}$ is \Mist a \subdom of it.
Then we claim: \Mit{$\MS$ is \Mist a \subdom of
$\Mprod{\MT}{G}{\Mmp}$, where $G$ is the isomorphic image of $\Macca{\MS}$},
which is trivial because $\MS$ is isomorphic to $\Mprod{\MQuot{\MS}}{\Macca{\MS}}{\Mmp}$, \Mit{
and $\Mprod{\MT}{G}{\Mmp}$ is proper}, which is trivial because $a-a$ is
$0$ in the product \Miff the first component is.
\end{proof}

\begin{MclaimE} \label{Moro}
\Mobda either $\Meg{\Mres{\MS}{\Minf}}$ is the whole
$\MQuot{\Mres{\MS}{\ge\Minf}}$ or it is dense in it.
\end{MclaimE}
\begin{proof}[Proof of \Mclaim \ref{Moro}]
If $\Mres{\MS}{\ge\Minf}$ is of the second type nothing should be proved. Otherwise,
by \Membeddingthirdtype, $\Mres{\MS}{\ge\Minf}$ is a \subdom of $\Mcuts{\MG}$
for some densely ordered group $\MG$, and by \Mlemmino $\MS$ may be embedded in the
\dom $\Mres{\MS}{0}\sqcup\Mcuts{\MG}$ (with appropriate definitions for
operations and order).
Now, by a general fact, if $\MG$ is dense,
then $\Macca{\Mcuts{\MG}}$ is dense in $\MQuot{\Mcuts{\MG}}$.
\end{proof}

\begin{proof}[Proof of \Mlemma \ref{Municipio}]
We give an explicit construction of a \dom $\MT$ \Mst~$\MS$ is a \subdom
of $\MT$, under the additional hypothesis stated in \Mclaims \ref{Mpulce} and \ref{Moro}.

Consider the exact sequence
\[
\Mzerog \rightarrow \Mker{\Mrho} \rightarrow \Mres{\MS}{0}
\xrightarrow{\Mrho} \Meg{\Mres{\MS}{\ge\Minf}}
\]
where $\Mrho$ is the function defined by
\[\begin{array}{r@{\,}c@{\,}c@{\,}l}
\Mrho :& \Mres{\MS}{0} &\to &\Meg{\Mres{\MS}{\ge\Minf}}\\
       & x &\mapsto &[x+\Minf].
\end{array}\]
%\Mrho : \Mres{\MS}{0} \ni x \mapsto [x+\Minf] \in \Meg{\Mres{\MS}{\ge\Minf}}
By \Mthatcorollary the group $\Mres{\MS}{0}$ (which is a group because
of \Mclaim \ref{Mpulce}) can be embedded in a \Mstrangeprod
\def\Mscp{\Mstrangetimes{\Meg{\Mres{\MS}{\ge\Minf}}}{\Mker{\Mrho}}{\factor}}
$\MP=\Mscp$
\def\Mscp{\MP}
for an appropriate choice of~$\factor$.
We know that the function
\[\begin{array}{r@{\,}c@{\,}c@{\,}l}
\Mrho' :& \Mscp &\to      & \Meg{\Mres{\MS}{\ge\Minf}}\\
        & (x,y) & \mapsto & x
%\Mrho' :
%\Mscp \ni (x,y) \mapsto x\in \Meg{\Mres{\MS}{\ge\Minf}}
\end{array}\]
consistently extends $\Mrho$, therefore, in the following,
$\Mres{\MS}{0}$ will be identified with a subgroup of $\Mscp$ and
$\Mrho$ will be identified with the restriction of $\Mrho'$ to
$\Mres{\MS}{0}$.

We claim that the structure $\MT :=
\Mdag\Mscp{\Mres{\MS}{\ge\Minf}}{\Mtheta_{\Minf}}$, where
$\Mtheta_{\Minf}(x) = [\Mrho(x) + \Minf]_{\Minf}$,
%$\forall i\in\Mord{\MS}\setminus\{0\}\;\Mtheta_i(x)=[\Mtheta(x)+i]$
is a super-\dom of $\MS$
(which is trivial because
$\MS = \Mdag{\Mres{\MS}{<\Minf}}{\Mres{\MS}{\ge\Minf}}{\Mtheta_\Minf} < \MT$), and
moreover it is proper and it verifies the hypothesis of \Mclaim~\ref{Moro}.

The only non-trivial point arises in proving the properness of $\MT$ when
$\Mres{\MS}{\ge\Minf}$ is of the third type.
In that case we need the density of $\Meg{\Mres{\MS}{\Minf}}$: given
$x, y \in \Mres{\MS}{\ge\Minf}$ \Mst $x<y$, \Mit{either} $x=\Mminus{z}$ and
$y = \Mplus{z}$ for some $z$, \Mit{or} there exists $z$ \Mst $[x]<[z]<[y]$,
\Mit{therefore}, given such a $z$, the inequality 
$x<(z, t) < y$ holds for an arbitrary choice of
$t \in \Mker{\Mrho}$.
\end{proof}

\begin{McorollaryE} \label{Mappropriamento}
Let $\MS$ be a \dom: if $\Mord{\MS}$ has finite cardinality
then $\MS$ is \Mist a \subdom of some proper \dom~$\MT$.
\end{McorollaryE}
The \Mcorollary is a trivial consequence of
\begin{MclaimE} \label{Mpappa}
Let $\MS$ be a \dom: \Mit{if} $\Mres{\MS}{\ge k}$ is proper for some
$k\in\Mord{\MS}$, \Mit{and}, for that $k$, $\{x\in\Mord{\MS}\;|\;x < k\}$
has finite cardinality, \Mit{then} $\MS$ is \Mist a \subdom of some
proper \dom~$\MT$.
\end{MclaimE}

\begin{proof}[Proof of \Mclaim \ref{Mpappa}]
%Consider a \dom $\MS$ \Mst
%$n=\#\{x\in\Mord{\MS}\;|\;x<k\}$ is minimal. Suppose, by contradiction,
%$n\ne 0$, 
Suppose, for contradiction, that there exists a \dom $\MS$ and $k\in\Mord{\MS}$  that satisfy the hypothesis of the claim, but not the conclusion.
We can assume that $n=\#\{x\in\Mord{\MS}\;|\;x<k\}$ is minimal.
Let
$\{k_0,\dotsc,k_{n-1}\}=\{x\in\Mord{\MS}\;|\;x<k\}$ with
$0=k_0<\dotsb<k_{n-1}$.
By \Mlemma \ref{Municipio}, $\Mres{\MS}{\ge k_{n-1}}$ is a \subdom of some
proper \dom~$K$, and, by \Mlemmino, $H=\Mres{\MS}{<k_{n-1}}\sqcup K$ is a
\dom extending~$\MS$. Moreover $\Mres{H}{\ge k_{n-1}} = K$ is proper by
construction, but $\{x\in\Mord{H}\;|\;x<k_{n-1}\}=\{k_0,\dotsc,k_{n-2}\}$
has cardinality $n-1$, contradicting the minimality of~$n$.
\end{proof}

\begin{McorollaryE} \label{Mimmfinord}
%Let $\MS$ be a \dom of the second or the third type: 
%if $\Mord{\MS}$ has finite cardinality, 
Let $\MS$ be a \dom, such that $\Mord{\MS}$ has finite cardinality.
Then, there exists an ordered group $\MG$ \Mst:
\begin{enumerate}
\item if $\MS$ is of the fist type, then $\MS$ is a \subdom of~$\vgt$;
\item
if $\MS$ is of the second or the third type, then $\MS$ is a \subdom of $\Mcuts{\MG}$.
\end{enumerate}
\end{McorollaryE}

\begin{proof}
Obvious because $\MS$ is a \subdom of some proper \dom $\MT$, and $\MT$
can be embedded in an appropriate $\Mcuts{\MG}$ or $\vgt$ either by
\Membeddingsecondtype or \Membeddingthirdtype.
\end{proof}

\begin{MtheoremE} \label{Munipart}
\begin{enumerate}
\item
The (first-order) theory of \doms of first type (axioms of \doms, plus $\zO = \zD$) is the universal part of the theory of the structures $\vgt$, as $\vg$ varies among ordered (Abelian) groups.
Moreover, every \dom of the first type is a \subdom of $\vgt$, for some ordered  group~$\vg$.
\item
The theory of \doms of second type (\predoms plus
\MMA with the strict inequality, \MMB and \MMC, and $\zO \ar \zO > \zO$) is the universal part of the theory of cuts of discretely ordered groups.
Moreover, every \dom of the second type is a \subdom of $\vgh$, for some discretely ordered  group~$\vg$.
\item
The theory of \doms of third type (\predoms plus
\MMA with the strict inequality, \MMB and \MMC, and $\zO \ar \zO = \zO$) is the universal part of the theory of cuts of densely ordered groups.
Moreover, every \dom of the second type is a \subdom of~$\vgh$, for some densely ordered  group~$\vg$.
\end{enumerate}
\end{MtheoremE}

\begin{comment}
We will give a proof of the first statement: the second one is identical.
Let $\MTS$ be the first theory and $\MTT$ the second one. Clearly it is
sufficient to prove that if
$\psi=\forall x_1,\dotsc,x_n\phi$ is an universal formula \Mst\
$\MTT\vdash\psi$, then $\MTS\vdash\psi$. For, if not, we consider a model
$\MMM$ of $\MTS+\neg\psi$ and interpretations $\bar{x_1},\dotsc,\bar{x_n}$
of the free variables of $\phi$ in that model \Mst\
$\MMM\models\phi[\bar{x_1}/x_1,\dotsc,\bar{x_n}/x_n]$: by restriction of
$\MMM$ to the submodel generated by $\bar{x_1},\dotsc,\bar{x_n}$ we get a model $\MMM'$ of
$\MTS$ satisfying the hypothesis of \Mcorollary \ref{Mimmfinord} but
still satisfying $\neg\psi$, which is a contradiction because, for
\Mcorollary \ref{Mimmfinord}, $\MMM'$ is a submodel of some model of $\MTT$.
\end{comment}
\begin{proof}
We will give a proof of the statement for \doms of the third type: the other cases are similar.

Let $T$ be the theory of proper \doms of the third type, $S$ be~the theory of \doms of the third type, and $\Cfam$ the class of of cuts of densely ordered Abelian groups.
By Corollary~\ref{Mappropriamento} every finitely generated model of $S$ can be embedded in a model of~$T$, and by Lemma~\ref{LEM:EMBEDDING-THIRD}, every model of $T$ can be embedded in a structure in~$\Cfam$.
Proposition~\ref{Marchesato} implies the conclusion.
\end{proof}
%%
\begin{comment}
\begin{MtheoremE} \label{Membedding}
Every \dom %$\MS$ 
of the second or third type is a \subdom of some $\Mcuts{\MG}$.
Every \dom of the first type is a \subdom of some $\vgt$.
\end{MtheoremE}
\begin{proof}
By Theorem~\ref{Munipart}, every \dom 
%of second or third type 
can be embedded in a proper \dom of the same type.
By Theorem~\ref{THM:M-CANONICAL}, every proper \dom of the first type is a \subdom of $\vgt$ for some ordered group~$\vg$.
By Theorem~\ref{THM:M-CANONICAL}, every proper \dom of the second type is a \subdom of $\Mcuts{\MG}$ for some discretely ordered group $\MG$.
By Lemma~\ref{LEM:EMBEDDING-THIRD}, every proper \dom of the third type is a \subdom of $\Mcuts{\MG}$ for some densely ordered group~$\MG$.
\end{proof}
\end{comment}
\endgroup

%%%
\subsection{Embedding a \dom in a collapse}
\label{SUBSEC:EMBEDDING-COLLAPSE}
%%
%We are now ready to give an analogue of Theorem~\ref{THM:M-CANONICAL} for \doms of the third type.
\abstractfont{%
We now give a different kind of embedding for \doms of the third type.%
}

Let $M$ be a \dom of the third type, $\vg := \GM$, and $\vh := \HM$.
If $\vg$ is densely ordered, define $\vk := \vg$.
Otherwise, let $1$ be the minimal positive element of $\vg$, and $\vk := \vg + 1 \cdot \Raz$, the subgroup of  $\vg \otimes \Raz$ (the divisible hull of $\vg$) generated by $\vg$ and $\nicefrac{1}{n}$, as $n$ varies in $\Nats$.
Since $\vk$ has no minimal positive element, $\vk$ is densely ordered.
Consider the \dom $\vkh$ of the Dedekind cuts of $\vk$; since $\vk$ is  densely ordered, $\vkh$ is a \dom of the third type.
We have seen that the group $\vk$ embeds into $\G{\vk}$, via the map sending $\gamma$ to $[\gp]$, hence $\vh$ also embeds into $\G{\vk}$.
Note also that the whole $\vk$ is contained in $\HL$, hence in particular $\vh$ is contained in $\HL$; call $\iotahat$ such embedding.
Therefore, we can define the \dom $\coll{\vkh}{\vh}$.
\begin{thm}
If $M$ is a strongly proper \dom of the third type, then, 
with the above definitions of $\vh$, $\vk$ and $\iotahat$, and $N := \coll{\vkh}{\vh}$, there exists a unique \domhom $\psihat: M \to N$ such that the following diagram commutes:
\begin{diagram}
M  ¤ \Edotar{\psihat} ¤ N ¤¤ %%
\Sepi\pi ¤ ¤ \sepI\pi ¤¤
\Mquot ¤ \Emono{\iotahat} ¤ \Lquot ¤¤
\end{diagram}
Moreover, $\psihat$ is injective.
\end{thm}
We recall that $N$ is $\Lquot \fprod[\vh] \fdom{2}$.
\begin{proof}
It is easy to see that $M$ and $\vk$ satisfy the condition (*) of Remark~\ref{REM:LAMBDA-COFINAL}.
Hence, the map $\Lambda_K : \notmzero \to \vkh$ sending $w$ to $\set{\lambda \in \vk: \exists y \in \Mgroup: \lambda \leq [y] \et y < w}^+$ is an injective \qdomhom.
Moreover, the only possible value for $\psihat(w)$ for $w \in \notmzero$ is $\Lambda_K(w)$.

It remains to define $\psihat(z)$ for $z \in \Mgroup$.

Since $z$ is in $\Mgroup$, then $[z]\in \vk$, hence $[z]^+$ and $[z]^-$ are elements of $\vkh$.
Let $B := \pi^{-1}(H) \subseteq \Mgroup$, $C := \Mgroup \setminus B$,
$B' := \pi^{-1}(H) \subseteq \zerogroup{\vkh}$, $C' := \zerogroup{\vkh} \setminus C$.

If $z = x \in C$, then $[x]^+$ and $[x]^-$ are identified in $N$ to the same element, which we will also call $[x]$, and we must define $\psihat(x) := [x] \in C' $.

If $z = y \in B$, then $y$ has multiplicity $2$, namely $[y] = \set{\FM(y), \FP(y)}$.
Moreover, $[y]^-$ and $[y]^+$ are distinct elements of $B'$.
If we want to preserve the signature, we must define $\psihat\Pa{\FM(y)} := [y]^-$, and $\psihat\Pa{\FP(y)} := [y]^+$.

It remains to prove that $\psihat$ is indeed a \domhom.

Let us prove that $\psihat(z + z') = \psihat(z) + \psihat(z')$.
Since the diagram commutes, $\pi\Pa{\psihat(z + z')} = \pi\Pa{\psihat(z) + \psihat(z')}$, hence the only case when we might not have equality is when $\psihat(z + z') \in B'$, \ie\ when $z + z' \in B$.
Note also that, by definition, $\psihat$ preserves the signature of every element of $\Mgroup$, and that $\sig y \neq 0$ for every $y \in B$.
It is enough to prove that $\psihat(z) + \psihat(z')$ and $z + z'$ have the same signature.
There are two possible cases: $z = x$ and $z' = x'$ are both in $C$, or $z = y$ and $z' = y'$ are both in $B$.
In the first case, $\sig x = \sig x' = 0$, hence, by Proposition~\ref{PROP:SIG-RULE}, $\sig{x + x'} = -1$. The same is true for $\psihat(x)$ and $\psihat(x')$, and we have the conclusion.
In the second case, by Proposition~\ref{PROP:SIG-RULE}, the signature of $y + y'$ depends only on the sign of $y$ and of $y'$, and the same is true for $\psihat(y)$ and $\psihat(y')$, and we can conclude.
The fact that $\psihat$ is injective is now trivial.
\end{proof}
Note that in the above proof we could not use, instead of $\vk$, an arbitrary densely ordered group containing $\vg$. More precisely,
every $w \in \notmzero$ determines the following partition of~$\vk$:
\begin{multline*}
\Lambda'(w) := 
\bigl(
{\set{\lambda \in \vk: \exists y \in \Mgroup: \lambda \leq [y] \et y < w}},\\
{\set{\lambda \in \vk: \exists y \in \Mgroup: \lambda \geq [y] \et y > w}}
\bigr).
\end{multline*}
The problem lies in the fact that $\Lambda'(w)$ is \emph{not} a cut in general (\cf\ Lemma~\ref{LEM:LAMBDA}).

For instance, let $\vg := \Zed \times \Real$, $M := \vgh$, $\vk := \Raz \times \Real$, $a := \set{(0,q): q \in \Real}^+  = \set{(1,q): q \in \Real}^-\in M$. 
Note that $\vg$ is already dense, and that $\GM = \vg$.
Then, 
\[\begin{aligned}
{\Lambda'(a)}^L &= \set{(0,q): q \in \Real}^+,\\
{\Lambda'(a)}^R &= \set{(1,q): q \in\Real}^-.
\end{aligned}\]
Hence, $\Lambda'(a)^L< (\unmezzo, y) < \Lambda'(a)^R$ for every $y \in \Real$.

\section{Axiomatisation of \doms}\label{SEC:AXIOMS}
\begin{lemma}
The axioms~\ref{AX:M1}, \ref{AX:M4}, \ref{AX:M5-RIGHT} and~\ref{AX:M5-LEFT} are independent. %, even in the case of finite structures.
That is, if we choose any one of them, we can find a \predom satisfying the other ones, but not the chosen axiom.
\end{lemma}
Example~\ref{EX:GROUP-DELTA} shows that axioms~\ref{AX:M1} (if we choose $\delta > 0$) and~\ref{AX:M4} (choosing $\delta < 0$ and different from $-1$) are independent.

For Axiom~\ref{AX:M5}, the examples are given below,
via the addition tables of some finite \predoms.

The \predom with $n$ element will be given as the set of the first $n$ natural element $\set{0,2,3,\dotsc, n-1}$, ordered in the usual way.
Note that if we want our structure to be a \predom, the only possible definition  of $-i$ is $-i := n \md i$ for every $i < n$ (where we denoted with $\md$ the minus on the integers).
Axiom~\ref{AX:PA} is equivalent to the fact that the sum increases as we go from the left to the right on the same row.
The commutativity of the sum is equivalent to the fact that the table is symmetric around the principal diagonal.

Moreover, axioms~\ref{AX:M1} and \ref{AX:M4} will be satisfied iff the neutral element is $\zO := \ceiling{n/2}$, namely $n/2$ if $n$ is even, $(n+1)/2$ if $n$ is odd, and hence $\zD = \floor{n/2}$.
Axiom~\ref{AX:M5-LEFT} is equivalent to $-x + x \leq \zD$, namely every element in the addition table in the anti-diagonal is less or equal to $\zD$.
Axiom~\ref{AX:M5-RIGHT} is equivalent to $x > -y \rightarrow y + x > \zD$, namely every element below the anti-diagonal is $> \zD$.
%\end{answer}
\begin{enumexamples}
\item A \predom with $3$ elements satisfying axioms~\ref{AX:M1}, \ref{AX:M4} and~\ref{AX:M5-RIGHT}, but not~\ref{AX:M5-LEFT} (the neutral element is in bold):
%\begin{center}
\[
\begin{array}{|c|ccc|}
%\toprule
\hline
\rule{0pt}{11pt}
+ & 0 & \textbf{1} & 2 \\[1.5pt]
\hline
\rule{0pt}{11pt}
0 & 0 & 0 & 2 \\
\textbf{1} & 0 & 1 & 2 \\
2 & 2 & 2 & 2\\[1.5pt]
\hline
%\bottomrule
\end{array}
\]
%\end{center}
\item Two \predoms with $4$ elements satisfying axioms~\ref{AX:M1}, \ref{AX:M4} and~\ref{AX:M5-LEFT}, but not~\ref{AX:M5-RIGHT}:
%\begin{center}
\[
\begin{array}{|c|cccc|}
%\toprule
\hline
\rule{0pt}{11pt}
+ & 0 & 1 & \textbf{2} & 3\\[1.5pt]
\hline
\rule{0pt}{11pt}
0 & 0 & 0 & 0 & 0 \\
1 & 0 & 1 & 1 & 1 \\
\textbf{2} & 0 & 1 & 2 & 3 \\
3 & 0 & 1 & 3 & 3\\[1.5pt]
\hline
%\bottomrule
\end{array}
\qquad
\begin{array}{|c|cccc|}
%\toprule
\hline
\rule{0pt}{11pt}
+ & 0 & 1 & \textbf{2} & 3 \\[1.5pt]
\hline
\rule{0pt}{11pt}
%\midrule
0 & 0 & 0 & 0 & 0 \\
1 & 0 & 0 & 1 & 1 \\
\textbf{2} & 0 & 1 & 2 & 3 \\
3 & 0 & 1 & 3 & 3\\[1.5pt]
\hline
%\bottomrule
\end{array}
\]
%\end{center}
\item For confrontation, here are the addition tables of the (trivial) \doms with $3$, $4$ and $5$ elements respectively:
%\begin{center}
\[
\begin{array}{|c|ccc|}
%\toprule
\hline
\rule{0pt}{11pt}
+ & 0 & \textbf{1} & 2 \\[1.5pt]
\hline
\rule{0pt}{11pt}
%\midrule
0 & 0 & 0 & 0 \\
\textbf{1} & 0 & 1 & 2 \\
2 & 0 & 2 & 2\\[1.5pt]
\hline
%\bottomrule
\end{array}
\qquad
\begin{array}{|c|cccc|}
\hline
\rule{0pt}{11pt}
%\toprule
+ & 0 & 1 & \textbf{2} & 3 \\[1.5pt]
\hline
\rule{0pt}{11pt}
%\midrule
0 & 0 & 0 & 0 & 0 \\
1 & 0 & 1 & 1 & 3 \\
\textbf{2} & 0 & 1 & 2 & 3 \\
3 & 0 & 3 & 3 & 3\\[1.5pt]
\hline
%\bottomrule
\end{array}
\qquad
\begin{array}{|c|ccccc|}
\hline
\rule{0pt}{11pt}
%\toprule
+ & 0 & 1 & \textbf{2} & 3 & 4 \\[1.5pt]
\hline
\rule{0pt}{11pt}
%\midrule
0 & 0 & 0 & 0 & 0 & 0\\
1 & 0 & 1 & 1 & 1 & 4 \\
\textbf{2} & 0 & 1 & 2 & 3 & 4\\
3 & 0 & 1 & 3 & 3 & 4\\
4 & 0 & 4 & 4 & 4 & 4\\[1.5pt]
\hline
%\bottomrule
\end{array}
\]
%\end{center}
\end{enumexamples}
The above tables were obtained using the Alloy program.%
\footnote{\url{http://alloy.mit.edu/}}
Thanks to Ivan Lanese for explaining to me how to use Alloy.
Note that the difficult part in such tables is checking whether the addition is associative: for instance, the following one is not.
\[\begin{array}{|c|cccc|}
\hline
\rule{0pt}{11pt}
%\toprule
+ & 0 & 1 & \textbf{2} & 3 \\[1.5pt]
\hline
\rule{0pt}{11pt}
%\midrule
0 & 0 & 0 & 0 & 0 \\
1 & 0 & 0 & 1 & 3 \\
\textbf{2} & 0 & 1 & 2 & 3 \\
3 & 0 & 3 & 3 & 3\\[1.5pt]
\hline
%\bottomrule
\end{array}\]
Finally, note that by Proposition~\ref{PROP:DOM}(\ref{EN:DOM-20}), $x - y$ can be defined in terms of the plus and order alone.
In particular, $- y = \zD - y = \max\set{z \in M: y + z \leq \zD}$.
Hence, the minus can be defined in terms of the plus, the order and $\zD$.

Therefore, if $\struct(M, \leq, \zO, +)$ is an Abelian ordered monoid, and $\zD \in M$, then there is at most one minus on $M$ such that
$\struct(M, \leq, \zO, +, -)$ is a \dom.
The necessary and sufficient conditions for the existence of such minus are the following:
\begin{itemize}
\item $\zD \leq \zO$;
\item the interval $(\zD, \zO)$ is empty;
\item for every $y \in M$, $- y := \max\set{z \in M: y + z \leq \zD}$ exists;
\item $-(-y) = y$.
\end{itemize}
%%
%%
%% \begin{question}
%% Let $M$ be a \predom, satisfying axioms~\ref{AX:M1} and~\ref{AX:M4}, plus Proposition~\ref{PROP:DOM}(\ref{EN:DOM-9}).
%% Is $M$ a \dom? And what happens if we also add the condition that $M$ is finite?
%% \end{question}
%% %%
%% \begin{answer}
%% A structure $M$ as above must also satisfy Axiom~\ref{AX:M5-LEFT}.
%% On the other hand, while we think that $M$ will not be a \dom, we do not have any counter-example.
%% \end{answer}
%%
Call~\ref{AX:MCP} the following axiom (which is equivalent to
Proposition~\ref{PROP:DOM}(\ref{EN:DOM-9})):
%We remind that Proposition~\ref{PROP:DOM}(\ref{EN:DOM-9}) says
\begin{enumerate}[label={\ref{AX:M5}'.}, ref={\ref{AX:M5}'}]
\item%[\ref{AX:M5}']
\label{AX:MCP}
$\displaystyle (x + y) - z \geq x + (y - z)$.%\tag*{\ref{AX:M5}'}.
\end{enumerate}
\begin{lemma}
Let $M$ be a \predom, satisfying axioms~\ref{AX:M1}, \ref{AX:M4} and \ref{AX:MCP}.
Then, $M$~is a \dom.
\end{lemma}
\begin{proof}
\begin{pseudoclaim}
$x \ar \zO \geq x$.
\end{pseudoclaim}
In fact, $x \ar \zO = -\Pa{(-x) + \zD }$.
Since $\zD \leq \zO$, the latter is greater or equal to \mbox{$-\Pa{(-x) + \zO} = x$}.
\begin{claim}
$x \ar z \geq x + z$.
\end{claim}
In fact, 
\[
x + z \leq (x \ar \zO) + z \leq  x \ar (\zO + z) = x \ar z,
\]
where the second inequality is obtained from Axiom~\ref{AX:MCP},
using~$y = \zO$.

We remind that, by Remark~\ref{REM:MC-EQ}, in a \predom Axiom~\ref{AX:M5-LEFT}
is equivalent to the following claim.
\begin{claim}
$\ig x \geq \zO$.
\end{claim}
In fact, $- \ig x = x \dl x \leq x - x = \ig x$.
The conclusion follows from Axiom~\ref{AX:M4}.

Therefore, it remains to prove the following claim.
\begin{claim}
If $x < y$, then $x - y < \zO$.
\end{claim}
In fact, $x < y \leq y + \ig x$.
By substituting $x := y$, $y := -x$, and $z := x$ in Axiom~\ref{AX:MCP}, we obtain that the latter is less or equal to $(y \dl x) \ar x$.
Hence,
\[
x \ar \zD = x < x \ar (y \dl x).
\]
Therefore, $y \dl x > \zD$, which is equivalent to the conclusion.
\end{proof}
Note that the \predom $M$ in Example~\ref{EX:GROUP-DELTA} satisfies Axiom~\ref{AX:MCP} for every choice of $\delta \in \vg$ (even in the stronger form $(x + y) - z = x + (y - z)$).
Moreover, for suitable choice of~$\delta$, $M$ will satisfy \ref{AX:M1} and not \ref{AX:M4}, or  \ref{AX:M4} and not \ref{AX:M1}.
Hence, the axioms \ref{AX:M1}, \ref{AX:M4}, and \ref{AX:MCP} are independent.
%%%

\section{Valuations}\label{SEC:VALUATION}
\begin{definizione}[Valued \doms]
A \intro{valued \dom} is a triple $\struct(M,v,C)$, where $M$ is a \dom, $C$ an ordered set with a minimum $-\infty$, and $v: M \to C$ (the \intro{valuation}) is a surjective map satisfying the following conditions: for every $x$, $y \in M$,
\begin{enumerate}[label={\textbf{V}\arabic*.}, ref={\textbf{V}\arabic*}]
\item\label{AX:VAL-ZERO} $v(\zO) = - \infty$;
\item\label{AX:VAL-MINUS} $v(-x) = v(x)$;
\item\label{AX:VAL-PLUS} $v(x + y) \leq \max \set{v(x), v(y)}$.
\end{enumerate}
For every $c \in C$, define
\[
\vring(c) := \set{x \in M: v(x) \leq c}.
\]
The valuation $v$ is \intro{convex} if it satisfies the following condition:
\begin{enumerate}[label={\textbf{V}4.}, ref={\textbf{V}4}]
\item\label{AX:VAL-CONVEX} $\abs x \leq \abs y$ implies $v(x) \leq v(y)$.
\end{enumerate}
The valuation $v$ is \intro{strong} if in the Axiom~\ref{AX:VAL-PLUS} equality holds.
\end{definizione}
For the rest of this section, $\struct(M,v,C)$ is a valued \dom.
\begin{remark}
For every $x, y \in M$,
\[v(x \ar y) \leq \max\set{v(x), v(y)},\]
and the same for $v(x - y)$ and $v(x \dl y)$.
Moreover, for every $c \in C$, $\vring(c)$ is a \subdom of~$M$.
\end{remark}
\begin{proof}\hspace{2ex}
%\begin{multline*}
$v(x \ar y) = v\Pa{-((-x) + (-y))}
= v\Pa{(-x) + (-y)} \leq$\\[0.7ex]
\mbox{}\hfill $\leq\max\set{v(-x), v(-y)} = \max \set{v(x), v(y)}.$
\hspace{2ex}
%\end{multline*}
\\[1ex]
The fact that $\vring(c)$  is a \subdom is now trivial.
\end{proof}
Note that the family $\Pa{\vring(c)}_{c\in C}$ is an increasing family of \subdoms of~$M$.
\begin{definizione}
Given another valuation $v': M \to C'$ on $M$, we say that $\struct(M, v', C')$ is a \intro{coarsening} of $v$, or that $v$ is a \intro{refinement} of $v'$, iff there exists an order-preserving map $\chi: C \to C'$ such that the following diagram commutes:
\begin{diagram}
M  ¤ \Eepi{v} ¤ C ¤¤ %%
¤ \Seepi{v'} ¤ \Sdotar{\chi} ¤¤
¤ ¤ C' ¤¤
\end{diagram}
\end{definizione}
\begin{remark}\label{REM:COARSE}
The function $\chi$ in the above definition, if it exists, is unique and surjective.
Moreover, the existence of $\chi$ is equivalent to:
\begin{itemize}
\item[(*)] for every $x$, $y \in M$, if $v(x) \leq v(y)$, then $v'(x) \leq v'(y)$.
\end{itemize}
\end{remark}
\begin{proof}
The only possible definition of $\chi$ is $\chi(v(x)) = v'(x)$.
It is evident that the above defined map $\chi$ is well-defined and order-preserving iff (*) is true.
Finally, $\chi$ is surjective, because $v$ and $v'$ are.
\end{proof}
\begin{enumexamples}
\item The map sending every element of $M$ to $-\infty$ is a valuation, the \intro{trivial} valuation.
\item If $M$ is not of the second type, the map $v: M \to \set{-\infty, 1}$ sending $\pm \zO$ to $- \infty$, and everything else to $1$ is a valuation.
\item On an ordered group~$\vg$, a valuation in the group-theoretic sense is also a  valuation in our sense.
On the other hand, a \dom-valuation $v$ is a group-valuation iff $0$ is the only element such that $v(x) = -\infty$.
\item The triple $\struct(M, \ig{\ }, \orders)$%
\footnote{Where $\ig{\ }(x) = x - x$}
is a strong valuation on~$M$, called the \intro{\width{}} valuation.
\item For every $x$, $y\in M$, define $x \archleq y$ if there exists $n \in \Nat$ such that
\[
\abs{x} \leq \underbrace{\abs{y} \ar \dotsb \ar \abs{y}}_{n \text{ times}}.
\]
Define $x \archeq y$ if $x \archleq y$ and $y\archleq x$.
It is obvious that $\archleq$ is a total pre-order, hence $\archeq$ is an equivalence relation on~$M$, and $\archleq$ induces an ordering on~$\Marcheq$.
Moreover, the equivalence class of $\zO$ is the minimum of~$\Marcheq$.
Finally, if we call $\archv$ the quotient map from $M$ to~$\Marcheq$,
then $\struct(M, \archv, \Marcheq)$ is a valued \dom.
The map $\archv$ is the \intro{natural} valuation on~$M$.
\end{enumexamples}
Note that the concept of strong valuation is trivial in the case when $M$ is an ordered group (namely, the only strong valuation is the trivial one).
\begin{lemma}\label{LEM:VAL-NATURAL}
A valuation $v$ is convex iff it is a a coarsening of the natural valuation, iff for every $c\in C$ the \subdom $\vring(c)$ is convex.
\end{lemma}
\begin{proof}
Immediate from Remark~\ref{REM:COARSE}.
\end{proof}
\begin{lemma}\label{LEM:VAL-ORDER}
A valuation is strong iff it is a coarsening of the \width valuation.
Moreover, in that case $v(x \ar y) = \max\set{v(x),v(y)}$, and the same for $v(x - y)$ and $v(x \dl y)$.
\end{lemma}
\begin{proof}
If $v$ is strong, then $v(\ig x) = v(x)$.
The conclusion now follows from Remark~\ref{REM:COARSE}.
\end{proof}
\begin{lemma}
If $\struct(M,v,C)$ is either strong or convex, then, for every $x, y \in M$,
$v(x \ar y) = v(x + y)$, and $v(x \dl y) = v(x - y)$.
Moreover, if $v(x) < v(y)$, then $v(x + y) = v(y)$.
\end{lemma}
\begin{proof}
The conclusion is true for the natural and the \width valuations.
A valuation $v$ satisfying the hypothesis is a coarsening of the natural or the trivial valuations, by Lemmata~\ref{LEM:VAL-NATURAL} and~\ref{LEM:VAL-ORDER}.
The conclusion now follows from Remark~\ref{REM:COARSE}.
\end{proof}
\begin{deflemma}
For every $x \in M$, define
\[
w(x) := \set{ \ig y \in \orders: y + \ig x = x \vel y - \ig x = x}^-,
\]
where the lower edge is taken in $\Dedekind{\orders}$.
Then, $w$~is a valuation on~$M$.
\end{deflemma}
Note that, by definition, $w(x) \leq \ig x$, that $w(\ig x) = \zO^-$,
and that $w(x) = \zO^-$ iff there exists $y \in \Mgroup$ such that
$x = y \pm \ig x$.
% (where $0^-$ is the minimum value in the image of $w$).
\begin{proof}
Let us prove that $w(-x) = w(x)$.
Let $y \in M$ such that $y \pm \ig x = x$.
Hence, \mbox{$(-y) \mp \ig x = -x$}, thus $w(-x) \leq w(x)$, and we are done.

Let $x$, $x'\in M$ such that $\ig x' \leq \ig x$.
Let us prove that $w(x + x') \leq w(x)$.
Choose $y$ and $y' \in M$ such that $y \pm \ig x = x$, $y' \pm \ig x' = x'$, and $\ig y \leq \ig x$, and $\ig y' \leq \ig x'$.
\begin{claim}
There exists $z \in M$ such that $\ig z \leq \ig{y + y'}$, and $x + x' = z \pm \ig x$.
\end{claim}
If $\ig y = \ig x$, or $\ig y' = \ig x$, then choose $z = x + x'$, and we are done.
Otherwise, let $z = y + y'$: in this case, $\ig z < \ig x$.
Then, $x + x' \leq z + \ig x$.
Moreover, $z - \ig x  = y + y' \dl \ig x = (y - \ig x) + (y' - \ig x) \leq x + x'$.
Let $t := x + x'$: note that $\ig t = \ig x$.
Suppose, for contradiction, that, $z - \ig t < t < z + \ig t$.
Thus, $(z + \ig t) \dl (z - \ig t) > t - t $, hence $\ig t \dl \ig z > \ig t$, absurd.

The conclusion now follows from the claim, since for every \width $o$ such that $o > \max\set{w(x), w(x')}$, there exist $\ig y$ and $\ig {y'}$ such that $o \geq \ig y > w(x)$ and $o \geq \ig{y'} > w(x')$.
By the claim, we have that $\ig{y + y'} > w(x + x')$, therefore $o > w(x + x')$, and we are done.
\end{proof}
\begin{enumexample}
Let $\vg$ be an ordered group.
It might happen that there exists $\Lambda \in \vgh$ such that $w (\Lambda) > \zO^-$.
For instance,
let $\vg$ be the group $\Raz \times \Real$,
with the lexicographic ordering.
Let $\Lambda := \set{(q, t): \sqrt 2 >q \in \Raz, t \in \Real}^+$.
Then, $w(\Lambda) > \zO^-$.
\end{enumexample}
%%%

%%
\section{Conclusion}
We have shown the fact that the theory of \doms (plus the axiom $- \zO < \zO$) is the universal part of the theory of Dedekind cuts of ordered (Abelian) groups.
This means that if a universal sentence for cuts is true, then it can be proven using the axioms for \doms alone and $\zO < \zO$.
%A typical example is Corollary~\ref{COR:D-K-M-J}, which was the sentence that the first author originally wanted to prove when he undertook the study of \doms.
Moreover, a \dom is nothing else than a substructure of $\vgh$ or $\vgt$ for some ordered group $\vg$, and the axioms of \doms characterise the class of such substructures.
Some natural questions we left open are the following:
\begin{itemize}
\item Is there a ``nice'' (e.g.\ recursive) axiomatisation for the (first-order) theory of cuts of ordered groups?
\item What is the model-completion of the theory of \doms (if it exists)?
\end{itemize}

Using results of Baur~\cite{BAUR:1976}, it is not difficult to see that the theory of \doms is undecidable.
Some related questions are:
\begin{itemize}
\item Does it exists an algorithm to decide which \emph{universal} formulae follows from the theory of \doms?
\item Is the theory of of cuts of ordered groups decidable?
\end{itemize}
\paragraph{Bibliographical notes.}
Dedekind cuts of an ordered Abelian group $\vg$  have been extensively studied, especially for the purpose of building the Dedekind completion of~$\vg$ \cite{COH-GOFF:1949,BIRKHOFF:1967,FUCHS:1963}.
A more detailed study of the arithmetic properties of the set of Dedekind cuts of $\vg$ has also been undertaken by several authors, especially for the case when $\vg$ is the additive group of an ordered field~\cite{GONSHOR:1985,WEHRUNG:1996,PESTOV:2001,TRESSL:2005,KUH?22}.
The concepts of \width of an element and signature were already introduced by Gonshor~\cite{GONSHOR:1985}.
For the reader's convenience, we include a ``translation'' between Gonshor's notation~\cite{GONSHOR:1985} and ours:
\begin{itemize}
\item $\mathrm{ab}(x)$, the absorption number of $x$, is $\ig x$, the \width of $x$;
\item $x$ is a positive idempotent iff $x$ is a \width of~$M$;
\item if $0 < k \in \Orders{M}$, then:\vspace{-1ex}
\[\begin{aligned}
x R y \mod k &&\text{iff}&&& x + k = y + k;\\
x S y \mod k &&\text{iff}&&& [x + k]_k = [y + k]_k;\\
x T y \mod k &&\text{iff}&&& \abs{x - y} < k;\\[-1ex]
\end{aligned}\]
\item $x$ has type 1 iff $\sig x = 1$, while $x$ has type 1A iff $\sig x = -1$.
\end{itemize}

\bibliography{valuation_theory,my_bibliography}

\def\cprime{$'$}
\begin{thebibliography}{10}

\bibitem{BAUR:1976}
W.~Baur.
\newblock Undecidability of the theory of abelian groups with a subgroup.
\newblock {\em Proc. Amer. Math. Soc.}, 55(1):125--128, 1976.

\bibitem{BIRKHOFF:1967}
G.~Birkhoff.
\newblock {\em Lattice theory}.
\newblock Third edition. American Mathematical Society Colloquium Publications,
  Vol. XXV. American Mathematical Society, Providence, R.I., 1967.

\bibitem{COH-GOFF:1949}
L.~W. Cohen and C.~Goffman.
\newblock The topology of ordered {A}belian groups.
\newblock {\em Trans. Amer. Math. Soc.}, 67:310--319, 1949.

\bibitem{FUCHS:1963}
L.~Fuchs.
\newblock {\em Partially ordered algebraic systems}.
\newblock Pergamon Press, Oxford, 1963.

\bibitem{FUCHS:1970-73}
L.~Fuchs.
\newblock {\em Infinite abelian groups. {V}ol. {I} \& {II}}.
\newblock Pure and Applied Mathematics, Vol. 36 \& 36-II. Academic Press, New
  York, 1970 \& 1973.

\bibitem{GONSHOR:1985}
H.~Gonshor.
\newblock Remarks on the {D}edekind completion of a nonstandard model of the
  reals.
\newblock {\em Pacific J. Math.}, 118(1):117--132, 1985.

\bibitem{GULDENBERG:2004}
T.~G{\"u}ldenberg.
\newblock Elementare invarianten von dedekindschnitten angeordneter k{\"o}rper.
\newblock Master's thesis, Universit{\"a}t Regensburg, 2004.

\bibitem{KUH?22}
F.-V. Kuhlmann.
\newblock Invariance group and invariance valuation ring of a cut.
\newblock Unpublished, 6 2004.

\bibitem{MACLANE:1967}
S.~Mac~Lane.
\newblock {\em Homology}.
\newblock Springer-Verlag, Berlin, first edition, 1967.
\newblock Die Grundlehren der mathematischen Wissenschaften, Band 114.

\bibitem{MACLANE:1998}
S.~Mac~Lane.
\newblock {\em Categories for the working mathematician}, volume~5 of {\em
  Graduate Texts in Mathematics}.
\newblock Springer-Verlag, New York, second edition, 1998.

\bibitem{PESTOV:2001}
G.~G. Pestov.
\newblock On the theory of cuts in ordered fields.
\newblock {\em Sibirsk. Mat. Zh.}, 42(6):1350--1360, iii, 2001.
\newblock [English translation: Siberian Math. J. 42 (2001), no. 6,
  1123--1131].

\bibitem{SCOTT:1969}
D.~Scott.
\newblock On completing ordered fields.
\newblock In {\em Applications of Model Theory to Algebra, Analysis, and
  Probability (Internat. Sympos., Pasadena, Calif., 1967)}, pages 274--278.
  Holt, Rinehart and Winston, New York, 1969.

\bibitem{TRESSL:2005}
M.~Tressl.
\newblock Model completeness of o-minimal structures expanded by {D}edekind
  cuts.
\newblock {\em J. Symbolic Logic}, 70(1):29--60, 2005.

\bibitem{WEHRUNG:1996}
F.~Wehrung.
\newblock Monoids of intervals of ordered abelian groups.
\newblock {\em J. Algebra}, 182(1):287--328, 1996.

\end{thebibliography}
\bibliographystyle{abbrv}

\end{document}